\documentclass[reqno]{amsart}

\usepackage[utf8]{inputenc}
\usepackage[T1]{fontenc}
\usepackage[english]{babel}
\usepackage[backend=biber, style=numeric, sorting=nyt, maxbibnames=10, doi=false,isbn=false,url=false,
eprint=false, giveninits=true]{biblatex}
\usepackage{csquotes}

\usepackage{latexsym}
\usepackage{amssymb}
\usepackage{dsfont}
\usepackage{mathtools}
\usepackage[hidelinks]{hyperref}
\usepackage{comment}
\usepackage{enumitem} 
\usepackage{float}
\usepackage{multicol}
\usepackage{stmaryrd} 
\usepackage{tikz}
\usepackage{tikz-3dplot}
\usetikzlibrary{patterns,calc,decorations.markings,matrix,backgrounds,fit,positioning,3d, arrows.meta}

\allowdisplaybreaks[4]

\theoremstyle{plain}
\newtheorem{theorem}{Theorem}[section]
\newtheorem*{proposition*}{Proposition}
\newtheorem{proposition}[theorem]{Proposition}
\newtheorem{lemma}[theorem]{Lemma}
\newtheorem{corollary}[theorem]{Corollary}
\newtheorem{definition}[theorem]{Definition}
\theoremstyle{definition}

\newtheorem{remark}[theorem]{Remark}

\numberwithin{equation}{section}

\newcommand{\R}{\mathbb{R}} 
\newcommand{\N}{\mathbb{N}} 
\newcommand{\ran}{\operatorname{ran}}
\renewcommand{\ker}{\operatorname{ker}}

\newcommand{\bH}{\mathbb{H}}

\newcommand{\Rsym}{\R^{3\times3}_{\mathrm{sym}}}

\newcommand{\cC}{\mathcal{C}}

\newcommand{\cG}{\mathcal{G}}
\newcommand{\cH}{\mathcal{H}}

\newcommand{\cN}{\mathcal{N}}

\newcommand{\frakf}{\mathfrak{f}}
\newcommand{\frakg}{\mathfrak{g}}

\renewcommand{\div}{\operatorname{div}}
\DeclareMathOperator{\curl}{curl}
\newcommand{\tr}{\operatorname{tr}}
\newcommand{\dist}{\operatorname{dist}}

\newcommand*\dd{\mathop{}\!\mathrm{d}}

\newcommand{\pt}{\partial_t}

\newcommand{\1}{\mathds{1}} 
\newcommand{\e}{\mathrm{e}}

\let\latexepsilon\epsilon
\let\latexvarepsilon\varepsilon
\renewcommand{\epsilon}{\latexvarepsilon}
\renewcommand{\varepsilon}{\latexepsilon}

\newcommand{\supp}{\operatorname{supp}}

\let\oldexists\exists \let\exists\relax \DeclareMathOperator{\exists}{\oldexists} 
\let\oldforall\forall \let\forall\relax \DeclareMathOperator{\forall}{\oldforall}

\newcommand{\inv}[1]{{#1}^{-1}}

\newcommand{\ol}{\overline}

\let\latexchi\chi
\makeatletter
\renewcommand\chi{\@ifnextchar_\sub@chi{{\latexchi}}}
\newcommand{\sub@chi}[2]{
	{\@ifnextchar^{\subsup@chi{#2}}{{\latexchi}^{}_{#2}}}
}
\newcommand{\subsup@chi}[3]{
	{{\latexchi}_{#1}^{#3}}
}

\makeatletter
\def\moverlay{\mathpalette\mov@rlay}
\def\mov@rlay#1#2{\leavevmode\vtop{%
		\baselineskip\z@skip \lineskiplimit-\maxdimen
		\ialign{\hfil$\m@th#1##$\hfil\cr#2\crcr}}}
\newcommand{\charfusion}[3][\mathord]{
	#1{\ifx#1\mathop\vphantom{#2}\fi
		\mathpalette\mov@rlay{#2\cr#3}
	}
	\ifx#1\mathop\expandafter\displaylimits\fi}
\makeatother
\newcommand{\cupdot}{\charfusion[\mathbin]{\cup}{\cdot}}

\newcommand{\abs}[1]{ \left\lvert#1\right\rvert} 
\newcommand{\norm}[1]{\left\lVert#1\right\rVert} 
\newcommand{\smallnorm}[1]{\lVert#1\rVert} 
\newcommand{\bignorm}[1]{\big\lVert#1\big\rVert} 
\newcommand{\vertiii}[1]{{\left\vert\kern-0.25ex\left\vert\kern-0.25ex\left\vert #1
		\right\vert\kern-0.25ex\right\vert\kern-0.25ex\right\vert}}

\AtBeginDocument{
}

\newcommand\restr[2]{{
		\left.\kern-\nulldelimiterspace 
		#1 
		\vphantom{\big|} 
		\right|_{#2} 
}}

\newcommand{\seq}{\hphantom{={}}}

\newcommand{\epsilond}{\epsilon_{\mathrm{d}}}

\newcommand{\epsilonk}{\hat{\epsilon}_k}
\newcommand{\tildeepsilond}{\tilde{\epsilon}_{\mathrm{d}}}
\newcommand{\mud}{\mu_{\mathrm{d}}}

\newcommand{\muk}{\hat \mu_k}
\newcommand{\tildemud}{\tilde{\mu}_{\mathrm{d}}}

\newcommand{\Tmax}{T_{\mathrm{max}}}
\newcommand{\Tmin}{T_{\mathrm{min}}}

\newcommand{\epsilondnl}{\epsilon_{\mathrm{d,nl}}}

\newcommand{\mudnl}{\mu_{\mathrm{d,nl}}}

\hypersetup{final} 

\renewbibmacro{in:}{%
 \ifentrytype{article}{}{\printtext{\bibstring{in}\intitlepunct}}}
\DeclareFieldFormat[unpublished]{title}{\emph{#1}}
\DeclareFieldFormat[article]{pages}{#1}
\AtEveryBibitem{
	\ifentrytype{article}{}{\clearfield{pages}}
	\clearfield{edition}
	\clearfield{series}
	\ifentrytype{unpublished}{}{\clearfield{note}}
	\ifentrytype{book}{\clearfield{volume}}{}
	\ifentrytype{unpublished}{}{\clearfield{url}}
}
\begin{document}

\title[Decay of the quasilinear Maxwell system caused by conductivity]{Exponential decay of the quasilinear Maxwell system
	caused by conductivity on a collar}
\author{Richard Nutt and Roland Schnaubelt}
\address{RN \& RS, Karlsruhe Institute of Technology, Department of Mathematics,
	Englerstraße~2, 76131 Karlsruhe, Germany}
\email{richard.nutt@kit.edu, \quad schnaubelt@kit.edu }




\begin{abstract}
	We study the anisotropic, quasilinear Maxwell system on a bounded domain $\Omega$ with perfectly conducting boundary conditions.
	It is damped via a conductivity $\sigma$ which is strictly positive on its support containing a collar of the boundary.
	We prove that solutions exist globally and decay exponentially to $0$ if there are no magnetic charges on $\Omega$, no electric charges
	off the support of $\sigma$, and the initial fields are small. Our approach relies on a splitting of the solution via adapted
	Helmholtz decompositions and an observability-type estimate for a related system without charges, shown by means of Morawetz multipliers.
	Moreover, one has to control the space derivatives of the solution by its time derivatives, exploiting the structure of the system.
	The discontinuity of $\sigma$ introduces
	an interface which poses significant challenges in all these steps, besides those coming from the quasilinear terms.
\end{abstract}

\maketitle

\section{Introduction}

Electromagnetic theory is based on the Maxwell system. It contains constitutive relations which describe the response of the material
to the fields. The problem becomes a quasilinear hyperbolic system if these relations are nonlinear and instantaneous, as assumed here.
Moreover, a non-zero conductivity $\sigma\ge0$ causes energy loss and may lead to decay of solutions. In the quasilinear case
this topic was first studied for
strictly positive $\sigma$ in \cite{lasiecka2019}, see \cite{nutt2024} and \cite{pokojovy2020} for strictly positive
boundary damping. In this setting for small initial fields one obtains a global solution decaying to 0 exponentially.
In our recent contribution \cite{nutt2026}, we showed exponential decay in the
linear case if the support $\sigma$ contains a collar of the boundary and is strictly positive on it. In this case the stability analysis
becomes much more involved than in \cite{lasiecka2019}. On the other hand, the quasilinear setting poses significant new difficulties
compared to \cite{nutt2026}. These are addressed in this paper providing exponential decay for small data in the quasilinear case
with damping near the boundary. In addition, we now have to treat an interface problem also admitting composite materials.

\smallskip

We study the quasilinear Maxwell system
\begin{alignat}{-1}
	\begin{split}\label{eq:maxwell}
		\pt (\epsilon(x, E(t,x)) E(t,x)) & = \curl H(t,x) - \sigma(x)E(t,x), \\
		\pt (\mu(x, H(t,x)) H(t,x))      & = - \curl E(t,x),
	\end{split}
	\qquad                    &                                                    & t\geq 0, x \in \Omega_\pm,                                                       \\
	\begin{split}\label{eq:maxwell2}
		\div(\mu(x)H(t,x)) & = 0,
	\end{split}
	\qquad                    &                                                    & t\geq 0, x \in \Omega_\pm,                                                       \\
	E(t,x) \times \nu(x) = 0, & \quad \nu(x) \cdot \mu(x,H(t,x))H(t,x) = 0, \qquad &                            & t\geq 0, x \in \partial \Omega, \label{eq:maxwell3} \\
	E(0,x) = E_0(x),          & \quad H(0,x)= H_0(x), \qquad                       &                            & x \in \Omega_\pm, \label{eq:maxwell4}
\end{alignat}
for the electric and magnetic fields $E(t,x),H(t,x)\in\R^3$
on a smooth bounded domain $\Omega$, which is separated into subdomains $\Omega_+=\omega$ and $\Omega_-=\upsilon$ by a smooth interface
$\Gamma$ with positive distance to the boundary $\partial\Omega$. We allow for a multiply connected $\Omega$ and disconnected boundary
$\partial\Omega$, in contrast to \cite{lasiecka2019} and the related literature, except for \cite{nutt2026} and parts of \cite{phung2000}
dealing with linear problems.

The coefficients are required to be $C^3$ on each subdomain and to take values in symmetric
real $3\times 3$-matrices, where $\epsilon(x,\xi),\mu(x,\xi)\ge \eta>0$ and $\sigma(x)\ge0$. We thus admit anisotropic
and composite materials whose
parameters jump at $\Gamma$. Moreover, we require that $\omega$ contains a collar of $\partial\Omega$ and that
$\sigma$ is strictly positive on $\omega$ and vanishes on $\upsilon$. This behavior is crucial for our analysis. Strict positivity of
$\sigma$ on its support is often assumed in papers on stability of the linear Maxwell system, see \cite{eller2019},
\cite{nicaise2025}, \cite{nutt2026}, \cite{phung2000}, or \cite{phung2013}. Moreover, we suppose that there are no magnetic
charges (as given by physics) and no electric charges on $\upsilon$, which excludes stationary solutions prohibiting
exponential decay. The permittivity $\epsilon$ and permeability $\mu$ have to satisfy the (quite explicit) non-trapping conditions
\eqref{eq:non-trapping} and \eqref{eq:non-trapping-interface} on and near $\upsilon$ and at $\Gamma$, respectively, which heuristically
prevent backreflections and ensure that the waves reach the damping region. They are quite common in studies of boundary damping and
observability of the Maxwell system with nonconstant coefficients already in the linear case, see \cite{eller2007}, \cite{nicaise2005},
\cite{nutt2024}, \cite{nutt2026}, \cite{pokojovy2020}, and for piecewise constant coefficients in \cite{kapitonov1994}.
Our precise assumptions are stated in Section~\ref{sec:ass} and collected in (H). We show exponential decay for small data in
Theorem~\ref{thm:exp-decay}.

\vspace*{0.1cm}

Already the local wellposedness theory for the quasilinear Maxwell system is quite challenging, see \cite{schnaubelt2022}
for \eqref{eq:maxwell}--\eqref{eq:maxwell4} and also \cite{spitz2017}, \cite{spitz2018}, \cite{spitz2022} for the case without interface.
It is known that in quasilinear Maxwell systems for non-small data blow-up may occur, see \cite{dancona2018} and the references therein.
The facts needed below are recalled in Section~\ref{sec:ass}. To deal with the quasilinear terms one has to work with solutions
in the space $\bigcap_{j=0}^3 C^j((\Tmin,\Tmax), H_\Gamma^{3-j}(\Omega_\pm))$ that belong to the Sobolev space $H^3$
only piecewise because of the interface. The presence of the interface causes severe challenges compared to the previous papers
\cite{lasiecka2019} and \cite{nutt2026} (co-)authored by us, which lead to significant changes in all core steps. Besides the linear
controllability result \cite{kapitonov1994}, we are not aware of studies on the longtime behavior of the Maxwell system with interfaces.
The higher regularity level and the quasilinear setting cause additional difficulties compared to our linear paper \cite{nutt2026},
e.g., the need of intricate approximation arguments based on \cite{nutt2026-2}.

Compared to the damped wave equation, see e.g.\ \cite{bardos1988}, \cite{bardos1992}, \cite{lasiecka2006}, \cite{lasiecka1992-2},
\cite{rauch1974}, and the references therein,
the theory of decay properties of the Maxwell system is far less developed. In this context, we note that \eqref{eq:maxwell}
is a highly coupled system, in particular in the presence of anisotropic coefficients which can describe nontrivial crystal structures.
Even for scalar $\epsilon$ and $\mu$ in \eqref{eq:maxwell}, we have to pass to matrix coefficients in the nonlinear analysis,
see \eqref{eq:def-epsilon_nl}. Moreover, the system has quite complicated boundary and interface conditions which are much harder to treat in
the anisotropic case. Finally, the huge kernel of the curl operator is an obstacle to regularity and stability theory. In principle
Gau{\ss}' charge laws \eqref{eq:div-evolution} may help here, but in zeroth order they contain a time integral spoiling the longtime
behavior.

Still, the results and methods for the wave equation are the starting point for the treatment of the Maxwell system.
For constant coefficients $\epsilon, \mu \in \mathbb{R}_+$ Theorem~5.5 in \cite{phung2000} establishes exponentially decay to $0$
provided that the conductivity $\sigma$ is strictly positive on its support and satisfies the geometric control condition
derived from the wave equation framework, see \cite{bardos1988} or \cite{bardos1992}. In this case, the Maxwell system is still rather
close to the wave equation, and at the moment it is not clear how to extend these microlocal methods to problems like
\eqref{eq:maxwell}--\eqref{eq:maxwell4}.
We instead use the more robust approach of Morawetz multipliers whose variants have also been used successfully for wave systems, see
e.g.\ \cite{lasiecka2006} or \cite{lasiecka1992-2}.

\vspace*{0.1cm}

We use ideas from \cite{nicaise2005} dealing with observability for isotropic Maxwell systems based on Morawetz multipliers. However, these
arguments only work in the fully charge-free case, which prohibits a direct application to stability. To overcome this problem, we employ
Helmholtz decompositions as in \cite{phung2000} for constant coefficients. In \cite{nutt2026} we have extended this approach to linear
anisotropic problems on domains with nontrivial topology, which required significant changes compared to \cite{nicaise2005} and
\cite{phung2000}. In particular, we have to deal with obstructions to the invertibility of the curl caused by the topology,
e.g., by passing to somewhat more complicated state spaces that exclude certain stationary solutions,
see Lemma~\ref{lem:curl-with-topology}. In the anisotropic case with strictly positive $\sigma$, exponential decay
has been shown only recently in \cite{eller2019} and \cite{lasiecka2019}, see \cite{eller2007} also for boundary observability.
Weaker decay properties of Maxwell have been treated in \cite{nicaise2025}, \cite{phung2000}, \cite{phung2013}, and \cite{racke1992}.

In the quasilinear case we need to control higher derivatives of solutions.
As in \cite{lasiecka2019}, cf.\ also the references therein, we first treat the time derivatives of the fields $\pt^k (E,H)$ with $k\le 3$
which satisfy the Maxwell system with small time-depending error terms and the same boundary and interface conditions,
see e.g.\ \eqref{pde:pt-maxwell}--\eqref{pde:pt-maxwell_BC}. These error terms interfere with all our calculations and call for robust methods.
A core part of our analysis are the observability-type estimates in Corollary~\ref{cor:obs-est}. They are derived in lengthy and demanding
calculations by means of Morawetz multipliers using integration by parts, vector analysis, and the non-trapping conditions.
In contrast to \cite{lasiecka2019} the damping only acts on $\omega$, and we have to deal with an interface which causes significant trouble
and leads to new terms on the right-hand side of the estimate (which are absorbed
later). Another fundamental difference to \cite{lasiecka2019} is the need to remove charges by a Helmholtz-type correction terms
$(\nabla p^k,\nabla q^k)$ for $\pt^k (E,H)$ which solve appropriate elliptic interface problems. It is crucial to bound $(\nabla p^k,\nabla q^k)$
and their time derivatives by dissipation terms (plus small errors), see Lemma~\ref{lem:reg_of_pkqk}, but this only works for $k\ge1$.
Fortunately, we can control the energy of $(E,H)$ by that of $\pt (E,H)$ in Lemma~\ref{lem:cE-estimate} taken from \cite{nutt2026}.
The observability analysis is thus performed for the fields $(V^k,W^k)\coloneq \pt^k (E,H)- (\nabla p^k,\nabla q^k)$, where we
start from related computations for a second-order formulation in \cite{nutt2026}. But now the presence of the interface and the
perturbations coming from the nonlinear terms require various new ideas and make the analysis much harder.

Theorem~\ref{thm:obs-est} then controls the energy of the solutions $(E,H)$ and their time derivatives uniformly and $L^1$ in time
by norms of initial data, dissipation terms, and small errors that contain space derivatives. Here the smallness comes from the size of
the data or is enforced by our estimates. As in \cite{lasiecka2019} and \cite{pokojovy2020}, it remains to control space by time derivatives
of $(E,H)$ uniformly and $L^1$ in time with constants independent of $t$, see Theorem~\ref{thm:regularity-analysis}.
Having achieved this, global existence and exponential decay
for small data can be shown as in \cite{lasiecka2019} in a rather standard way.

For the remaining regularity analysis we split the solution
into parts with support in $\upsilon$, in $\omega$, near $\partial \Omega$, and near $\Gamma$. The first three can be treated
similarly to \cite{pokojovy2020}, respectively \cite{lasiecka2019}. The last one requires substantial
additional efforts. Since tangential derivatives of the solutions satisfy
the same interface conditions up to lower-order terms, they can be controlled on $\omega$ by means of energy estimates.
Normal derivatives are then reduced to tangential ones exploiting the structure of the Maxwell system, where we adapt ideas from
\cite{lasiecka2019}. These estimates require the strict positivity of $\sigma$. On $\upsilon$ the conductivity $\sigma$ vanishes
which leads to the absence of charges, up to lower-order error terms. This facilitates the use of div-curl estimates, see
Theorem~\ref{thm:div-curl-perfect}, but they require regularity of the traces on $\partial\upsilon=\Gamma$. Here we can
invoke the estimates for the fields on $\omega$ thanks to the interface conditions of the Maxwell system.
In all these inequalities we iteratively trade space into time derivatives, which
amounts to an intricate iteration scheme also relating normal to tangential derivatives, the different spatial areas, and the electric
to the magnetic fields (which are easier to handle via Lemma~\ref{lem:H-est}).


In the next three sections and the appendix we discuss background material and preliminary results. Section~\ref{sec:helm}
deals with our Helmholtz decomposition, Section~\ref{sec:obs} establishes the observability-type estimate, and
Section~\ref{sec:regularity-analysis} is devoted to the regularity analysis. The main decay theorem is stated and proven in the
last section.

\section{Notation}
\label{sec:notation}

We list notation used throughout this paper. The set $\Omega \subset \R^3$ is a bounded
domain (i.e., non-empty, open, and connected) with at least $C^5$-regularity. It is
\emph{multiply connected} if there are disjoint $C^k$-surfaces $\Sigma_1,\dots, \Sigma_L$ such that
$\Omega\setminus \bigcup_{l=1}^{L} \Sigma_l$ is simply connected, see \cite{ciarlet18}, \cite{cessenat1996} or \cite{dautray1990}.
The connected components of the boundary $\partial \Omega$ are denoted by $(\partial\Omega)_0 ,\dots, (\partial\Omega)_{K}$. Let
\[
	\cN_{\!a} \coloneqq \{ x \in \overline\Omega \mid \mathrm{dist}(x, \partial\Omega)<a\}
\]
be the collar of width $a>0$ around $\partial\Omega$ within $\Omega$.

The interior part of the boundary of $\supp \sigma$ is considered as an interface $\Gamma \subset \Omega$, and it
is assumed to have
same regularity as $\partial\Omega$ and a positive distance to it, i.e., $\dist(\Gamma, \partial\Omega)>0$.
So $\Gamma$ splits $\Omega$ into two subdomains. The one bordering $\partial\Omega$ is called
$\omega \subseteq \Omega$ and the other $\upsilon \subseteq \Omega$. For the union of these we write
$\Omega_\pm \coloneqq \omega \cupdot \upsilon$ and thus
\[
	\Omega = \Omega_\pm \cupdot \Gamma = \upsilon \cupdot \Gamma \cupdot \omega.
\]
The collar within $\Omega$ around $\Gamma$ of width $a>0$ is denoted by
\[
	\cN_a^\Gamma \coloneqq \{ x \in \overline\Omega \mid \mathrm{dist}(x, \Gamma)<a\}.
\]
We write $\nu^\Gamma$ for the unit normal at the interface pointing from $\upsilon$ to $\omega$.
The surface measure on $\partial\Omega$ or $\Gamma$ is designated by $\varsigma$.

At first, for $u \in C(\overline\Omega)$ we define the vector-valued traces

\vspace{-0.1cm}
\noindent
\begin{minipage}[t]{.48\textwidth}
	\vspace{0.2cm}
	\begin{alignat*}{-1}
		\tr_n u    & = \tr_{n, \Xi} u   &  & \coloneqq (u \cdot \nu^\Xi) \nu^\Xi,         \\
		\tr_t u    & = \tr_{t, \Xi} u   &  & \coloneqq u \times \nu^\Xi ,                 \\
		\tr_\tau u & = \tr_{\tau,\Xi} u &  & \coloneqq \nu^\Xi \times (u \times \nu^\Xi),
	\end{alignat*}
\end{minipage}
\begin{minipage}[t]{.51\textwidth}
	\vspace{-0.3cm}
	\resizebox{\linewidth}{!}{%
		\tdplotsetmaincoords{73}{90}

		\begin{tikzpicture}[scale=2,tdplot_main_coords]

			\coordinate (O) at (0,0,0);

			\coordinate (N) at (0,0,0.85);       
			\coordinate (U) at (0.6,0.6,0.6); 

			\coordinate (Unorm) at (0,0,0.6);

			\coordinate (UxN) at (0.6,-0.6,0);

			\coordinate (NxUxN) at (0.6,0.6,0);

			\draw[->, thin] (O) -- (N) node[anchor=south]{$\nu^\Xi$};

			\draw[->, thick] (O) -- (U) node[anchor=west]{$u$};

			\draw[->, thick] (O) -- (Unorm) node[anchor= east]{$(u \cdot \nu^\Xi)\nu^\Xi$};

			\draw[->, thick] (O) -- (UxN) node[anchor=south]{$u \times \nu^\Xi$};

			\draw[->, thick] (O) -- (NxUxN) node[anchor= west]{$\nu^\Xi \times (u \times \nu^\Xi)$};

			\draw[dashed] (U) -- (Unorm);
			\draw[dashed] (U) -- (NxUxN);

			\draw (-45:0.3) arc (-45:45:0.3);
			\fill (0:0.15) circle[radius=0.02]; 

			\tdplotsetrotatedcoords{-135}{-90}{45}
			\draw[tdplot_rotated_coords] (-45:0.3) arc (-45:45:0.3);
			\fill[tdplot_rotated_coords] (0:0.15) circle[radius=0.02]; 

			\shade[ball color=gray!10, opacity=0.3] (0,0,0) circle (1.1);

		\end{tikzpicture}
	}
	\vspace{-0.2cm}
\end{minipage}
where $\Xi$ is usually $\Gamma= \partial \upsilon$, $\partial \omega$ or $\partial \Omega$ and $\nu=\nu^\Xi$ denotes the corresponding
outer normal. Note that $\nu^\Gamma = \nu^{\partial \upsilon} = - \nu^{\partial \omega}$ on $\Gamma$.
We set $\llbracket f \rrbracket \coloneqq \tr_{\partial \upsilon} f - \tr_{\partial \omega} f$ for the jump
across $\Gamma$. Here and below we often omit the arguments of functions.

\subsection*{Function spaces}

We use the real, $L^2$-based Sobolev spaces $H^s$ on the domains and their boundaries. We sometimes
abbreviate $H^s(\Omega) = H^s$.
Set
\begin{align*}
	H(\div_\alpha) & \coloneqq \left\{f \in \big(L^2(\Xi)\big)^3 \bigm\vert \div (\alpha f)
	\eqqcolon \div_\alpha f \in L^2(\Xi)\right\},                                           \\
	H(\curl_\beta) & \coloneqq \left\{f \in \big(L^2(\Xi)\big)^3 \bigm\vert \curl(\beta f)
	\eqqcolon \curl_\beta f \in \big(L^2(\Xi)\big)^3\right\},
\end{align*}
for $\Xi \in \{\upsilon,\omega,\Omega\}$ and suitable matrix-valued $\alpha,\beta \colon \Xi \to \R^{3 \times 3}$.
If $\alpha=\beta = I_3$, we simply write $H(\div)$ and $H(\curl)$.
These are the maximal domains of the distributional operators $\div$ and $\curl$ in $L^2(\Xi)$ equipped with the graph norm,
respectively. (Occasionally we simplify $\big(L^2(\Xi)\big)^3$ to $L^2(\Xi)$ etc.) The kernels of these operators are written as
$H(\div 0)$ and $H(\curl 0)$ and are endowed with the $L^2$-norm.

On the boundary $\partial\Xi$, we introduce the negative fractional Sobolev space $H^{-1/2}(\partial\Xi)$ as the continuous
dual of $H^{1/2}(\partial\Xi)$ with respect to the pivot space $L^2(\partial\Xi)$. The normal trace $\tr_n$ and the tangential
traces $\tr_t, \tr_\tau$ can be extended to maps from $H(\curl)$ and $H(\div)$ to $H^{-1/2}(\partial \Xi)^3$. We denote
the kernels of these trace mappings by $H_{n0}(\div)$ and $H_{t0}(\curl)$, respectively, which coincide with the closures of
the space of test functions $C_c^\infty(\Xi)^3$ under the respective graph norms. The intersection of the kernels are abbreviated by
\[
	H_{t0}(\curl0)\coloneqq H_{t0}(\curl)\cap H(\curl0) \quad \text{and} \quad H_{n0}(\div0)\coloneqq H_{n0}(\div) \cap H(\div0),
\]
and analogously for $\div_\alpha$ and $\curl_\beta$.
Furthermore, we define the refinements
\begin{alignat*}{-1}
	H^{\partial\Xi}(\div_\alpha 0) & \coloneqq \big\{ u \in L^2(\Xi)^3 \bigm|   &  & \div_\alpha u =0,
	\langle \tr_n (\alpha u), \1 \rangle_{H^{-\frac12}((\partial\Xi)_k)}=0, \; 0\le k \le K \big\},                        \\
	H^\Sigma_{n0}(\div_\alpha0)    & \coloneqq \big\{u \in L^2(\Omega)^3 \bigm| &  & \div_\alpha u =0, \tr_n(\alpha u) =0,
	\langle \tr_n (\alpha u), \1 \rangle_{H^{-\frac12}(\Sigma_l)}=0,                                                       \\
	                               &                                            &  & l\!\in\!\{ 1, \dots, L\}\big\}.
\end{alignat*}
where the boundary $\partial\Xi$ has connected components $(\partial\Xi)_0 ,\dots, (\partial\Xi)_{K}$. (We will need the space
$H^\Sigma_{n0}(\div_\alpha0)$ only for $\Omega$. For $\Xi = \upsilon$ we write
$H^{\partial\Xi}(\div_\alpha 0) = H^{\Gamma}(\div_\alpha 0)$.)
For $\alpha= I$ we simplify $H^\Sigma_{n0}(\div 0)$ and $H^{\partial\Xi}(\div 0)$.
The divergence theorem yields $H^{\partial\Xi}(\div_\alpha 0)=H(\div_\alpha 0)$ if $\partial\Xi$ is connected,
and for simply connected $\Omega$ we have $H^\Sigma_{n0}(\div_\alpha 0)= H_{n0}(\div_\alpha 0)$.
See \cite{ciarlet18}, \cite{cessenat1996}, or \cite{dautray1990} for these and related facts.

To describe regularity on either side of an interface we introduce the space
\begin{align*}
	H_\Gamma^k = H_\Gamma^k(\Omega) & \coloneqq \big\{f \in L^2(\Omega) \bigm\vert f|_\omega \in H^k(\omega)
	\text{ and } f|_\upsilon \in H^k(\upsilon)\big\}.
\end{align*}
For an interval $J\subset \R$, we consider solutions in
\begin{align*}
	\cG_\Gamma^k = \cG_\Gamma^k(J \times \Omega) & \coloneqq \bigcap\nolimits_{j=0}^k C^j((\Tmin,\Tmax), H_\Gamma^{k-j}(\Omega))^6.
\end{align*}
If $J$ is compact, this space has the norm
$\norm{u}_{\cG_\Gamma^k} \coloneqq \sup_{t \in J} \sum_{j=0}^k \, \smallnorm{\partial_t^j u(t)}_{H^{m-j}_\Gamma}$.
Finally, ``$\lesssim$'' means that an inequality holds up to a constant not depending on time or the explicit norm of the
solution.

\section{Problem setting}\label{sec:ass}

We start by introducing the standing assumptions throughout the paper. They are summarized at the end of the section, see \ref{hypothesis-quasi}.

\subsection*{Assumptions on the coefficients}

The material parameters satisfy
\begin{equation}\label{eq:ass-parameters-1}
	\epsilon, \mu, \partial_{\xi_k}\epsilon, \partial_{\xi_k}\mu \in C^{3}(\overline \Omega_\pm \times \R^3, \Rsym),
\end{equation}
for $k \in \{1,2,3\}$, where $C^{3}(\overline{\Omega_\pm} \times \R^3, \Rsym)$ is the space
\[\big\{ \alpha \in C^{3}(\Omega_\pm \times \R^3, \Rsym)\,\big|\, \alpha|_\upsilon \in C^{3}(\overline\upsilon\times\R^3,\Rsym),
	\alpha|_\omega \in C^{3}(\overline \omega \times \R^3, \Rsym)\big\}. \]
So the coefficients are $C^3$ on either side of the interface.
They are assumed to be uniformly positive definite for vanishing fields in the sense that
\begin{equation}
	\label{eq:ass-parameters-2}
	\epsilon(x,0),\mu(x,0) \geq 2\eta>0
\end{equation}
for some $\eta> 0$ and all $x\in\Omega_\pm$. For the conductivity we require that
\begin{equation}\label{eq:positivity-condition-conductivity}
	\omega \subseteq \Omega \text{ is open}, \partial\omega\in C^5, \sigma \in C^3(\overline\omega, \Rsym)
	\text{ with } \omega \supseteq \cN_a, \text{ and } \sigma|_{\omega} \geq \sigma_0>0,
\end{equation}
for some $a > 0$.
\vspace*{-0.5cm}
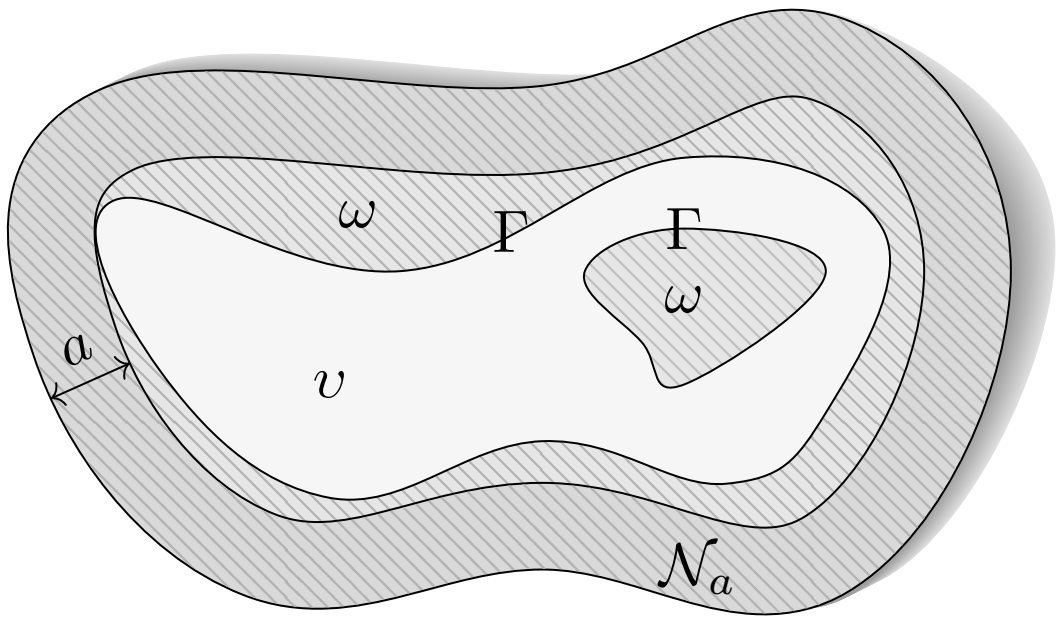
\begin{figure}[H]
	\centering
	\resizebox{0.6\textwidth}{!}{
		\def\nodeScale{1.2}
		\newcommand\decWidth{8.45pt}

		\begin{tikzpicture}[decoration={markings,
						mark= between positions 0 and 1 step 2pt with{
								\path[pattern=north west lines, pattern color=black!30] (0,-\decWidth) rectangle (3pt,\decWidth);
							}}]


			\foreach \step in {25, ..., 0} {
					\pgfmathsetmacro{\xoff}{\step*0.010-0.25}
					\pgfmathsetmacro{\yoff}{\step*0.005-0.2}
					\pgfmathtruncatemacro{\col}{75 - \step*2}
					\fill[gray!\col] plot [smooth cycle,tension=0.7, scale=1.12] coordinates {
							(\xoff+2,\yoff+1) (\xoff+4,\yoff+1.2) (\xoff+5.78,\yoff+0.85) (\xoff+6.8,\yoff+3)
							(\xoff+5.5,\yoff+4.35) (\xoff+4,\yoff+4) (\xoff+1,\yoff+4.01) (\xoff+0.7,\yoff+2.5)
						};
				}

			\draw[fill=gray!20] plot [smooth cycle,tension=0.7] coordinates {(2,1) (4,1.2) (6,1) (7,3) (6,4.5) (4,4) (1,4) (0.7,2.5)};

			\draw[double=gray!30, double distance= 17pt, postaction={decorate},pattern=north west lines, pattern color=black!30] plot [smooth cycle,tension=0.7] coordinates {(2,1) (4,1.2) (6,1) (7,3) (6,4.5) (4,4) (1,4) (0.7,2.5)};

			\draw[fill=gray!7] plot [smooth cycle,tension=0.7] coordinates {(2.5,1.4) (4,1.8) (5.3,1.5) (6,2) (6.4,3.3) (5,3.8) (3,3) (0.98,3.5) (1.3,2.3)};

			\draw[preaction={fill, gray!20}, pattern=north west lines, pattern color=black!25] plot [smooth cycle,tension=0.7] coordinates {(4.7,2.5) (5,2.2) (6,3) (5,3.3) (4.3,3) };

			\node[scale=\nodeScale] at (5.1,0.9) {$\mathcal{N}_a$};
			\node[scale=\nodeScale] at (2.5,2.2) {$\upsilon$};
			\node[scale=\nodeScale] at (3.78,3.28) {$\Gamma$};
			\node[scale=\nodeScale] at (5,3.3) {$\Gamma$};
			\node[scale=\nodeScale] at (2.7,3.4) {$\omega$};
			\node[scale=\nodeScale] at (5,2.8) {$\omega$};

			\draw[<->] (0.53,2.1)--(1.09,2.35)node[scale=\nodeScale, midway,above, sloped] {$a$};
		\end{tikzpicture}
	}
	\vspace*{-0.5cm}
	\caption{$2d$-cut through $\Omega \subset \R^3$. Note that $\supp \sigma=\overline \omega$ contains $\cN_a$.}
\end{figure}

We extend $\sigma$ discontinuously by $0$ to $\upsilon$, leading to an interface problem on $\Gamma = \partial \upsilon$
with outer unit normal $\nu^\Gamma$ of $\upsilon$. The Maxwell equations \eqref{eq:maxwell}--\eqref{eq:maxwell4} then imply
in the absence of surface currents the interface conditions (see Equation~(1.2) in \cite{schnaubelt2022})
\begin{equation}\label{eq:interface-cond}
	\begin{aligned}
		\llbracket E \times \nu^\Gamma \rrbracket & =0,                                                                                                                     & \quad\llbracket H\times\nu^\Gamma \rrbracket & =0, \\
		\llbracket \epsilon E \cdot \nu^\Gamma \rrbracket
		                                          & = \llbracket\epsilon E(t_0)\cdot\nu^\Gamma\rrbracket - \int_{t_0}^t\llbracket\sigma E\cdot\nu^\Gamma\rrbracket\dd\tau ,
		                                          & \llbracket \mu H \cdot \nu^\Gamma \rrbracket                                                                            & =0.
	\end{aligned}
\end{equation}

Note that the nonlinearities $\epsilon(x,\xi)$ and $\mu(x,\xi)$ are not required to be bounded in $\xi \in \R^3$.
Later on, however, the arguments $E$ and $H$ will be bounded and so will $\epsilon(x,E)$ as well as $\mu(x,H)$.
Typical examples are polynomial nonlinearities, see Example~2.1 in \cite{lasiecka2019} and also Chapter~1.5 in \cite{boyd2008}.

\subsection*{Assumptions on the initial values}

We take initial fields
\begin{equation*}
	(E_0,H_0) \in \left(H_\Gamma^3(\Omega)\right)^6,
\end{equation*}
satisfying the charge conditions
\begin{equation}\label{eq:div-cond-initial}
	\mu(H_0)H_0 \in H^\Sigma_{n0}(\div 0) \quad \text{and} \quad \epsilon(E_0) E_0|_{\upsilon} \in H^\Gamma(\div 0, \upsilon).
\end{equation}
The spaces in \eqref{eq:div-cond-initial} ensure that the charge conditions and magnetic boundary condition are met.
Furthermore, they exclude certain stationary solutions for non-trivial topologies of $\Omega$. (Heuristically, the subset
$\upsilon$ does not enclose electric charges and there are no currents running through loops of $\Omega$.)

Additionally, the initial values have to satisfy the \textit{compatibility conditions} of order $3$, see
\eqref {eq:compatibility-trace} and \eqref{eq:compatibility-jump}.
In order to obtain solutions of higher regularity the
time derivatives at $t=0$ have to fit to the spatial interface and boundary conditions.
Note that these conditions stay invariant under the evolution of \eqref{eq:maxwell}--\eqref{eq:maxwell4}.
In order to express the compatibility conditions we let
\begin{equation}\label{eq:def-epsilon_nl}
	\begin{alignedat}{2}
		\epsilond(x,\xi) & \coloneqq \epsilon(x,\xi) + \Big(\sum_{l=1}^3\partial_{\xi_j}\epsilon_{il}(x,\xi)\xi_l\Big)_{ij}
		                 &                                                                                                  & \eqqcolon \epsilon(x,\xi) + \epsilondnl(x,\xi), \\
		\mud(x,\xi)      & \coloneqq \mu(x,\xi) + \Big(\sum_{l=1}^3\partial_{\xi_j}\mu_{il}(x,\xi)\xi_l\Big)_{ij}
		                 &                                                                                                  & \eqqcolon \mu(x,\xi)+ \mudnl(x,\xi),
	\end{alignedat}
\end{equation}
for $x \in \Omega_\pm$, $\xi \in \R^3$, i.e., $\epsilond$ equals $\epsilon$ up to a perturbation $\epsilondnl$,
which is small for small $\xi$. On the initial values we further impose the definiteness assumption
\begin{equation}
	\label{eq:definiteness}
	\epsilon(\cdot,E_0), \epsilond(\cdot,E_0), \mu(\cdot,H_0), \mud(\cdot,H_0) \geq \tfrac32 \eta
\end{equation}
on $\Omega_\pm$.
By \eqref{eq:ass-parameters-1} and \eqref{eq:ass-parameters-2} this is true if $\norm{E_0,H_0}_{H^3_\Gamma} \leq \delta_\eta^0$
for a sufficiently small $\delta_\eta^0>0$. We set
\begin{equation}
	\label{eq:def-derivatives-for-CC}
	\begin{aligned}
		E_0^1 & \coloneqq\epsilond\left(E_0\right)^{-1}\left[\curl H_0-\sigma E_0\right],                       \\*
		E_0^2 & \coloneqq\epsilond\left(E_0\right)^{-1}\left[\curl H_0^1-\sigma E_0^1
		-\Big(\sum_{\ell=1}^3 \partial_{\xi_\ell}{\epsilond}_{i,j}(E_0)(E^{1}_0)_\ell \Big)_{i,j} E^1_0\right], \\
		H_0^1 & \coloneqq-\mud\left(H_0\right)^{-1} \curl E_0,                                                  \\*
		H^2_0 & \coloneqq - \mud\left(H_0\right)^{-1} \left[\curl E^1_0
			- \Big(\sum_{\ell=1}^3 \partial_{\xi_\ell}{\mud}_{i,j}(H_0)(H^{1}_0)_\ell \Big)_{i,j} H^1_0 \right].
	\end{aligned}
\end{equation}
Requiring the above fields to be in $H(\curl)$, we assume the boundary
\begin{equation}\label{eq:compatibility-trace}
	\tr_{t,\partial\Omega} E_0 = \tr_{t,\partial\Omega} E_0^1 = \tr_{t,\partial\Omega} E_0^2 = 0,
\end{equation}
and the interface conditions
\begin{equation}\label{eq:compatibility-jump}
	\begin{alignedat}{3}
		\llbracket E_0 \times \nu^\Gamma\rrbracket
		 & = \llbracket E^1_0 \times \nu^\Gamma\rrbracket &  & = \llbracket E^2_0 \times \nu^\Gamma\rrbracket &  & = 0, \\
		\llbracket H_0 \times \nu^\Gamma\rrbracket
		 & = \llbracket H^1_0 \times \nu^\Gamma\rrbracket &  & = \llbracket H^2_0 \times \nu^\Gamma\rrbracket &  & = 0.
	\end{alignedat}
\end{equation}
Surface currents are thus excluded (i.e., we let $J_\Sigma = 0$ in \cite{schnaubelt2022}).

\subsection*{Existence of solutions}

Under these conditions, Theorem~7.3 of \cite{schnaubelt2022} yields a unique solution
\begin{equation}
	\label{eq:solution}
	u=(E,H) \in \cG_\Gamma^3((\Tmin,\Tmax) \times \Omega),
\end{equation}
to \eqref{eq:maxwell}--\eqref{eq:maxwell4}, which takes values in $B_{H^3_\Gamma}(0,\delta_\eta^0)$ and depends
continuously on data. Hence, as $(E_0,H_0)=0$ has the maximal existence time $\Tmax=\infty$,
for fixed $T \in\left(0, \infty\right)$ there is a number $\tilde r>0$ such that for all data
$\tilde{u}_0 \in B_{H_\Gamma^3}(0,\tilde r)$, satisfying the compatibility/interface conditions above, we have
$\Tmax\left(\tilde{u}_0\right)>T$.

Let $0 < \delta_z \leq \delta_\eta \leq \min\{1, \delta_\eta^0\}$. The number $\delta_\eta$ will be fixed in
\eqref{eq:non-trapping-tilde1} and \eqref{eq:non-trapping-tilde2} below. Since flow map
$\tilde{u}_0 \mapsto u(\cdot ; \tilde{u}_0)$ is continuous from
$ B_{H_\Gamma^3}(0, \tilde r)$ to $\cG_\Gamma^3\left((0, T) \times \Omega\right),$
for sufficiently small $r(\delta_z)\in(0 ,\tilde r)$ and data with
\begin{equation}\label{eq:initial_data_bounded}
	\norm{(E_0,H_0)} _{H_\Gamma^3(\Omega)} \leq r(\delta_z)
\end{equation}
satisfying the compatibility and interface conditions, there exists a unique solution $(E,H) \in \cG_\Gamma^3([0,\Tmax))$
such that $\Tmax > T$,
\begin{equation}\label{eq:sol-bound}
	\max_{0\leq j \leq 3} \norm{\pt^j(E(t),H(t))}_{H_\Gamma^{3-j}(\Omega)}\leq \delta_z \le \delta_\eta
	\quad \text{for } 0\leq t\leq T.
\end{equation}
We set
\begin{equation*}
	T_* \coloneqq \sup \big\{\hat T \in [0,\Tmax] \,\big|\, \eqref{eq:sol-bound} \text{ holds for } t \in [0,\hat T]\big\} .
\end{equation*}
From the blow-up condition in Theorem~7.3 of \cite{schnaubelt2022} we infer that $T_* < \infty$ implies
\begin{equation}\label{eq:sol-normBounded}
	\max_{0\leq j \leq 3}\norm{\pt^j(E(T_*),H(T_*))}_{H_\Gamma^{3-j}(\Omega)} = \delta_z.
\end{equation}

Furthermore, the condition \eqref{eq:div-cond-initial} on the initial values is invariant under the evolution of the system,
which follows from \eqref{eq:maxwell} and Proposition~IX.1.3 in \cite{dautray1990} on $\upsilon$ for the electric field and
Proposition~IX.1.4 in \cite{dautray1990} for the magnetic field.
\begin{lemma}
	\label{lem:invariant-under-evo}
	Let $(E_0,H_0) \in \left(H_\Gamma^3(\Omega)\right)^6$ satisfy \eqref{eq:div-cond-initial}. Then the solution $(E,H)$ of
	\eqref{eq:solution} satisfies
	\[
		\epsilon E|_{\upsilon} \in H^\Gamma(\div 0, \upsilon) \quad \text{and} \quad\mu H \in H^\Sigma_{n0}(\div 0)
		\qquad \text{for all } t \in (\Tmin,\Tmax).
	\]
\end{lemma}



\subsection*{Non-Trapping conditions}

In order to show decay, we impose the following \textit{non-trapping condition} in the non-dissipative region and on a
neighborhood of $\Gamma$. Let $x_0 \in \R^3$ and $m(x)\coloneqq x-x_0$. We assume
\begin{equation}\label{eq:non-trapping}
	\begin{alignedat}{4}
		\tilde\epsilon(x,\xi)|_{\xi=0} & \coloneqq (\epsilon +(m(x)\cdot\nabla_x)\epsilon) & (x,\xi)|_{\xi=0}                &
		                               &                                                   & \geq 2\tilde\eta \epsilon(x,0),   \\
		\tilde\mu(x,\xi)|_{\xi=0}      & \coloneqq (\mu + (m(x)\cdot \nabla_x)\mu)         & (x,\xi)|_{\xi=0}                &
		                               &                                                   & \geq 2\tilde\eta \mu(x,0),
	\end{alignedat}
\end{equation}
for $x \in \upsilon \cup \cN_b(\Gamma)$ and some constants $\tilde \eta, b> 0$.
This condition conceptually limits back reflection of the solutions towards $x_0$, allowing the fields to move into
$\cN_a$, where we assumed dissipation to occur. We transfer \eqref{eq:non-trapping} to $\epsilond$ and $\mud$ using
\begin{equation}\label{eq:def-tilde-epsilond}
	\tildeepsilond(x,\xi) \coloneqq (\epsilond + m \cdot \nabla_x \epsilond) (x,\xi), \quad
	\tildemud(x,\xi) \coloneqq (\mud + m \cdot \nabla_x \mud)(x,\xi).
\end{equation}
The non-trapping condition~\eqref{eq:non-trapping} implies
\begin{alignat*}{4}
	\tildeepsilond(x,0) & = (\tilde \epsilon + \epsilondnl + m \cdot \nabla_x \epsilondnl )(x,0) &
	                    & = \tilde \epsilon(x,0)                                                 &   & \geq 2 \tilde\eta\epsilond(x,0), \\
	\tildemud(x,0)      & = (\tilde \mu + \mudnl + m \cdot \nabla_x \mudnl )(x,0)                &
	                    & = \tilde \mu(x,0)                                                      &   & \geq 2 \tilde\eta \mud(x,0),
\end{alignat*}
at $(x,0)$ for $x \in \upsilon \cup \cN_b(\Gamma)$.

Furthermore, we require the \textit{non-trapping interface condition}
\begin{align}
	\label{eq:non-trapping-interface}
	(m \cdot \nu^\Gamma)\llbracket \epsilond(\cdot,E) \rrbracket & \leq 0, &
	(m \cdot \nu^\Gamma)\llbracket \mud(\cdot,H) \rrbracket      & \leq 0, &   &
	\quad \text{on }[0,T_*)\times\Gamma,
	\intertext{
		to also limit reflections at the interface. Note that this condition involves the solution and might therefore be hard to verify.
		However, it is satisfied for small solutions if we assume
	}
	\label{eq:special-non-trapping-interface}
	(m \cdot \nu^\Gamma)\llbracket \epsilon(\cdot,0) \rrbracket  & < 0,    &
	(m \cdot \nu^\Gamma)\llbracket \mu(\cdot,0) \rrbracket       & < 0,    &   &
	\quad \text{on } \Gamma.
\end{align}
If this condition is violated back reflections can occur, see Example~4.1.3 in \cite{nutt2026-diss}.

To simplify some expressions, we define
\begin{align*}
	\epsilonk(t, x) & \coloneqq \begin{cases}
		                            \epsilon(x, E(t,x)),  & \quad k =0,      \\
		                            \epsilond(x, E(t,x)), & \quad k = 1,2,3,
	                            \end{cases}\qquad t \in [0,\Tmax), \ x \in \Omega_\pm, \\
	\muk(t, x)      & \coloneqq \begin{cases}
		                            \mu(x, E(t,x)),  & \quad k =0,      \\
		                            \mud(x, E(t,x)), & \quad k = 1,2,3,
	                            \end{cases}\qquad t \in [0,\Tmax), \ x \in \Omega_\pm.
\end{align*}

We assume without loss of generality that $\eta\leq \tilde\eta$, cf.\ \eqref{eq:ass-parameters-2},
and fix $\delta_\eta\in\left(0,\min\{1, \delta_\eta^0\}\right)$, see \eqref{eq:definiteness}, such that
$(E,H)$ from \eqref{eq:sol-bound} satisfies
\begin{align}\label{eq:non-trapping-tilde1}
	 & \hspace{8mm}\epsilonk \geq \eta, \hspace{23mm}\muk \geq \eta,
	 &                                                                                                              & \text{on } [0,T_*) \times \Omega_\pm,                              \\
	\label{eq:non-trapping-tilde2}
	 & \begin{aligned}
		   \tilde \epsilon(\cdot,E) & \geq \eta \epsilon(\cdot,E),  & \tilde \mu(\cdot,H) & \geq \eta \mu(\cdot,H),  \\
		   \tildeepsilond(\cdot,E)  & \geq \eta \epsilond(\cdot,E), & \tildemud(\cdot,H)  & \geq \eta \mud(\cdot,H),
	   \end{aligned}
	 &                                                                                                              & \text{on } [0,T_*) \times \bigl(\upsilon \cup \cN_b(\Gamma)\bigr).
\end{align}

\smallskip

\noindent We summarize our assumptions:
\begin{enumerate}[leftmargin=25pt, label=\textup{($\text{H}$)}]
	\item \label{hypothesis-quasi} \quad \parbox[m]{\dimexpr\linewidth-1em\relax}{
		      Let $\Omega\subset \R^3$ be open, bounded and multiply connected, and $\partial\Omega\in C^5$ have connected components
		      $(\partial \Omega)_0,\dots, (\partial \Omega)_K$. Let $\omega \subseteq \Omega$ with $C^5$-boundary such that the
		      \textit{conductivity} $\sigma$ is supported on $\omega$ and satisfies \eqref{eq:positivity-condition-conductivity}.
		      We assume the \textit{permittivity} $\epsilon$ and \textit{permeability} $\mu$ to satisfy \eqref{eq:ass-parameters-1}
		      and \eqref{eq:ass-parameters-2}. Let the initial values $(E_0,H_0) \in \left(H_\Gamma^3(\Omega)\right)^6$ satisfy
		      \eqref{eq:div-cond-initial} as well as the \textit{compatibility} and \textit{smallness conditions}
		      \eqref{eq:compatibility-trace}, \eqref{eq:compatibility-jump}, and \eqref{eq:initial_data_bounded}.
		      In \eqref{eq:initial_data_bounded} we have $0 < \delta_z\leq \delta_\eta \leq \min\{1,\delta_\eta^0\}$ so that the coefficients
		      are positive definite, see \eqref{eq:non-trapping-tilde1}, and the \textit{non-trapping condition} \eqref{eq:non-trapping-tilde2}
		      on $\upsilon \cup \cN_b(\Gamma)$ is satisfied. In addition, we assume the \textit{non-trapping interface condition}
		      \eqref{eq:non-trapping-interface}.}
\end{enumerate}

\section{Reformulations, auxiliary results and energy estimates}

\subsection*{Reformulations of Maxwell equations}

Let $k\in \{0,1,2,3\}$ and (H) be true.
We will study the time derivatives of the solutions $\pt^k (E,H) \in \cG^{3-k}_\Gamma$ of
\eqref{eq:maxwell}--\eqref{eq:maxwell4} separately. To this end, we define the commutator terms
\begin{equation*}
	\begin{aligned}
		f_0 & =f_1=0, & f_2                                                                                  & = (\pt \epsilond(\cdot,E))\pt E,
		    & f_3     & =(\partial_t^2 \epsilond(\cdot,E))\pt E + 2 (\pt \epsilond(\cdot,E)) \partial_t^2 E,                                    \\
		g_0 & =g_1=0, & g_2                                                                                  & = (\pt \mud(\cdot,H))\pt H,
		    & g_3     & =(\partial_t^2 \mud(\cdot,H))\pt H + 2 (\pt \mud(\cdot,H)) \partial_t^2 H.
	\end{aligned}
\end{equation*}
By Lemma~\ref{lem:error-estimate} even the maps $f_3$ and $g_3$ belong to $\cG_\Gamma^1$ and are small.
With the definitions above, the time-derived solutions then satisfy
\begin{equation}\label{eq:commute-pt-epsilon}
	\begin{aligned}
		\partial_t^k (\epsilon E ) = \epsilonk \partial_t^{k} E + f_k, \quad
		\partial_t^k (\mu H )= \muk \partial_t^{k} H + g_k \quad \text{on } \Omega, \\
	\end{aligned}
\end{equation}
For $k\in\{0,1,2,3\}$ they thus solve the Maxwell system
\begin{align}\label{pde:pt-maxwell}
	\begin{split}
		\pt (\epsilonk \pt^k E) & = \curl \pt^k H - \sigma \pt^k E - \pt f_k, \\
		\pt( \muk \pt^k H)      & = -\curl \pt^k E - \pt g_k,
	\end{split}\qquad &                                                          & t\geq 0, x \in \Omega_\pm,                                                            \\
	\label{pde:pt-maxwell-IC} \llbracket \pt^k E \times \nu^\Gamma\rrbracket & = 0,
	\quad \llbracket \pt^k H \times \nu^\Gamma\rrbracket = 0, \qquad         &                                                          & t\geq 0, x \in \Gamma,         \\
	\label{pde:pt-maxwell_BC} \tr_t \pt^k E                                  & = 0, \quad \tr_n \hat\mu_k \pt^k H = - \tr_n g_k, \qquad
	                                                                         &                                                          & t\geq 0, x \in \partial\Omega.
\end{align}
\begin{remark}\label{rem:tr}
	Here we use that the tangential traces of functions $(u,v)$ in $C([0,T],L^ 2(\omega))$ with
	$(\partial_t(\epsilonk u)-\curl v, \partial_t(\muk v)+\curl u)\in\! L^2((0,T)\times \omega)$
	belong to $L^2((0,T),H^{-\frac12}(\partial \omega))$, and that this trace operator is continuous.
	The same applies to $\upsilon$. See Section~2.1 in \cite{spitz2017} and Proposition~2.1.60 in \cite{ciarlet18}.
\end{remark}
Below we work with regularizations of the solutions, see Remark~\ref{rem:approx-cG4}, mostly avoiding this trace result.
With the error terms
\begin{align*}
	\tilde f_0 \coloneqq 0, \quad \tilde f_k \coloneqq \sum_{j=1}^k \binom{k}{j} \pt^j \epsilond \pt^{k+1-j}E, \quad
	\tilde g_0 \coloneqq 0, \quad \tilde g_k & \coloneqq \sum_{j=1}^k \binom{k}{j} \pt^j \mud \pt^{k+1-j}H,
\end{align*}
one can also derive from \eqref{eq:maxwell}--\eqref{eq:maxwell4} the system
\begin{align}\label{pde:tilde1}
	\begin{split}
		\epsilond \pt (\pt^k E) & = \curl \pt^k H - \sigma \pt^k E - \tilde f_k, \\
		\mud \pt(\pt^k H)       & = -\curl \pt^k E - \tilde g_k,
	\end{split}\qquad &                                                                       & t\geq 0, x \in \Omega_\pm,                                                                                            \\
	\label{pde:tilde2}\llbracket \pt^k E \times \nu^\Gamma\rrbracket
	                                                                            & = 0, \quad \llbracket \pt^k H \times \nu^\Gamma\rrbracket = 0, \qquad &                            & t\geq 0, x \in \Gamma,         \\
	\label{pde:tilde3}\tr_t \pt^k E                                             & = 0, \qquad                                                           &                            & t\geq 0, x \in \partial\Omega.
\end{align}
Both are needed below.
We apply the divergence to \eqref{eq:maxwell} and \eqref{pde:tilde1}, obtaining
\begin{equation}\label{eq:div-evolution}
	\begin{aligned}
		\div (\epsilon E)        & = \div(\epsilon(E_0) E_0) - \int_0^t \div(\sigma E) \dd\tau, & \div(\mu H)               & = 0, \\
		\div (\epsilond \pt^k E) & = -\div(\sigma \pt^{k-1}E +\tilde f_{k-1}),
		                         & \div(\mud \pt^k H)                                           & = - \div \tilde g_{k-1} ,
	\end{aligned}
\end{equation}
on $\Omega_\pm$ for $k \in\{1,2,3\}$. We observe that for $k=3$ the right-hand side of \eqref{eq:div-evolution} lies in
$L^2(\Omega)$ by Lemma~\ref{lem:error-estimate}. Therefore, the term $(\epsilonk \pt^k E,\muk \pt^k H)$ is contained in
$H(\div, \Omega_\pm)^2$. By the results mentioned in Section~\ref{sec:notation}, 
its normal trace on either side of the interface exists in $H^{-1/2}(\Gamma)$, allowing us to derive interface conditions
for the normal jump in addition to \eqref{eq:interface-cond}. Let $k \in \{1,2,3\}$. For $\psi \in C^\infty_c(\Omega)$,
the divergence theorem, \eqref{pde:tilde1}, and \eqref{pde:pt-maxwell-IC} yields
\begin{align*}
	\Big\langle & \llbracket (\epsilond\pt^k E + \sigma \pt^{k-1} E +\tilde f_{k-1}) \cdot
	\nu^\Gamma \rrbracket, \psi \Big\rangle_{H^{-1/2}(\Gamma)\times H^{1/2}(\Gamma)}                                                           \\
	            & = \int_{\Omega_\pm} \left( \psi \div (\epsilond\pt^k E + \sigma \pt^{k-1} E +\tilde f_{k-1})
	+ \nabla \psi \cdot (\epsilond\pt^k E + \sigma \pt^{k-1} E +\tilde f_{k-1}) \right) \dd x                                                  \\
	            & = \int_{\Omega_\pm} \left( \psi \div (\curl \pt^{k-1} H) + \nabla \psi \cdot (\curl \pt^{k-1} H) \right) \dd x
	= \int_{\Omega_\pm} \nabla \psi \cdot (\curl \pt^{k-1} H) \dd x                                                                            \\
	            & = \langle \nabla\psi, \llbracket\pt^{k-1} H \times \nu^\Gamma\rrbracket\rangle_{H^{-1/2}(\Gamma)\times H^{1/2}(\Gamma)} = 0.
\end{align*}
This means that
\begin{equation}\label{eq:normal-jump}
	\left\llbracket \left(\epsilond\pt^k E\right) \cdot \nu^\Gamma \right\rrbracket
	= \tr_{\partial\omega}\left( \left(\sigma \pt^{k-1} E \right)\cdot \nu^\Gamma\right)
	- \left\llbracket \tilde f_{k-1}\cdot \nu^\Gamma\right\rrbracket \qquad \text{for }k \geq 1,
\end{equation}
where the equality holds in $H^{1/2}(\Gamma)$ due to the regularity of the right-hand side.
Likewise, one can show the analogous interface condition for the jump in normal direction of the magnetic field
\begin{equation*}
	\left\llbracket \left(\mud\pt^k H\right) \cdot \nu^\Gamma \right\rrbracket
	= - \left\llbracket \tilde g_{k-1}\cdot \nu^\Gamma\right\rrbracket \qquad \text{for }k \geq 1.
\end{equation*}
Equations~\eqref{pde:pt-maxwell_BC}, \eqref{pde:tilde1}, and Lemma~\ref{lem:characterization-H-div-H-curl} thus result in
\begin{equation}
	\label{eq:div-condition-H-kgeq1}
	\mud\pt^k H + \tilde g_{k-1} \in H_{n0}(\div 0) \qquad \text{for }k \geq 1.
\end{equation}

\subsection*{Energies and error terms}

We set
\begin{align}\label{eq:def-edz}
	\notag e_k(t)                                                    & \coloneqq \tfrac12 \max_{0\leq j \leq k} \Big(\norm{\hat \epsilon_j^{1/2}\partial_t^jE(t)}_{L^2(\Omega)}^2
	+\norm{\hat \mu_j^{1/2}\partial_t^j H(t)}_{L^2(\Omega)}^2 \Big), & e \coloneqq e_3,                                                                                                             \\
	d_k(t)                                                           & \coloneqq \max_{0\leq j \leq k} \norm{\sigma^{1/2}\partial_t^jE(t)}_{L^2(\Omega)}^2,                       & d\coloneqq d_3, \\
	\notag z_k(t)                                                    & \coloneqq \max_{0\leq j \leq k} \Big(\norm{\partial_t^jE(t)}_{H_\Gamma^{k-j}(\Omega)}^2
	+\norm{\partial_t^j H(t)}_{H_\Gamma^{k-j}(\Omega)}^2 \Big),      & z \coloneqq z_3.
\end{align}
We stress that $z$ can be made small by \eqref{eq:sol-normBounded} and that the dissipation $d$
is controlled via the energy $e$ in Lemma~\ref{lem:energy-decreasing}.
The following estimates provide smallness of the commutator terms and
can be shown as in \cite{lasiecka2019}, see also Lemma~4.2.1 in \cite{nutt2026-diss}.

\begin{lemma}[Error estimates]\label{lem:error-estimate}
	Assume (H). We have
	\begin{align*}
		\norm{\epsilonk(E)(t)}_{\infty},\norm{\muk(t)}_{\infty},\bignorm{\inv{\epsilonk(E)}(t)}_{\infty},
		\bignorm{\inv{\muk(H)}(t)}_{\infty}                      & \leq c,                                             \\
		\norm{\partial^\alpha \epsilonk(E)(t)}_{L^2(\Omega)},\norm{\partial^\alpha \muk(H)(t)}_{L^2(\Omega)}
		                                                         & \leq c\big(z_k^{1 / 2}(t)+\delta_{\alpha_0=0}\big), \\
		\max _{k \in\{0,1,2,3\}, j \in\{0,1\}}\left(\bignorm{\partial_t^j f_k(t)}_{H_\Gamma^{4-j-k}}
		+\bignorm{\partial_t^j g_k(t)}_{H_\Gamma^{4-j-k}}\right) & \leq c z(t),                                        \\
		\max _{k \in\{0,1,2,3\}}\Big(\bignorm{\tilde f_k(t)}_{H_\Gamma^{3-k}}
		+\bignorm{\tilde g_k(t)}_{H_\Gamma^{3-k}}\Big)           & \leq c z(t),                                        \\
		\max _{k \in\{0,1,2\}}\Big(\bignorm{\pt \tilde f_k(t)}_{H_\Gamma^{2-k}}
		+\bignorm{\partial_t\tilde g_k(t)}_{H_\Gamma^{2-k}}\Big) & \leq c z(t).
	\end{align*}
	Here the multiindex $\alpha \in \N_0 \times \N_0^3$ satisfies $\abs{\alpha} = k \leq 3$ with $\alpha_0\in \N_0$
	corresponding to the time derivative $\partial^{(\alpha_0,0,0,0)} = \pt^{\alpha_0}$,
	and $(\alpha_1,\alpha_2,\alpha_3)\in \N_0^3$ to the spatial derivatives.
\end{lemma}

These inequalities are often applied implicitly in combination with the following lemma,
see Lemma~2.1 in \cite{schnaubelt2022} and the proof of Lemma~2.22 in \cite{spitz2017}.
It is shown by combining the usual Sobolev embeddings with Hölder's inequality.
\begin{lemma}\label{lem:sobolev-algebra-property}
	Let $m_1,m_2 \in \N_0$ with $m_1 \geq m_2$ and $m_1 \geq 2$.

	\smallskip

	a) Let $k\in\{0,\dots,m_1\}$, $f \in H_\Gamma^{m_1-k}(\Omega)$ and $g \in H_\Gamma^{k}(\Omega)$. Then
	$fg \in L^2(\Omega)$ satisfies $\norm{fg}_{L^2} \lesssim \norm{f}_{H_\Gamma^{m_1-k}}\norm{g}_{H_\Gamma^{k}}.$

	\smallskip

	b) Let $f \in H_\Gamma^{m_1}(\Omega)$ and $g \in H_\Gamma^{m_2}(\Omega)$. Then $fg \in H_\Gamma^{m_2}(\Omega)$ satisfies
	$\norm{fg}_{H_\Gamma^{m_2}}\lesssim\norm{f}_{H_\Gamma^{m_1}}\norm{g}_{H_\Gamma^{m_2}}.$
\end{lemma}

Integration by parts shows the next result where for $k=3$ one needs an additional approximation argument
based on Theorem~\ref{thm:approx}), cf.\ Section~4 in \cite{lasiecka2019}.

\begin{lemma}[Energy estimates]\label{lem:energy-decreasing}
	Assume (H).	The energies $e_k$ are decreasing up to an error term, namely
	\begin{equation*}
		e_k(t) + \int_s^t d_k(\tau) \dd \tau \leq e_k(s) + c \int_s^t z^{3/2}(\tau) \dd \tau, \qquad
		t\in[0,T_*), \	k \in \{0,1,2,3\}.
	\end{equation*}
\end{lemma}

Because of the time integral in \eqref{eq:div-evolution}, for instance, we cannot treat $(E,H)$ directly.
Instead we use the next inequality that controls the usual energy through that of the time derivatives.
The proof is the same as that of Lemma~4.3 in \cite{nutt2026}. It needs our assumption
that $\sigma$ either vanishes or is uniformly positive.

\begin{lemma}\label{lem:cE-estimate}
	Assume (H).	For all $0<t<T_*$ we have $e_0(t) \lesssim e_1(t)$.
\end{lemma}

\section{Helmholtz decomposition}
\label{sec:helm}

We split our solution $\pt^k(E,H) \in \cG^{3-k}_\Gamma$ into a $\div_{\epsilond}$-free and a $\curl$-free part,
adjusting the interface and boundary conditions in such a way that the curl operator behaves well with the fields.
Since we only show the observability estimate for $k \geq 1$, it suffices to consider $k\in \{1,2,3\}$, simplifying the
notations slightly. For the case of $k=0$, we then rely on Lemma~\ref{lem:cE-estimate} in the proof of exponential decay.

\subsection*{Electric case}
For the $\curl$-free part, we consider the elliptic interface problem
\begin{align}
	\label{eq:def-pk-eq} \div(\epsilond\nabla p^k)                                                                         & = \div(\epsilond \pt^k E ) = -\div(\sigma \pt^{k-1}E + \tilde f_{k-1})
	\quad                                                                                                                  &                                                                        & \text{ on } \Omega_\pm,                                    \\*
	\label{eq:def-pk-IC0} \llbracket p^k \rrbracket                                                                        & =0                                                                     &                         & \text{ on } \Gamma,              \\*
	\label{eq:def-pk-IC1}\llbracket \tr_n (\epsilond\nabla p^k) \rrbracket                                                 & = \llbracket \tr_n (\epsilond \pt^k E) \rrbracket
	                                                                                                                       &                                                                        & \text{ on } \Gamma,                                        \\*
	\label{eq:def-pk-BC0} p^k                                                                                              & = 0                                                                    &                         & \text{ on } (\partial \Omega)_0, \\*
	\label{eq:def-pk-BC1}
	\begin{split} p^k & \text{ constant on } (\partial \Omega)_\ell,\text{ \ and \ }                                  \\
              0   & =\langle\tr_n\epsilond (\nabla p^k-\pt^k E),\1\rangle_{H^{-\frac12}((\partial \Omega)_\ell)},
	\end{split} \qquad & \forall\+ \ell                                                         & \in \{1,\dots,K\}.
\end{align}
In \eqref{eq:def-pk-eq} we use \eqref{pde:tilde1} for $k-1$. Observe that the right-hand side of \eqref{eq:def-pk-eq} belongs to
$L^2(\Omega)$, so that the first-order conditions on $\Gamma$ and $(\partial \Omega)_\ell$ in \eqref{eq:def-pk-IC1} and
\eqref{eq:def-pk-BC1} are well-defined in $H^{-1/2}$. For the weak formulation, we use the space
\[
	\cH_c \coloneqq \big\{\varphi \in H^1(\Omega) \bigm| \varphi = 0 \text{ on }(\partial \Omega)_0, \
	\varphi \text{ constant on } (\partial \Omega)_\ell, \ \forall\+ \ell\in\{1,\dots K\}\big\},
\]
endowed with the norm of $H^1(\Omega)$.
We introduce the forms
\begin{align*}
	B[u,v]  & = B_{E,k}[u,v] \coloneqq \int_\Omega \nabla u \cdot \epsilond \nabla v \dd x,                  \\
	\ell(v) & =\ell_{E,k}(v) \coloneqq -\int_\Omega (\sigma \pt^{k-1}E + \tilde f_{k-1})\cdot \nabla v \dd x
\end{align*}
for $u,v \in \cH_c$. Then the weak formulation of \eqref{eq:def-pk-eq}--\eqref{eq:def-pk-BC1} reads as
\begin{equation}\label{eq:pk-weak-formulation}
	B_{E,k}[p^k,\varphi] = \ell_{E,k}(\varphi) \quad \text{for all }\varphi \in \cH_c.
\end{equation}

Clearly, $B$ and $\ell$ are bounded on $\cH_c$, and $B$ is coercive by Theorem~13.6.9 in \cite{tucsnak2009}.
Lax--Milgram thus yields a unique solution $p^k=p^k(t) \in \cH_c$ of \eqref{eq:pk-weak-formulation} with
\begin{align}\label{eq:est-pk}
	\smallnorm{p^k}_{H^1(\Omega)} & \lesssim \norm{\ell}_{\cH_c^*} = \sup_{\norm{\varphi}_{\cH_c}=1} \abs{\ell_{E,k}(\varphi)}
	\lesssim \norm{\sigma \pt^{k-1}E + \tilde f_{k-1}}_{L^2(\Omega)}
	\lesssim d^\frac12 + z,
\end{align}
where we also use Lemma~\ref{lem:error-estimate}.
Since the map $\cG^3 \to \cH_c^*; E\mapsto \ell_{E,k}$, is continuous and bounded, we obtain a continuous and bounded mapping
$\cG^3 \to \cH_c$; $E\mapsto p^k$.

We verify that the solution to \eqref{eq:pk-weak-formulation} solves \eqref{eq:def-pk-eq}--\eqref{eq:def-pk-BC1}.
Since $p^k\in \cH_c$, it satisfies \eqref{eq:def-pk-BC0} and the first part of \eqref{eq:def-pk-BC1} in $H^{1/2}$, as well as
\eqref{eq:def-pk-IC0} by Lemma~\ref{lem:characterization-H-div-H-curl}. Note that \eqref{eq:def-pk-eq} is true and
$\epsilond \pt^k E - \epsilond \nabla p^k$
belongs to $H(\div 0,\Omega)$, since \eqref{eq:pk-weak-formulation} holds in particular for $H^1_0(\Omega) \subseteq \cH_c$.
Therefore, again by Lemma~\ref{lem:characterization-H-div-H-curl}, $p^k$ also fulfills \eqref{eq:def-pk-IC1} in $H^{-1/2}$.
In particular, the last boundary condition in \eqref{eq:def-pk-BC1} is well-defined.
It follows from \eqref{eq:pk-weak-formulation} by inserting
$\varphi_\ell \in \cH_c$ with $\varphi_\ell|_{(\partial \Omega)_j} = \delta_{\ell j}$ for $\ell\neq 0$ via
\[
	0= \int_\Omega (\epsilond \pt^k E - \epsilond \nabla p^k)\cdot \nabla \varphi_\ell \dd x
	= \langle\tr_n (\epsilond \pt^k E - \epsilond \nabla p^k),\1\rangle_{H^{-\frac12}((\partial \Omega)_\ell)}.
\]
Moreover, based on elliptic regularity, Lemma~\ref{lem:pkqk-H^2} below shows that $p^k \in H^2_\Gamma(\Omega)$.

\subsection*{Magnetic case}
For $H$ we solve the following Neumann problem to correct for the jump over the interface
\begin{equation}\label{eq:def-qk}
	\begin{aligned}
		\div(\mud\nabla q^k)                         & = \div(\mud\pt^k H) = - \div(\tilde g_{k-1})\qquad &   & \text{ on } \Omega, \\
		\llbracket q^k \rrbracket                    & = 0                                                &   & \text{ on } \Gamma, \\
		\llbracket \tr_n (\mud\nabla q^k) \rrbracket & = \llbracket \tr_n (\mud \pt^k H) \rrbracket       &
		                                             & \text{ on } \Gamma,                                                          \\
		\tr_n (\mud \nabla q^k)                      & = \tr_n(\mud \pt^k H) =-\tr_n \tilde g_{k-1}       &
		                                             & \text{ on } \partial\Omega.
	\end{aligned}
\end{equation}
We define the space $\cH_n(\Omega) \coloneqq \{ \varphi \in H^1(\Omega) \mid \int_\Omega \varphi \dd x = 0 \}$ of functions
with zero mean. As for \eqref{eq:def-pk-eq}--\eqref{eq:def-pk-BC1}, we define the forms
\begin{align*}
	B[u,v] & = B_{H,k}[u,v] \coloneqq \int_\Omega \nabla u \cdot \mud \nabla v \dd x,\qquad
	\ell(v)=\ell_{H,k}(v) \coloneqq -\int_\Omega \tilde g_{k-1}\cdot \nabla v \dd x
\end{align*}
for $u,v \in H^1(\Omega)$. The weak formulation of \eqref{eq:def-qk} is given by
\begin{equation}\label{eq:qk-weak-formulation}
	B_{H,k}[q^k,\varphi] = \ell_{H,k}(\varphi) \quad \text{for all }\varphi \in H^1(\Omega).
\end{equation}
By Poincar\'e's inequality, see for example Section~5.8.1 in \cite{evans1998}, the form $B$ is coercive on $\cH_n(\Omega)$. Thus,
as before Lax--Milgram provides a unique weak solution $q^k \in \cH_n(\Omega)$ satisfying \eqref{eq:qk-weak-formulation} with
\begin{align}\label{eq:est-qk}
	\smallnorm{q^k}_{H^1(\Omega)} & \lesssim \norm{\ell}_{(\cH_n)^*} = \sup_{\norm{\varphi}_{\cH_n}=1} \abs{\ell_{H,k}(\varphi)}
	\lesssim z.
\end{align}
Again, the map $\cG^3 \to \cH_n$; $H\mapsto q^k$, is continuous and bounded, and $q^k$ solves \eqref{eq:def-qk}.

\subsection*{Estimates for \texorpdfstring{$p^k$}{p\textasciicircum k} and \texorpdfstring{$q^k$}{q\textasciicircum k}}

The following Lemma guarantees that the correction terms belong to $H^2_\Gamma(\Omega)$, since
$\div(\epsilond\pt^k E) \in L^2(\Omega)$ by Maxwell.

\begin{lemma}\label{lem:pkqk-H^2}
	Assume (H).	Let $k\in\{1,2,3\}$.
	The weak solutions $(p^k,q^k)$ of \eqref{eq:pk-weak-formulation} respectively \eqref{eq:qk-weak-formulation}
	belong to $H^2_\Gamma(\Omega)$.
\end{lemma}

\begin{proof}
	1) We focus on $p^k$ as the proof for $q^k$ is similar.
	Let $\chi_j \in C^\infty(\overline\Omega)$ be smooth cutoffs satisfying $\chi_j|_{(\partial \Omega)_\ell} = \delta_{j\ell}$
	and $\chi_j|_\Gamma = 0$ for $j \in \{0, \dots , K\}$. We set
	\[ h\coloneqq \sum\nolimits_{j=0}^K \chi_j c_j \in C^\infty(\overline\Omega) \subseteq H^2(\Omega), \]
	where $p^k|_{(\partial \Omega)_j} = c_j$.
	Let $\tilde p^k \coloneqq p^k - h$. Then $\tilde p^k$ is a weak solution of
	\begin{equation*}
		\begin{aligned}
			\div(\epsilond\nabla \tilde p^k)                        & = \div(\epsilond(\pt^k E - \nabla h)) = -\div(\sigma \pt^{k-1}E + \tilde f_{k-1}
			+ \nabla h) \eqqcolon f                                 &                                                                                  & \text{ on } \Omega,                               \\
			\llbracket \tilde p^k \rrbracket                        & = 0                                                                              &                     & \text{ on } \Gamma,         \\
			\llbracket \tr_n (\epsilond\nabla \tilde p^k)\rrbracket & = \llbracket \tr_n (\epsilond \pt^k E) \rrbracket
			= \tr_{\partial\omega}\left(\sigma \pt^{k-1}E \cdot \nu^\Gamma \right)
			- \llbracket \tr_n \tilde f_{k-1}\rrbracket \eqqcolon g &                                                                                  & \text{ on } \Gamma,                               \\
			\tilde p^k                                              & = 0                                                                              &                     & \text{ on } \partial\Omega,
		\end{aligned}
	\end{equation*}
	where we used \eqref{pde:tilde1} as well as \eqref{eq:normal-jump}. Lemma~\ref{lem:error-estimate} yields
	$f\in L^2(\Omega)$ and	$g \in H^{1/2}(\Gamma).$ By elliptic regularity for transmission problems, see Theorem~5.3.7
	in \cite{costabel2010} or Theorem~A.11 in \cite{breschPHD2026}, the weak solution $\tilde p^k$ lies in $H^2_\Gamma(\Omega)$
	and satisfies
	\[ \bignorm{\tilde p^k}_{H^2_\Gamma} \lesssim \norm{f}_{L^2(\Omega)} + \norm{g}_{H^{1/2}(\partial\Omega)}
		+ \bignorm{\tilde p^k}_{H^1_\Gamma}. \]
	Therefore, also $p^k = \tilde p^k + h$ belongs to $H^2_\Gamma(\Omega)$.

	2) We can apply Theorem~5.3.7 in \cite{costabel2010} directly to $q^k$ as the solution of the transmission problem
	\eqref{eq:qk-weak-formulation} with Neumann boundary conditions. Then $q^k$ also lies in $H^2_\Gamma(\Omega)$ and satisfies a similar estimate.
\end{proof}

We can bound the curl-free part of the decomposition by the dissipation and small errors.
Compare Proposition~5.2 in \cite{phung2000} for constant scalar $\epsilon$ and $\mu$.
\begin{lemma}\label{lem:reg_of_pkqk}
	Assume (H).	Let $k \in \{1,2,3\}$ and $\tau \in [0,T_*)$. Then $(\nabla p^k, \nabla q^k)$ belongs to
	$C^1([0,T_*), L^2(\Omega))^6$ with
	\[ \smallnorm{\pt^j \nabla p^k(\tau)}_{L^2(\Omega)}^2,\smallnorm{\pt^j \nabla q^k(\tau)}_{L^2(\Omega)}^2
		\lesssim d(\tau) + z^2(\tau)\qquad \text{for } j\in\{0,1\}. 	\]
\end{lemma}

\begin{proof}
	We only consider $p^k$, since the other case follows similarly. First, we show that
	$\epsilond \nabla p^k + \tilde f_{k-1}$ is differentiable. Let $h \in L^2(\Omega)$. We decompose $h$ into
	\begin{equation} \label{eq:h}
		h = \epsilond^{1/2}\nabla \varphi \oplus \epsilond^{-1/2} \Psi
	\end{equation}
	for $\varphi \in \cH_c$ and $\Psi \in H^{\partial\Omega}(\div 0)$ satisfying $\llbracket \nu^\Gamma \cdot \Psi \rrbracket = 0$
	on $\Gamma$, with orthogonality in $L^2(\Omega)$. To establish this decomposition, we solve the elliptic problem
	\begin{equation} \label{eq:proof:elliptic-pb-for-helmholtz}
		\int_\Omega \epsilond \nabla \varphi \cdot \nabla v \dd x = \int_\Omega \epsilond^{1/2} h \cdot \nabla v \dd x \quad
		\text{for all } v \in \cH_c.
	\end{equation}
	By Lax--Milgram is has a unique solution $\varphi \in \cH_c$ with
	$\smallnorm{\epsilond^{1/2}\nabla \varphi}_{L^2} \lesssim \norm{h}_{L^2}$, analogous to the weak formulation for $p^k$ above.
	Defining $\Psi \coloneqq \epsilond^{1/2} h - \epsilond \nabla \varphi$, we immediately see that
	$\int_\Omega \Psi \cdot \nabla v \dd x = 0$ for all $v \in \cH_c$, which implies the stated divergence and trace properties
	of $\Psi$ as above. Furthermore, the difference satisfies $\smallnorm{\epsilond^{-1/2}\Psi}_{L^2} \lesssim \norm{h}_{L^2}$.
	By setting $v = \varphi$ in \eqref{eq:proof:elliptic-pb-for-helmholtz}, we obtain $\int\Psi \cdot \nabla \varphi \dd x=0$ and
	the orthogonality in \eqref{eq:h} follows. This decomposition leads to
	\begin{align}\label{eq:pk-C1-1}
		\notag\int_\Omega & (\epsilond^{1/2} \nabla p^k + \epsilond^{-1/2} \tilde f_{k-1}) \cdot h \dd x
		= \int_\Omega (\epsilond \nabla p^k + \tilde f_{k-1}) \cdot \nabla \varphi \dd x
		+ \int_\Omega (\epsilond^{-1} \tilde f_{k-1}) \cdot \Psi \dd x                                   \\
		\notag            & = - \int_\Omega \div (\epsilond \nabla p^k + \tilde f_{k-1}) \varphi \dd x
		+ \oint_{\partial\Omega} \nu\cdot (\epsilond \nabla p^k + \tilde f_{k-1}) \varphi \dd \varsigma  \\
		                  & \seq + \oint_\Gamma \llbracket\nu^\Gamma\cdot (\epsilond \nabla p^k
		+ \tilde f_{k-1}) \rrbracket \varphi \dd \varsigma + \int_\Omega (\epsilond^{-1} \tilde f_{k-1}) \cdot \Psi \dd x
		\eqcolon I.
	\end{align}
	Since $\varphi$ is constant, say $c_\ell$, on each $(\partial \Omega)_\ell$ and $c_0=0$,
	equation \eqref{eq:def-pk-BC1} yields
	\[ \oint_{(\partial \Omega)_\ell} \nu\cdot \epsilond \nabla p^k \dd \varsigma
		= \oint_{(\partial \Omega)_\ell} \nu\cdot \epsilond \pt^k E \dd \varsigma, \qquad \ell\in\{1,\dots,K\}. \]
	Due to Stokes' theorem, $\oint_{(\partial \Omega)_\ell} \nu\cdot \curl(\pt^{k-1}H) \dd \varsigma $
	vanishes as $(\partial \Omega)_\ell$ has no boundary in $\partial\Omega$. Hence the Maxwell system \eqref{pde:tilde1} leads to
	\begin{align*}
		\oint_{\partial\Omega} & \nu\cdot (\epsilond \nabla p^k + \tilde f_{k-1}) \varphi \dd \varsigma
		= \sum_{\ell=1}^K c_\ell \oint_{(\partial \Omega)_\ell} \nu\cdot (\epsilond \pt^k E + \tilde f_{k-1}) \dd \varsigma               \\
		                       & = \sum_{\ell=1}^K c_\ell \oint_{(\partial \Omega)_\ell} \nu\cdot (\curl(\pt^{k-1}H) - \sigma \pt^{k-1} E
		-\tilde f_{k-1} + \tilde f_{k-1}) \dd \varsigma                                                                                   \\
		                       & =- \oint_{\partial\Omega} \nu\cdot (\sigma \pt^{k-1} E ) \varphi \dd \varsigma.
	\end{align*}
	Inserting these identities, \eqref{eq:def-pk-eq}, \eqref{eq:def-pk-IC1} and \eqref{eq:normal-jump} into \eqref{eq:pk-C1-1}, we derive
	\begin{align}\label{eq:pk-C1-2}
		\notag I & = \int_\Omega \div (\sigma \pt^{k-1} E) \varphi \dd x
		- \oint_{\partial\Omega} \nu\cdot (\sigma \pt^{k-1} E) \varphi \dd \varsigma
		- \oint_\Gamma \llbracket\nu^\Gamma\cdot (\sigma \pt^{k-1} E) \rrbracket \varphi \dd \varsigma \\
		\notag   & \seq + \int_\Omega (\epsilond^{-1} \tilde f_{k-1}) \cdot \Psi \dd x                 \\
		         & = - \int_\Omega (\sigma \pt^{k-1} E) \cdot \nabla\varphi \dd x
		+ \int_\Omega (\epsilond^{-1} \tilde f_{k-1}) \cdot \Psi \dd x,
	\end{align}
	by means of the divergence theorem on $\upsilon$ and $\omega$. Since the right-hand side is continuously differentiable
	in time, the differences and Cauchy sequences of difference quotients in time of the left-hand side converge in
	$L^2(\Omega)^* \simeq L^2(\Omega)$. (Recall that $\norm{\nabla \varphi}_{L^2},\norm{\Psi}_{L^2} \lesssim \norm{h}_{L^2}$.)
	It then easily follows that $\epsilond^{1/2} \nabla p^k + \epsilond^{-1/2} \tilde f_{k-1}$ belongs to
	$C^1([0,T_*), L^2(\Omega))$ and thus also $\nabla p^k$ by the properties of $\epsilond$ and $\tilde f_{k-1}$.
	The differentiated version of \eqref{eq:pk-C1-2} further yields
	\begin{align*}
		\norm{\pt(\epsilond \nabla p^k)}_{L^2(\Omega)} & \lesssim \norm{\pt(\epsilond^{1/2} \nabla p^k
		+ \epsilond^{-1/2}\tilde f_{k-1})}_{L^2(\Omega)} +\norm{(\pt\epsilond^{1/2}) \epsilond^{1/2}\nabla p^k}_{L^2(\Omega)}                                   \\*
		                                               & \seq + \norm{\pt(\epsilond^{-1/2}\tilde f_{k-1})}_{L^2(\Omega)}                                        \\
		                                               & \lesssim \norm{\sigma \pt^k E}_{L^2(\Omega)} +\norm{\pt (\epsilond^{-1} \tilde f_{k-1})}_{L^2(\Omega)} \\
		                                               & \seq +\norm{(\pt\epsilond^{1/2}) \epsilond^{1/2}\nabla p^k}_{L^2(\Omega)}
		+ \norm{\pt(\epsilond^{-1/2}\tilde f_{k-1})}_{L^2(\Omega)}                                                                                              \\
		                                               & \lesssim d^{1/2}(\tau)+ z(\tau),
	\end{align*}
	where we used \eqref{eq:est-pk} and $\norm{\pt\epsilond^{1/2}}_{H^2_\Gamma} \lesssim z^{1/2}(\tau)$ in the final inequality.
\end{proof}

Next, we subtract $(p^k,q^k)$ from $\pt^k(E,H)$, obtaining the divergence-free part of our Helmholtz decomposition.

\begin{definition}\label{def:VkWk}
	Assume (H).	For $k\in \{1,2,3\}$ we set
	\begin{equation}\label{eq:def-VkWk}
		\begin{alignedat}{2}
			V^k & \coloneqq \pt^k E -\nabla p^k &  & \in C([0,T_*),H(\div_{\epsilond}0)), \\
			W^k & \coloneqq \pt^k H -\nabla q^k &  & \in C([0,T_*),H_{n0}(\div_{\mud}0)).
		\end{alignedat}
	\end{equation}
\end{definition}

\begin{remark}\label{rem:approx-cG4}
	In the following, for $k=3$ several trace terms might not exist in $L^2$. To justify the calculations, we first approximate
	the data and coefficients in the Maxwell system \eqref{pde:pt-maxwell}--\eqref{pde:pt-maxwell_BC} for $k=3$ to obtain solutions
	$(E^3_n,H^3_n) \in \cG_\Gamma^1([0,T_*))$ employing Theorem~\ref{thm:approx}.
	By Lemma~\ref{lem:pkqk-H^2} the corresponding weak solutions $(p_n^3,q_n^3)$ of \eqref{eq:pk-weak-formulation} and
	\eqref{eq:qk-weak-formulation} already lie in $H_\Gamma^2(\Omega)$ for $t \in [0,T_*)$. Thus, also the approximation
	$(V^3_n, W^3_n) \coloneqq (E^3_n,H^3_n) - (\nabla p_n^3, \nabla q_n^3)$ belongs to $H_\Gamma^1(\Omega)$ for $t \in [0,T_*)$.

	The observability estimate in Proposition~\ref{prop:hom_obsI} for $k=3$ only depends on
	$(V^3_n, W^3_n)$ in the $\cG^0$-norm and therefore also holds for the original solution $(V^3,W^3) \in \cG^0$
	after taking the limit $n \to \infty$. Here we again use Theorem~\ref{thm:approx} and the above observations on
	$\nabla(p^k,q^k)$. In the following we drop the index to simplify the notation.
\end{remark}

The properties of $p^k$ and $q^k$ imply the core facts
\begin{equation}\label{eq:properties-Vk}
	\begin{alignedat}{2}
		\curl V^k & = \curl \pt^k E, \qquad & \epsilond V^k & \in H^{\partial\Omega}(\div 0), \\
		\curl W^k & = \curl \pt^k H, \qquad & \mud W^k      & \in H_{n0}(\div 0),
	\end{alignedat}
\end{equation}
Here the condition $\langle\tr_n\epsilond V^k,\1\rangle_{H^{-\frac12}((\partial\Omega)_i)}=0$
is clear by \eqref{eq:def-pk-BC1} for $i \neq 0$. For $i=0$ it then follows from $\div(\epsilond V^k)=0$ which yields
$\langle\tr_n\epsilond V^k,\1\rangle_{H^{-\frac12}(\partial\Omega)}=0$. Moreover, $(p^k,q^k)$ were constructed such that
\begin{equation}\label{eq:interface-cond-VW}
	\begin{alignedat}{3}
		\tr_{t,\partial\Omega} V^k       & =\llbracket V^k \times \nu^\Gamma \rrbracket          &
		                                 & =\llbracket W^k \times \nu^\Gamma \rrbracket          &   & = 0, \\
		\tr_{n,\partial\Omega}(\mud W^k) & =\llbracket \epsilond V^k \cdot \nu^\Gamma \rrbracket &
		                                 & =\llbracket \mud W^k \cdot \nu^\Gamma \rrbracket      &   & = 0,
	\end{alignedat}
\end{equation}
in $H^{1/2}$ preserving \eqref{pde:pt-maxwell-IC}, see also Remark~\ref{rem:approx-cG4} for $k=3$.

We state two results needed in the observability estimate.
Lemmas~\ref{lem:energy-decreasing} and \ref{lem:reg_of_pkqk} yield an energy estimate for $(V^k,W^k)$.
\begin{lemma}\label{lem:VkWk-est}
	Assume (H).		For $0 \leq s \leq t < T_*$ we have
	\begin{align*}\label{lem:energy-est-VkWk}
		\norm{(V^k,W^k)(t)}_{L^2(\Omega)}^2 & \lesssim e(t) + \norm{(\nabla p^k, \nabla q^k)(t)}_{L^2(\Omega)}^2
		\lesssim e(t) + z^2(t)                                                                                   \\
		                                    & \lesssim e(s) + \int_s^t z^{3/2}(\tau) \dd \tau + z^2(t).
	\end{align*}
\end{lemma}

Next, we compare $(V^k,W^k)$ to $\pt(V^{k-1}, W^{k-1})$ for $k\geq1$.
\begin{lemma}\label{lem:diff-ptVWk-1-VWk}
	Assume (H).		Let $k \in \{1,2,3\}$. Then for $0\leq \tau < T_*$
	\begin{alignat*}{-1}
		\norm{(V^k - \pt V^{k-1})(\tau)}_{L^2(\Omega)} & \lesssim \norm{((\pt\epsilond)\nabla p^{k-1})(\tau)}_{L^2(\Omega)} &  & \lesssim z^{3/2}(\tau), \\
		\norm{(W^k - \pt W^{k-1})(\tau)}_{L^2(\Omega)} & \lesssim \norm{((\pt\mud)\nabla q^{k-1})(\tau)}_{L^2(\Omega)}      &  & \lesssim z^{3/2}(\tau).
	\end{alignat*}
\end{lemma}

\begin{proof} We only consider $V^k$, since for $W^k$ the argument works in a similar way.
	Note that $V^k - \pt V^{k-1} = \nabla p^k - \pt \nabla p^{k-1}$ solves
	\[ \div_{\epsilond}(\nabla p^k - \pt \nabla p^{k-1})
		=\div_{\epsilond} \pt^k E - \div_{\epsilond} \pt^k E - \div_{\pt\epsilond} \nabla p^{k-1}. \]
	As in \eqref{eq:est-pk} we then deduce
	\begin{align*}
		\norm{(V^k - \pt V^{k-1})(\tau)}_{L^2(\Omega)} & \lesssim \norm{((\pt\epsilond)\nabla p^{k-1})(\tau)}_{L^2(\Omega)}
		\lesssim z^{3/2}(\tau). \qedhere
	\end{align*}
\end{proof}

\section{Observability-type estimate for \texorpdfstring{$\pt (V^k, W^k)$}{d/dt (V\^{}k,W\^{}k)}}
\label{sec:obs}

The observability-type estimate in Proposition~\ref{prop:hom_obsI} only works for charge-free fields. Hence, we first treat
the fields $(V_k,W_k)$ from Definition~\ref{def:VkWk}. Lemma~\ref{lem:reg_of_pkqk} then allows us to go back to $\pt^k(E,H)$
in Theorem~\ref{thm:obs-est}. So, we insert the Helmholtz decomposition $(V^k,W^k) = (\pt^k E-\nabla p^k, \pt^k H-\nabla q^k)$
from \eqref{eq:def-VkWk} into Maxwell equations~\eqref{pde:pt-maxwell}. Together with \eqref{eq:properties-Vk}, this results
in the system
\begin{align}
	\begin{split}\notag
		\pt(\epsilond V^k) & = \curl W^k - \sigma \pt^k E - \pt (\epsilond\nabla p^k) - \pt f_k, \\
		\pt(\mud W^k)      & = - \curl V^k- \pt (\mud \nabla q^k )- \pt g_k,
	\end{split}\  &                                                                         & t\geq 0, x \in \Omega_\pm,                                                                                                                \\
	\notag 
	\llbracket V^k \times \nu^\Gamma \rrbracket                                                 & =\llbracket W^k \times \nu^\Gamma \rrbracket = 0, \quad
	\llbracket \epsilond V^k \cdot \nu^\Gamma \rrbracket =\llbracket \mud W^k \cdot \nu^\Gamma \rrbracket = 0,
	                                                                                            &                                                                         & t\geq 0, x \in \Gamma,                                      \\
	\notag V^k \times \nu                                                                       & = 0,                                                                    &                            & t\geq 0, x \in \partial\Omega, \\
	V^k(0)                                                                                      & = \pt^k E(0) - \nabla p^k(0), \quad W^k(0) = \pt^k H(0)- \nabla q^k(0), &                            & x \in \Omega,\label{eq:VW}
\end{align}
for $k \in \{1,2,3\}$. (The case $k=0$ is handled by means of Lemma~\ref{lem:cE-estimate}.)

Our goal is to control space-time $L^2$-norms of $V^k$ and $W^k$ for $k>0$ by the initial data, dissipation terms supported on $\cN_a$ and
various small error terms. Our first estimate still includes the $L^2$-norm of $W^k$ on the right, which is not damped
through Maxwell. This term is treated afterwards. (One could establish a similar estimate for $k=0$, however here one lacks control over
the magnetic term on the right-hand side.) For $k=3$ we approximate $(V^k,W^k)$ as in Remark~\ref{rem:approx-cG4}, so that the trace
terms appearing in the subsequent proof are well-defined.

\begin{proposition}[Observability-type Estimate I]\label{prop:hom_obsI}
	Assume that \ref{hypothesis-quasi} holds. For $k \in \{1,2,3\}$ and $0\leq s\leq t < T_*$, we have
	\begin{align*}
		\int_s^t \!\! \int_{\upsilon \cup \cN_b^\Gamma} \left(\abs{V^k}^2 + \abs{W^k}^2\right) \dd x \dd \tau & \lesssim e(s) + z^2(s) + z^2(t) + \int_s^t \left(d(\tau) +z^{3/2}(\tau)\right) \dd \tau \\
		                                                                                                      & \seq +\int_s^t\!\!\int_{{\omega \setminus \cN_b^\Gamma}}\abs{W^k}^2 \dd x \dd \tau.
	\end{align*}
\end{proposition}

\begin{proof}
	\newcommand{\deltaA}{\delta_{1}}
	\newcommand{\deltaB}{\delta_{2}}
	\newcommand{\deltaN}{\delta_{3}}
	\newcommand{\deltaO}{\delta_{4}}
	\newcommand{\deltaP}{\delta_{5}}
	We first give a brief overview of the structure of the proof, using terms only defined later.
	Our reasoning is an involved energy-type estimate using a \textit{Morawetz multiplier}. Starting from \eqref{eq:VW},
	we obtain the identity
	\begin{equation}\label{eq:morawetz}
		0 = \int_s^t \!\! \int_\Omega \left(\pt( \epsilond V^k) - \curl W^k + \sigma \pt^k E + \pt (\epsilond\nabla p^k)
		+ \pt f_k\right)\cdot((\mud W^k) \times \chi m) \dd x \dd \tau,
	\end{equation}
	where $k\in \{1,2,3\}$, $m(x) = x-x_0$ for $x_0\in \R^3$ from \eqref{eq:non-trapping}, and
	$0\leq \chi \in C^1(\overline\Omega)$ is a cutoff to be varied in the calculations.
	Integrating by parts and several lengthy calculations lead to the desired observability estimate.
	We freely use standard formulas from vector analysis.
	\begin{itemize}
		\item[1)] We start with the Morawetz multiplier in Equation~\eqref{eq:morawetz} above.
		\item[2)] After integrating by parts (see Equation~(20) of Chapter~2.1 in \cite{cessenat1996} or Equation~(3.1)
		      in \cite{ciarlet18}), Equation~\eqref{eq:morawetz} leads to several terms we deal with individually. Namely,
		      \begin{itemize}
			      \item[2.1)] $I_t$ -- a term, where the time derivative falls upon $\curl V^k$ and
			      \item[2.2)] $I_c$ -- a term where $\curl$ is applied to $\chi m \times \curl V^k$.
			      \item[2.3)] We combine the previous two steps to one estimate.
			      \item[2.4)] Letting $\chi = 1$ on $\Gamma$ we show that, due to the non-trapping condition on the interface~\eqref{eq:non-trapping-interface},
			            the interface terms have a negative sign and can be neglected.
			      \item[2.5)] We estimate several error terms.
			      \item[2.6)] For $\chi \equiv 1$ we obtain our desired observability estimate up to boundary terms, which we treat in step
			      \item[2.7)] to finally complete the estimate using a cutoff with $\chi = 1$ on $\partial\Omega$.
		      \end{itemize}
	\end{itemize}

	\smallskip


	1) Starting from Equation~\eqref{eq:morawetz} and using \eqref{eq:VW} for $\pt(\mud W^k)$, we calculate
	\begin{align*}
		0 & =\! \int_s^t \!\! \int_\Omega\! \big(\pt( \epsilond V^k) - \curl W^k + \sigma \pt^k E + \pt (\epsilond\nabla p^k)
		+ \pt f_k\big)\cdot((\mud W^k)\! \times\! \chi m) \dd x \dd \tau                                                                               \\
		  & = \int_s^t \!\! \int_\Omega \pt \left(\epsilond V^k \cdot((\mud W^k) \times \chi m)\right) \dd x \dd \tau
		+ \int_s^t \!\! \int_\Omega \left(\epsilond V^k \cdot(\curl V^k \times \chi m)\right) \dd x \dd \tau                                           \\
		  & \seq + \int_s^t \!\! \int_\Omega \left(\epsilond V^k \cdot(\pt(\mud\nabla q^k) \times \chi m)\right) \dd x \dd \tau
		+ \int_s^t \!\! \int_\Omega \left(\epsilond V^k \cdot(\pt g_k \times \chi m)\right) \dd x \dd \tau                                             \\
		  & \seq - \int_s^t \!\! \int_\Omega \curl W^k\cdot((\mud W^k) \times \chi m) \dd x \dd \tau
		+ \int_s^t \!\! \int_\Omega \sigma \pt^k E\cdot((\mud W^k) \times \chi m) \dd x \dd \tau                                                       \\
		  & \seq + \int_s^t \!\! \int_\Omega \pt (\epsilond\nabla p^k)\cdot((\mud W^k) \times \chi m) \dd x \dd \tau
		+ \int_s^t \!\! \int_\Omega \pt f_k\cdot((\mud W^k) \times \chi m) \dd x \dd \tau                                                              \\
		  & \eqqcolon \left[\int_\Omega \epsilond V^k \cdot((\mud W^k) \times \chi m) \dd x\right]_s^t + I_t + I_q + I_g + I_c + I_\sigma + I_p + I_f.
	\end{align*}

	2) We consider the terms above individually.


	2.1) By the interface and boundary conditions~\eqref{eq:interface-cond-VW}, the field $V^k \times \nu$ vanishes at $\partial \Omega$ and
	$V^k \times \nu$ does not jump at the interface $\Gamma$. We integrate by parts on $\upsilon$ and $\omega$ separately. Then $I_t$ becomes
	\begin{align*}
		I_t & = \int_s^t \!\! \int_\Omega \curl V^k\cdot(\chi m \times (\epsilond V^k)) \dd x \dd \tau                                                        \\
		    & = \int_s^t \!\! \int_\Omega \!V^k\cdot\curl\!\left(\chi m\! \times\! (\epsilond V^k)\right) \dd x \dd \tau
		+ \int_s^t \!\! \oint_\Gamma \left\llbracket\left(\nu^\Gamma \!\times\! V^k\right) \cdot\left(\chi m \times (\epsilond V^k)\right)
		\right\rrbracket \dd \varsigma\dd \tau                                                                                                                \\
		    & = \int_s^t \!\! \int_\Omega \! V^k\cdot\curl\left(\chi m\! \times \! (\epsilond V^k)\right) \dd x \dd \tau
		+ \int_s^t \!\! \oint_\Gamma\! \left\llbracket\chi m \cdot\left((\epsilond V^k)\times \left(\nu^\Gamma \!\times \!V^k\right)\right)
		\right\rrbracket \dd \varsigma \dd \tau                                                                                                               \\
		    & = \int_s^t \!\! \int_\Omega V^k\cdot\curl\left(\chi m \times (\epsilond V^k)\right) \dd x \dd \tau                                              \\*
		    & \seq + \int_s^t\!\!\oint_\Gamma\left\llbracket\chi (m \cdot\nu^\Gamma)(\epsilond V^k\cdot V^k)-(\chi m \cdot V^k)(\epsilond V^k\cdot\nu^\Gamma)
		\right\rrbracket \dd \varsigma \dd \tau.
	\end{align*}
	Using the curl vector product rule
	\begin{equation}\label{eq:curl-product}
		\curl(a \times b) = (b\cdot\nabla) a - (a \cdot\nabla) b + a \div b - b \div a ,
	\end{equation}
	and $\div(\epsilond V^k)=0$, we obtain the pointwise identity
	\begin{align*}
		 & \curl\left(\chi m \times (\epsilond V^k)\right)
		= \left((\epsilond V^k) \cdot\nabla\right)(\chi m) - \chi (m \cdot \nabla)(\epsilond V^k) - (\epsilond V^k) \div(\chi m) \\
		 & \ = \left((\epsilond V^k)\cdot \nabla \chi\right)m + \chi(\epsilond V^k) - \chi (m \cdot \nabla)(\epsilond V^k)
		- (\nabla \chi \cdot m)(\epsilond V^k) - 3\chi(\epsilond V^k)                                                            \\
		 & \ = \left((\epsilond V^k)\cdot \nabla \chi\right)m -2 \chi(\epsilond V^k) - \chi (m \cdot \nabla)(\epsilond V^k)
		- (\nabla \chi \cdot m)(\epsilond V^k)
	\end{align*}
	almost everywhere on $\upsilon \cup \omega$. It follows that
	\begin{align*}
		I_t & = \int_s^t \!\! \int_\Omega\Big( ( V^k \cdot m )((\epsilond V^k) \cdot \nabla \chi) - 2 \chi \epsilond V^k\cdot V^k                              \\
		    & \seq - (\nabla \chi \cdot m) (\epsilond V^k \cdot V^k) - \left(\chi (m\cdot \nabla)(\epsilond V^k)\right)\cdot V^k \Big)\dd x \dd \tau           \\
		    & \seq + \int_s^t \!\! \oint_\Gamma \llbracket \chi (m \cdot\nu^\Gamma)(\epsilond V^k \cdot V^k) -(\chi m \cdot V^k)(\epsilond V^k\cdot\nu^\Gamma)
		\rrbracket \dd \varsigma \dd \tau .
	\end{align*}
	Notice that
	\begin{equation*}
		\left(\chi (m\cdot \nabla)(\epsilond V^k)\right)\cdot V^k = \frac12 \chi (m\cdot \nabla)\left(\epsilond V^k\cdot V^k\right)
		+\frac12 \left(\chi (m\cdot \nabla)\epsilond\right) V^k\cdot V^k.
	\end{equation*}
	The divergence theorem on $\upsilon \cup \omega$ further yields
	\begin{align*}
		 & - \frac12 \int_s^t \!\! \int_\Omega\chi (m\cdot \nabla)\left(\epsilond V^k \cdot V^k\right)\dd x \dd \tau                               \\
		 & \quad= \int_s^t \!\! \int_\Omega\Big(\frac12 \nabla \chi \cdot m (\epsilond V^k \cdot V^k) + \frac32 \chi(\epsilond V^k \cdot V^k)\Big)
		\dd x \dd \tau                                                                                                                             \\
		 & \quad\seq -\int_s^t \!\! \oint_\Gamma \frac12 (m\cdot\nu^\Gamma) \llbracket\chi \epsilond V^k \cdot V^k\rrbracket\dd \varsigma \dd \tau
		-\int_s^t \!\! \oint_{\partial \Omega} \frac12 \chi (m \cdot \nu) (\epsilond V^k \cdot V^k)\dd \varsigma \dd \tau.
	\end{align*}

	Involving $\tildeepsilond = \epsilond + (m\cdot \nabla) \epsilond$, which is positive definite on $\upsilon \cup \cN_b^\Gamma$ by
	\eqref{eq:def-tilde-epsilond}, we thus derive
	\begin{align*}
		I_t & = \int_s^t \!\! \int_\Omega\Big( ( V^k \cdot m )((\epsilond V^k) \cdot \nabla \chi) - 2 \chi \epsilond V^k\cdot V^k
		- (\nabla \chi \cdot m) (\epsilond V^k \cdot V^k)                                                                               \\
		    & \seq + \frac12 (\nabla \chi \cdot m) (\epsilond V^k \cdot V^k) + \frac32 \chi(\epsilond V^k \cdot V^k)
		- \frac12 \left(\chi (m\cdot \nabla)\epsilond\right) V^k\cdot V^k \Big)\dd x \dd \tau                                           \\
		    & \seq + \int_s^t \!\! \oint_\Gamma \left\llbracket \frac12 \chi (m \cdot\nu^\Gamma)(\epsilond V^k \cdot V^k)
		-(\chi m \cdot V^k)(\epsilond V^k\cdot\nu^\Gamma) \right\rrbracket \dd \varsigma \dd \tau                                       \\
		    & \seq -\int_s^t \!\! \oint_{\partial \Omega} \frac12 \chi (m \cdot \nu) (\epsilond V^k \cdot V^k) \dd \varsigma \dd \tau   \\
		    & = \int_s^t \!\! \int_\Omega\Big( ( V^k \cdot m )((\epsilond V^k) \cdot \nabla \chi) - \frac12 \chi \epsilond V^k\cdot V^k
		- \frac 12(\nabla \chi \cdot m) (\epsilond V^k \cdot V^k)                                                                       \\*
		    & \seq - \frac12 \left(\chi (m\cdot \nabla)\epsilond\right) V^k\cdot V^k \Big)\dd x \dd \tau
		-\int_s^t \!\! \oint_{\partial \Omega} \frac12 \chi (m \cdot \nu) (\epsilond V^k \cdot V^k) \dd \varsigma \dd \tau              \\*
		    & \seq + \int_s^t \!\! \oint_\Gamma \left\llbracket \frac12 \chi (m \cdot\nu^\Gamma)(\epsilond V^k \cdot V^k)
		-(\chi m \cdot V^k)(\epsilond V^k\cdot\nu^\Gamma) \right\rrbracket \dd \varsigma \dd \tau                                       \\*
		    & = \int_s^t \!\! \int_\Omega\Big( ( V^k \cdot m )((\epsilond V^k) \cdot \nabla \chi)
		- \frac 12(\nabla \chi \cdot m) (\epsilond V^k \cdot V^k) \Big)\dd x \dd \tau                                                   \\
		    & \seq + \int_s^t \!\! \oint_\Gamma (m \cdot\nu^\Gamma)\left\llbracket \frac12 \chi (\epsilond V^k \cdot V^k)
		-(\chi m \cdot V^k)(\epsilond V^k\cdot\nu^\Gamma) \right\rrbracket \dd \varsigma \dd \tau                                       \\
		    & \seq-\int_s^t \!\! \oint_{\partial \Omega} \frac12 \chi (m \cdot \nu) (\epsilond V^k \cdot V^k) \dd \varsigma \dd \tau    \\
		    & \seq - \int_s^t \!\! \int_{\upsilon \cup \cN_b^\Gamma} \frac12 \chi \tildeepsilond V^k\cdot V^k\dd x \dd \tau
		- \int_s^t \!\! \int_{\omega \setminus \cN_b^\Gamma} \frac12 \chi \tildeepsilond V^k\cdot V^k\dd x \dd \tau.
	\end{align*}
	We denote the interface term $I_t^\Gamma \coloneqq I_{t,1}^\Gamma +I_{t,2}^\Gamma$ and the boundary term by
	\begin{align}
		I_{t,1}^\Gamma       & \coloneqq \int_s^t \!\! \oint_\Gamma (m \cdot\nu^\Gamma)\left\llbracket \frac12 \chi (\epsilond V^k \cdot V^k)\right\rrbracket
		\dd \varsigma \dd \tau, \notag                                                                                                                        \\
		I_{t,2}^\Gamma       & \coloneqq \int_s^t \!\! \oint_\Gamma \llbracket-(\chi m \cdot V^k)(\epsilond V^k\cdot\nu^\Gamma) \rrbracket \dd \varsigma
		\dd \tau, \label{def:ItG}                                                                                                                             \\
		I_t^{\partial\Omega} & \coloneqq-\int_s^t\!\!\oint_{\partial \Omega} \frac12\chi(m\cdot\nu)(\epsilond V^k\cdot V^k)\dd \varsigma\dd\tau.\notag
	\end{align}

	2.2) We now turn our attention towards $I_c$. We set
	\begin{align*}
		I_{c,1}^{\partial \Omega} & \coloneqq \int_s^t \!\! \oint_{\partial \Omega} W^k \cdot\left( \nu \times ((\mud W^k) \times \chi m)\right)
		\dd \varsigma \dd \tau,                                                                                                                              \\
		I_{c,1}^{\Gamma}          & \coloneqq \int_s^t \!\! \oint_\Gamma \left\llbracket W^k \cdot\left( \nu^\Gamma \times ((\mud W^k) \times \chi m)\right)
		\right\rrbracket \dd \varsigma \dd \tau.
	\end{align*}
	By means of the curl product rule~\eqref{eq:curl-product} and $\div(\mud W^k)=0$, we derive
	\begin{align*}
		I_c & = \int_s^t \!\! \int_\Omega W^k\cdot \curl(\chi m \times (\mud W^k)) \dd x \dd \tau + I_{c,1}^{\partial \Omega} + I_{c,1}^{\Gamma}           \\
		    & = \int_s^t \!\! \int_\Omega W^k\cdot \Big(\chi m \div(\mud W^k) - \div (\chi m)(\mud W^k)                                                    \\*
		    & \seq + ((\mud W^k)\cdot \nabla)(\chi m) - (\chi m \cdot \nabla) (\mud W^k)\Big) \dd x \dd \tau + I_{c,1}^{\partial \Omega}+ I_{c,1}^{\Gamma} \\
		    & = \int_s^t \!\! \int_\Omega W^k\cdot \Big(- \div (\chi m)(\mud W^k)                                                                          \\*
		    & \seq + ((\mud W^k)\cdot \nabla)(\chi m) - (\chi m \cdot \nabla\mud) W^k- \mud(\chi m \cdot \nabla) W^k\Big) \dd x \dd \tau
		+ I_{c,1}^{\partial \Omega} + I_{c,1}^{\Gamma}.
	\end{align*}
	The last term on $\Omega$ is treated by the divergence theorem on $\upsilon \cup \omega$, resulting in
	\begin{align*}
		-\! & \int_s^t \!\! \int_\Omega \mud W^k\cdot ((\chi m \cdot \nabla) W^k) \dd x \dd \tau
		= - \sum_{j=1}^3\int_s^t \!\! \int_\Omega \chi (\mud W^k)_j m \cdot \nabla ( W^k)_j \dd x \dd \tau                                          \\
		    & = \sum_{j=1}^3\int_s^t \!\! \int_\Omega \div\left(\chi (\mud W^k)_j m\right) ( W^k)_j \dd x \dd \tau                                  \\
		    & \seq- \int_s^t \!\! \oint_\Gamma (m\cdot\nu^\Gamma)\llbracket\chi\mud W^k\cdot W^k\rrbracket \dd \varsigma \dd \tau
		- \int_s^t \!\! \oint_{\partial\Omega}\chi (m\cdot\nu)(\mud W^k\cdot W^k) \dd \varsigma \dd \tau                                            \\
		    & \ = \int_s^t \!\! \int_\Omega \Big[\div(\chi m) \mud W^k\cdot W^k + (\chi m \cdot \nabla) (\mud W^k )\cdot ( W^k) \Big]\dd x \dd \tau \\
		    & \seq - \int_s^t \!\! \oint_\Gamma (m\cdot\nu^\Gamma)\llbracket\chi\mud W^k\cdot W^k\rrbracket \dd \varsigma \dd \tau
		- \int_s^t \!\! \oint_{\partial\Omega}\chi (m\cdot\nu)(\mud W^k\cdot W^k) \dd \varsigma \dd \tau                                            \\
		    & = \!\int_s^t \!\! \int_\Omega \!\Big[\div(\chi m) \mud W^k\cdot W^k + (\chi m \!\cdot\! \nabla\mud ) W^k\cdot W^k
		+ \mud W^k \cdot (\chi m \!\cdot\! \nabla) W^k\Big]\mathrm{d} x \dd \tau                                                                    \\
		    & \seq- \int_s^t \!\! \oint_\Gamma (m\cdot\nu^\Gamma)\llbracket\chi\mud W^k\cdot W^k\rrbracket \dd \varsigma \dd \tau
		- \int_s^t \!\! \oint_{\partial\Omega}\chi (m\cdot\nu)(\mud W^k\cdot W^k) \dd \varsigma \dd \tau \eqcolon S.
	\end{align*}
	Bringing the term $\int_s^t \!\! \int_\Omega \mud W^k\cdot ((\chi m \cdot \nabla) W^k) \dd x \dd \tau$ to the left-hand side, we infer
	\begin{align*}
		S & = \frac12 \int_s^t \!\! \int_\Omega \Big[\div(\chi m) \mud W^k\cdot W^k + (\chi m \cdot \nabla\mud ) W^k\cdot W^k\Big]\dd x \dd \tau \\
		  & \seq - \frac12\int_s^t \!\! \oint_\Gamma (m\cdot\nu^\Gamma)\llbracket\chi\mud W^k\cdot W^k\rrbracket \dd \varsigma \dd \tau
		- \frac12\int_s^t \!\! \oint_{\partial\Omega}\chi (m\cdot\nu)(\mud W^k\cdot W^k) \dd \varsigma \dd \tau .
	\end{align*}

	With $I_1\coloneq I_{c,1}^{\partial \Omega} + I_{c,1}^{\Gamma}$, these identities lead to
	\begin{align*}
		I_c & =\! \int_s^t \!\! \int_\Omega\! W^k\!\cdot\! \Big[{-}\frac12 \div (\chi m)\mud W^k
		+ ((\mud W^k)\cdot \nabla)(\chi m) - \frac12(\chi m \cdot \nabla\mud) W^k\Big]\mathrm{d} x \dd \tau                                                     \\*
		    & \seq - \frac12\!\int_s^t \!\! \oint_\Gamma\! (m\!\cdot\!\nu^\Gamma)\llbracket\chi\mud W^k\!\cdot\! W^k\rrbracket \dd \varsigma \dd \tau
		- \frac12\!\int_s^t \!\! \oint_{\partial\Omega}\!\chi (m\!\cdot\!\nu)(\mud W^k\!\cdot\! W^k) \dd \varsigma \dd \tau + I_1                               \\
		    & = \int_s^t \!\! \int_\Omega\! W^k\!\cdot\! \Big[ {-} \frac12 (m\cdot\nabla \chi + 3\chi)(\mud W^k) + \mud W^k \nabla \chi \cdot m + \mud W^k \chi \\
		    & \seq\seq - \frac12(\chi m \cdot \nabla\mud) W^k\Big]\mathrm{d} x \dd \tau
		- \frac12\int_s^t \!\! \oint_\Gamma (m\cdot\nu^\Gamma)\llbracket\chi\mud W^k\cdot W^k\rrbracket \dd \varsigma \dd \tau                                  \\
		    & \seq - \frac12\int_s^t \!\! \oint_{\partial\Omega}\chi (m\cdot\nu)(\mud W^k\cdot W^k) \dd \varsigma \dd \tau +I_1                                 \\
		    & = \int_s^t \!\! \int_\Omega W^k\cdot \Big(- \frac12 (m \cdot\nabla \chi)(\mud W^k) + (\mud W^k \cdot \nabla\chi)m \Big)\dd x \dd \tau +I_1        \\*
		    & \seq - \frac12\int_s^t \!\! \oint_\Gamma (m\cdot\nu^\Gamma)\llbracket\chi\mud W^k \cdot W^k\rrbracket \dd \varsigma \dd \tau
		- \frac12\int_s^t \!\! \oint_{\partial\Omega}\chi (m\cdot\nu)(\mud W^k \cdot W^k) \dd \varsigma \dd \tau                                                \\*
		    & \seq - \int_s^t \!\! \int_{\upsilon \cup \cN_b^\Gamma} \frac12 \chi \tildemud W^k \cdot W^k\dd x \dd \tau
		- \int_s^t \!\! \int_{\omega \setminus \cN_b^\Gamma} \frac12 \chi \tildemud W^k \cdot W^k\dd x \dd \tau,
	\end{align*}
	where ${\tildemud} \coloneqq \mud + m \cdot \nabla \mud \gtrsim \eta \mud$, see \ Equation~\eqref{eq:def-tilde-epsilond}.
	The interface and boundary terms above are denoted by $I_{c,2}^\Gamma$ and $I_{c,2}^{\partial \Omega}$, and we abbreviate
	\pagebreak[4]
	\begin{align}\label{eq:def-I_c}
		I_c^\Gamma \coloneqq I_{c,1}^\Gamma + I_{c,2}^\Gamma
		 & =\int_s^t \!\!\oint_\Gamma \left\llbracket W^k\cdot\left( \nu^\Gamma \times ((\mud W^k)\times\chi m)\right)\right\rrbracket
		\dd\varsigma\dd\tau \notag                                                                                                       \\*
		 & \seq - \frac12\int_s^t \!\! \oint_\Gamma \chi(m\cdot\nu^\Gamma)\llbracket\mud W^k \cdot W^k\rrbracket \dd \varsigma \dd \tau, \\
		I_c^{\partial \Omega} \coloneqq I_{c,1}^{\partial \Omega} + I_{c,2}^{\partial \Omega}\notag
		 & = \int_s^t \!\! \oint_{\partial \Omega} W^k \cdot\left( \nu \times ((\mud W^k) \times \chi m)\right) \dd \varsigma \dd \tau   \\*
		 & \seq - \frac12\int_s^t \!\! \oint_{\partial\Omega}\chi (m\cdot\nu)(\mud W^k \cdot W^k) \dd \varsigma \dd \tau.\notag
	\end{align}
	(If we had $m = \nu^\Gamma$ the term in $I_{c,1}^\Gamma$ would correspond to the tangential projection.)

	\smallskip

	2.3) The previous two steps yield
	\begin{align*}
		\Big[\int_\Omega & (\epsilond V^k) (\chi m \times - W^k) \Big]_s^t = -I_t - I_c - I_\sigma - I_p - I_f - I_q - I_g                                       \\
		                 & = \int_s^t \!\! \int_{(\upsilon \cup \cN_b^\Gamma)\cupdot (\omega \setminus \cN_b^\Gamma)} \frac \chi 2 (\tildeepsilond V^k \cdot V^k
		+ \tildemud W^k \cdot W^k) \dd x\dd \tau                                                                                                                 \\*
		                 & \seq + \int_s^t \!\! \int_\Omega\Big(\frac12 (\nabla \chi \cdot m)(\epsilond V^k \cdot V^k + \mud W^k \cdot W^k)                      \\*
		                 & \seq\seq\hphantom{\int_s^t \!\! \int_\Omega} - ( V^k \cdot m)( (\epsilond V^k) \cdot \nabla \chi)
		- (m \cdot W^k)(\mud W^k \cdot \nabla \chi) \Big)\dd x \dd \tau                                                                                          \\*
		                 & \seq - I_t^\Gamma - I_t^{\partial \Omega} - I_c^\Gamma - I_c^{\partial \Omega} - I_\sigma - I_p - I_f - I_q - I_g.
	\end{align*}
	We rearrange the equation into
	\begin{align}\label{eq:proof-non-est-eq}
		\int_s^t \!\! & \int_{\upsilon \cup \cN_b^\Gamma} \frac \chi 2 (\tildeepsilond V^k \cdot V^k + \tildemud W^k \cdot W^k) \dd x\dd \tau           \\
		              & = \left[\int_\Omega (\epsilond V^k) (\chi m \times - W^k) \right]_s^t
		- \int_s^t \!\! \int_{\omega \setminus \cN_b^\Gamma} \frac \chi 2 (\tildeepsilond V^k \cdot V^k + \tildemud W^k \cdot W^k) \dd x\dd \tau \notag \\
		              & \seq - \int_s^t \!\! \int_\Omega\Big(\frac12 (\nabla \chi \cdot m)(\epsilond V^k \cdot V^k + \mud W^k \cdot W^k) \notag         \\
		              & \seq\seq\hphantom{\int_s^t \!\! \int_\Omega\Big(} - ( V^k \cdot m)( (\epsilond V^k) \cdot \nabla \chi)
		- (m \cdot W^k)(\mud W^k \cdot \nabla \chi) \Big)\dd x \dd \tau \notag                                                                          \\
		              & \seq + I_t^\Gamma + I_t^{\partial \Omega} + I_c^\Gamma + I_c^{\partial \Omega} + I_\sigma + I_p + I_f + I_q + I_g. \notag
	\end{align}

	2.4) The non-trapping condition at the interface $(m \cdot \nu^\Gamma)\llbracket \epsilond \rrbracket \leq 0$ yields
	\begin{align*}
		 & \llbracket \epsilond \rrbracket >0 \text{ and } (m \cdot \nu^\Gamma) \leq 0, \quad
		\text{or } \ \llbracket \epsilond \rrbracket = 0, \quad \text{or } \
		\llbracket \epsilond \rrbracket <0 \text{ and } (m \cdot \nu^\Gamma) \geq 0,
	\end{align*}
	on $\Gamma$, see \eqref{eq:non-trapping-interface}. We next show that the interface integrals $I_{t}^\Gamma$ and
	$I_{c}^\Gamma$ from \eqref{def:ItG} resp.\ \eqref{eq:def-I_c} are non-positive, so that they can be dropped
	when estimating \eqref{eq:proof-non-est-eq}.

	We take $\chi$ with $\chi = 1$ on the interface $\Gamma$. The integrand in $I_{t}^\Gamma$ then becomes
	\begin{equation}\label{eq:proof-integrand-I_t^Gamma}
		\tfrac{1}{2} (m\cdot\nu^\Gamma) \llbracket V^k \cdot \epsilond V^k \rrbracket - \llbracket (m \cdot V^k) (\nu^\Gamma\cdot\epsilond V^k)\rrbracket.
	\end{equation}
	As $\llbracket \tr_\tau V^k \rrbracket =0$ by \eqref{eq:interface-cond-VW}, the jump across $\Gamma$ only occurs in normal direction, i.e.,
	\begin{equation*}
		V^k|_{\upsilon} = V^k|_{\omega} + \alpha \nu^\Gamma
	\end{equation*}
	for some $\alpha(x) \in \R$. Furthermore, using also the continuity $\llbracket \nu^\Gamma \cdot \epsilond V^k \rrbracket = 0$ of the jump of the
	`displacement' field from \eqref{eq:interface-cond-VW}, we can determine $\alpha$ as
	\[ \nu^\Gamma \cdot \restr{\epsilond}{\upsilon} (V^k|_{\omega} + \alpha \nu^\Gamma) = \nu^\Gamma \cdot \restr{\epsilond}{\omega} V^k|_{\omega}
		\quad\implies\quad \alpha = - \frac{\nu^\Gamma\cdot \llbracket\epsilond\rrbracket V^k|_{\omega}}{\nu^\Gamma\cdot\epsilond|_{\upsilon}\nu^\Gamma}.\]
	Next, we insert this decomposition into \eqref{eq:proof-integrand-I_t^Gamma}. The first term is given by
	\begin{align*}
		\llbracket V^k\! \cdot\! \epsilond V^k \rrbracket
		 & = (V^k|_{\omega} + \alpha \nu^\Gamma) \cdot \restr{\epsilond}{\upsilon} (V^k|_{\omega} + \alpha \nu^\Gamma)
		- \restr{\left(V^k \cdot \epsilond V^k\right)\!}{\omega}                                                                                     \\
		 & = V^k|_{\omega} \cdot \restr{\epsilond}{\upsilon} V^k|_{\omega} + 2\alpha \nu^\Gamma \cdot \restr{\epsilond}{\upsilon} V^k|_{\omega}
		+ \alpha^2 \nu^\Gamma \cdot \restr{\epsilond}{\upsilon} \nu^\Gamma - \restr{\left(V^k \cdot \epsilond V^k\right)\!}{\omega}                  \\
		 & = V^k|_{\omega} \cdot \llbracket \epsilond\rrbracket V^k|_{\omega} + 2\alpha (\nu^\Gamma \cdot \restr{\epsilond}{\upsilon} V^k|_{\omega})
		+ \alpha^2 (\nu^\Gamma \cdot \restr{\epsilond}{\upsilon} \nu^\Gamma).
	\end{align*}
	Because of \eqref{eq:interface-cond-VW}, the second term in \eqref{eq:proof-integrand-I_t^Gamma} becomes
	\begin{align*}
		- \llbracket \nu^\Gamma \cdot (\epsilond V^k)(m \cdot V^k) \rrbracket
		 & = - \tr_{\partial \upsilon}(\nu^\Gamma \cdot\epsilond V^k) (m \cdot (V^k|_{\upsilon} - V^k|_{\omega}))        \\
		 & = -\alpha \tr_{\partial \upsilon}(\nu^\Gamma \cdot \epsilond V^k) (m \cdot \nu^\Gamma)                        \\
		 & = -\alpha \nu^\Gamma \cdot\restr{\epsilond}{\upsilon} (V^k|_{\omega} + \alpha \nu^\Gamma)(m \cdot \nu^\Gamma) \\
		 & = -(m \cdot \nu^\Gamma)\Bigl(\alpha \nu^\Gamma \cdot\restr{\epsilond}{\upsilon} V^k|_{\omega}
		+\alpha^2 \nu^\Gamma \cdot\restr{\epsilond}{\upsilon} \nu^\Gamma\Bigr).
	\end{align*}
	We factor out $\frac{1}{2} (m \cdot \nu^\Gamma)$ and insert
	$\alpha = - \frac{\nu^\Gamma \cdot \llbracket \epsilond \rrbracket V^k|_{\omega}}{\nu^\Gamma \cdot \epsilond|_{\upsilon} \nu^\Gamma}$, resulting in
	\begin{align*}
		I_{t}^\Gamma & = \int_s^t \!\! \oint_\Gamma \frac{1}{2} (m \cdot \nu^\Gamma) \Bigl( \llbracket V^k \cdot \epsilond V^k \rrbracket
		- 2 \alpha \nu^\Gamma \cdot\restr{\epsilond}{\upsilon} V^k|_{\omega} -2 \alpha^2 \nu^\Gamma \cdot\restr{\epsilond}{\upsilon} \nu^\Gamma\Bigr)
		\dd \varsigma \dd \tau                                                                                                                               \\
		             & = \int_s^t \!\! \oint_\Gamma\frac{1}{2} (m \cdot \nu^\Gamma) \Bigl( V^k|_{\omega} \cdot \llbracket \epsilond \rrbracket V^k|_{\omega}
		- \alpha^2 (\nu^\Gamma \cdot \restr{\epsilond}{\upsilon} \nu^\Gamma)\Bigr)\dd \varsigma \dd \tau                                                     \\
		             & =\int_s^t \!\! \oint_\Gamma\frac{1}{2} (m \cdot \nu^\Gamma) \Bigl(V^k|_{\omega} \cdot \llbracket \epsilond \rrbracket V^k|_{\omega}
		- \frac{(\nu^\Gamma \cdot \llbracket \epsilond \rrbracket V^k|_{\omega})^2}{\nu^\Gamma \cdot \epsilond|_{\upsilon} \nu^\Gamma}\Bigr)
		\dd \varsigma \dd \tau.
	\end{align*}

	We next discuss the different cases of the non-trapping condition at the interface \eqref{eq:non-trapping-interface}.
	First, if $\llbracket \epsilond \rrbracket = 0$ or $m \cdot \nu^\Gamma =0$, then the integrand of $I_{t}^\Gamma$ vanishes
	at these points.
	Second, let $\llbracket \epsilond \rrbracket >0$ and $m \cdot \nu^\Gamma \leq 0$. Since
	\[(\nu^\Gamma \cdot \llbracket \epsilond \rrbracket V^k|_{\omega})^2
		\leq (\nu^\Gamma \cdot \llbracket \epsilond \rrbracket \nu^\Gamma) \left(V^k|_{\omega} \cdot \llbracket \epsilond \rrbracket V^k|_{\omega}\right)\]
	by Cauchy-Schwarz, we estimate
	\begin{align*}
		V^k|_{\omega} \cdot \llbracket \epsilond \rrbracket V^k|_{\omega}
		- \frac{(\nu^\Gamma \cdot \llbracket \epsilond \rrbracket V^k_{\omega})^2}{\nu^\Gamma \cdot \epsilond|_{\upsilon} \nu^\Gamma}
		 & \geq \left(V^k|_{\omega} \cdot \llbracket \epsilond \rrbracket V^k|_{\omega}\right)
		\left[ 1 - \frac{\nu^\Gamma \cdot \llbracket \epsilond \rrbracket \nu^\Gamma}{\nu^\Gamma \cdot \epsilond|_{\upsilon} \nu^\Gamma} \right] \\
		 & = \left(V^k|_{\omega} \cdot \llbracket \epsilond \rrbracket V^k|_{\omega}\right)
		\frac{\nu^\Gamma \cdot \epsilond|_{\omega} \nu^\Gamma}{\nu^\Gamma \cdot \epsilond|_{\upsilon} \nu^\Gamma} > 0.
	\end{align*}
	As $m \cdot \nu^\Gamma \leq 0$, the integrand of $I_{t}^\Gamma$ is non-positive here.
	Finally, let $\llbracket \epsilond \rrbracket <0$ and $m \cdot \nu^\Gamma \geq 0$. As above we obtain
	\[ \Bigl(V^k|_{\omega} \cdot \llbracket \epsilond \rrbracket V^k|_{\omega}
		- \frac{(\nu^\Gamma \cdot \llbracket \epsilond \rrbracket V^k|_{\omega})^2}{\nu^\Gamma \cdot \epsilond|_{\upsilon} \nu^\Gamma}\Bigr) \leq 0,
		\qquad \text{and hence all together \ } I_{t}^\Gamma \leq 0.\]

	Next, we consider the integrand of $I_c^\Gamma$ given by
	\[ D\coloneq \left\llbracket W^k \cdot\left( \nu^\Gamma \times ((\mud W^k) \times m)\right)\right\rrbracket
		- \tfrac12 (m\cdot\nu^\Gamma)\llbracket\mud W^k \cdot W^k\rrbracket . \]
	To repeat the above arguments, we compute
	\begin{align*}
		\left\llbracket W^k \cdot\left( \nu^\Gamma \times ((\mud W^k) \times m)\right)\right\rrbracket
		= \left\llbracket\mud W^k \cdot W^k (\nu^\Gamma \cdot m) - ( W^k \cdot m)(\nu^\Gamma \cdot\mud W^k)\right\rrbracket.
	\end{align*}
	using the bac-cab-identity, and therefore
	\[ D=\tfrac12\left\llbracket(\nu^\Gamma \cdot m) \mud W^k \cdot W^k \right\rrbracket
		- \left\llbracket ( W^k \cdot m)(\nu^\Gamma \cdot\mud W^k)\right\rrbracket.\]
	The corresponding integral can be treated as in \eqref{eq:proof-integrand-I_t^Gamma} by means of the magnetic non-trapping condition
	$(m \cdot \nu^\Gamma)\llbracket \mud \rrbracket \leq 0$.


	\smallskip

	2.5) Next we take a look at the other terms not involving $\nabla \chi$.
	Employing Lemma~\ref{lem:error-estimate} and Lemma~\ref{lem:reg_of_pkqk}, for any $\delta>0$ we infer
	\begin{align*}
		\abs{I_\sigma} & = \abs{\int_s^t \!\! \int_\Omega \sigma \pt^k E\cdot((\mud W^k) \times \chi m) \dd x \dd \tau}
		\lesssim c_\delta \int_s^t d(\tau)\dd \tau + \delta \int_s^t \!\! \int_\Omega \abs{W^k}^2 \dd x \dd \tau,                                                       \\
		\abs{I_p}      & = \abs{\int_s^t \!\! \int_\Omega \pt (\epsilond\nabla p^k)\cdot((\mud W^k) \times \chi m) \dd x \dd \tau}                                      \\*
		               & \lesssim c_\delta\int_s^t\!\!\int_\Omega \left(d(\tau)+z^2(\tau)\right)\dd x \dd\tau + \delta\int_s^t\!\!\int_\Omega\abs{W^k}^2 \dd x \dd\tau, \\
		\abs{I_q}      & = \abs{\int_s^t \!\! \int_\Omega \epsilond V^k \cdot(\pt (\epsilond\nabla q^k) \times \chi m) \dd x \dd \tau}                                  \\*
		               & \lesssim c_\delta\int_s^t\!\!\int_\Omega \left(d(\tau)+z^2(\tau)\right)\dd x \dd\tau + \delta\int_s^t\!\!\int_\Omega\abs{V^k}^2\dd x \dd\tau,  \\
		\abs{I_f}      & = \abs{\int_s^t \!\! \int_\Omega \pt f_k\cdot((\mud W^k) \times \chi m) \dd x \dd \tau}
		\lesssim c_\delta \int_s^t z^2(\tau)\dd \tau + \delta \int_s^t \!\! \int_\Omega\abs{W^k}^2 \dd x \dd \tau ,                                                     \\
		\abs{I_g}      & = \abs{\int_s^t \!\! \int_\Omega \epsilond V^k \cdot(\pt g_k \times \chi m)\dd x \dd \tau}
		\lesssim c_\delta \int_s^t z^2(\tau)\dd \tau + \delta \int_s^t \!\! \int_\Omega \abs{V^k}^2 \dd x \dd \tau.
	\end{align*}

	2.6) Let $\chi=1$. In \eqref{eq:proof-non-est-eq} the terms involving $\nabla \chi$ thus vanish and we have seen in step 2.4)
	that we can drop $I_t^\Gamma$ and $I_c^\Gamma$. We then use the non-trapping condition \eqref{eq:non-trapping-tilde2} and absorb the integrals
	of $\abs{V^k}^2$ and $\abs{W^k}^2$ on $\upsilon \cup \cN_b^\Gamma$ by the right-hand side fixing a small $\delta>0$ in the previous step.
	This results in the inequality
	\begin{align}\label{eq:proof-chi-1-with-traces}
		\notag & \eta \int_s^t \!\! \int_{\upsilon \cup \cN_b^\Gamma} \left(\abs{V^k}^2 + \abs{W^k}^2\right) \dd x \dd \tau                 \\
		       & \lesssim \int_\Omega \abs{ (\epsilond V^k)(t)}\abs{W^k(t)} \dd x + \int_\Omega \abs{ (\epsilond V^k)(s)}\abs{W^k(s)} \dd x \\
		\notag & \seq + \int_s^t \!\! \int_{\omega \setminus \cN_b^\Gamma} \left(\abs{V^k}^2 + \abs{W^k}^2\right) \dd x \dd \tau
		+ \abs{I_t^{\partial \Omega}+ I_c^{\partial \Omega}} +\int_s^t \left(d(\tau)+z^2(\tau)\right) \dd \tau.
	\end{align}
	Thus, it remains to estimate the trace terms on $\partial \Omega$.

	\smallskip

	2.7) To do so, we choose $\chi$ with $\supp \chi \subseteq \omega \setminus \cN_b^\Gamma$ and $\chi =1$ on $\partial \Omega$. In particular,
	interface terms do not appear in \eqref{eq:proof-non-est-eq}. Equation~\eqref{eq:proof-non-est-eq} and Step~2.5) then imply
	\begin{align}\label{eq:proof-est-boundary-traces}
		\notag\abs{I_t^{\partial \Omega} + I_c^{\partial \Omega}}
		       & \lesssim \Big|\Big[\int_{\omega \setminus \cN_b^\Gamma} (\epsilond V^k) (\chi m \times \mud W^k) \Big]_s^t\Big|              \\
		\notag & \seq + \int_s^t\!\!\int_{\omega \setminus \cN_b^\Gamma} \left(\abs{V^k}^2 + \abs{W^k}^2\right) \dd x \dd \tau +\abs{I_\sigma
		+ I_p + I_f + I_q + I_g}                                                                                                              \\
		       & \lesssim \int_\Omega \abs{ (\epsilond V^k)(t)}\abs{W^k(t)} \dd x + \int_\Omega \abs{ (\epsilond V^k)(s)}\abs{W^k(s)} \dd x   \\*
		\notag & \seq +\int_s^t\!\!\int_{\omega \setminus \cN_b^\Gamma} \left(\abs{V^k}^2 + \abs{W^k}^2\right) \dd x \dd \tau                 \\*
		\notag & \seq +\delta \int_s^t \!\! \int_\Omega \left(\abs{V^k}^2 + \abs{W^k}^2\right) \dd x \dd \tau
		+ c_\delta\int_s^t \left(d(\tau)+z^2(\tau)\right) \dd \tau.
	\end{align}

	We now combine \eqref{eq:proof-chi-1-with-traces} and \eqref{eq:proof-est-boundary-traces}, again absorbing terms for a fixed, possibly decreased
	$\delta>0$. The space integral at times $t$ and $s$ are handled via Lemma~\ref{lem:VkWk-est}. Also inserting \eqref{eq:def-VkWk}, we derive
	\begin{align*}
		 & \eta\! \int_s^t \!\! \int_{\upsilon \cup \cN_b^\Gamma}\!\! \big(\abs{V^k}^2 + \abs{W^k}^2\big) \dd x \dd \tau
		\lesssim \int_\Omega \abs{V^k(t)}\abs{W^k(t)} \dd x + \int_\Omega \abs{V^k(s)}\abs{W^k(s)} \dd x                                           \\
		 & \hspace*{3cm} +\int_s^t\!\!\int_{{\omega \setminus \cN_b^\Gamma}} \left(\abs{V^k}^2
		+ \abs{W^k}^2\right) \dd x \dd \tau + c_\delta\int_s^t \left(d(\tau)+z^2(\tau)\right) \dd \tau                                             \\
		 & \lesssim e(s)+ z^2(s) + z^2(t) + (1+c_\delta)\int_s^t \left(d(\tau) +z^{3/2}(\tau)\right) \dd \tau                                      \\
		 & \seq +\int_s^t\!\!\int_{{\omega \setminus \cN_b^\Gamma}} \left(\abs{\pt^k E}^2 + \abs{\nabla p^k}^2 + \abs{W^k}^2\right) \dd x \dd \tau \\
		 & \lesssim e(s) + z^2(s) + z^2(t) + (1+c_\delta)\int_s^t \left(d(\tau) +z^{3/2}(\tau)\right) \dd \tau
		+\int_s^t\!\!\int_{{\omega \setminus \cN_b^\Gamma}}\abs{W^k}^2 \dd x \dd \tau,
	\end{align*}
	also using Lemma~\ref{lem:reg_of_pkqk}. \qedhere
\end{proof}

Next, we have to get rid of the $W^k$-term on ${\omega \setminus \cN_b^\Gamma}$ since the dissipation only acts on the electric field.
This is done for $k\geq 1$ in a similar way as in Step~3) in the proof of Proposition~3.2 in \cite{nutt2026}, however there we used a
second-order formulation of Maxwell which was then integrated by parts in time. It seems that in the present quasilinear case a second-order
approach fails because of the resulting interface terms. We thus keep the first-order system and apply a time derived multiplier.

\begin{proposition}[Observability-type Estimate II]\label{prop:Wk-est}
	Assume that \ref{hypothesis-quasi} holds. Let $\delta \in (0,1]$. For $k \in \{1,2,3\}$ and $0\leq s\leq t < T_*$, we have
	\begin{align*}
		\int_s^t \norm{W^k}_{L^2({\omega \setminus \cN_b^\Gamma})}^2 \dd x & \lesssim \int_s^t \left(\delta z(\tau) +c_\delta (z^{3/2}(\tau) + d(\tau))\right) \dd\tau \\
		                                                                   & \seq + \delta(z(s) + z(t)) + c_\delta (e(s) +z^2(s) + z^2(t)).
	\end{align*}
\end{proposition}

\begin{proof}
	Let $k \in \{1,2,3\}$.	Equations~\eqref{eq:def-VkWk}, \eqref{pde:tilde1}, and \eqref{eq:properties-Vk} imply
	\begin{align*}
		\mud W^k & = \mud \pt W^{k-1} +\mud (W^k - \pt W^{k-1})                                            \\
		         & = -\curl V^{k-1} - \mud \pt \nabla q^{k-1} - \tilde g_{k-1} + \mud (W^k - \pt W^{k-1}).
	\end{align*}
	Lemmas~\ref{lem:error-estimate}, \ref{lem:reg_of_pkqk}, and \ref{lem:diff-ptVWk-1-VWk} then lead to
	\begin{equation}\label{eq:W^k-on-omega}
		\int_s^t\!\norm{W^k}_{L^2({\omega \setminus \cN_b^\Gamma})}^2 \dd x
		\lesssim\!\int_s^t\!\norm{\curl V^{k-1}}_{L^2({\omega\setminus\cN_b^\Gamma})}^2\dd\tau + \int_s^t\!\big(d(\tau)+z^{\frac32}(\tau)\big)\mathrm{d}\tau.
	\end{equation}
	It remains to estimate the curl term. To this end, take $0 \leq \vartheta \in C^\infty(\overline \omega)$ with $\vartheta = 1$ on
	$\omega\setminus\cN_b^\Gamma$ and $\vartheta = 0$ on $\Gamma$.	Using these properties and also \eqref{eq:interface-cond-VW}, we start
	from \eqref{eq:VW}.
	Integrating by parts, we compute
	\begin{align*}
		0 & = \!\int_s^t \!\!\! \int_\omega\!\!\big(\pt(\epsilond V^{k-1})\! -\! \curl\! W^{k-1}\! +\! \sigma \pt^{k-1}E \!+ \!\pt(\epsilond \nabla p^{k-1})
		\!+ \!\pt f_{k-1}\big) \cdot \vartheta \pt V^{k-1} \mathrm{d} x \mathrm{d}\tau                                                                       \\
		  & = \int_s^t\!\!\int_\omega\big(\pt(\epsilond V^{k-1}) + \sigma \pt^{k-1}E + \pt(\epsilond \nabla p^{k-1}) + \pt f_{k-1}\big)
		\cdot \vartheta \pt V^{k-1}\dd x \dd\tau                                                                                                             \\
		  & \ + \int_s^t\!\!\int_\omega \pt (\inv\mud \mud W^{k-1})\cdot \curl (\vartheta V^{k-1}) \dd x \dd\tau
		- \left[ \int_\omega \!\!\mud W^{k-1}\cdot \curl (\vartheta V^{k-1}) \dd x \right]_s^t                                                               \\
		  & = \int_s^t\!\!\int_\omega\big(\pt(\epsilond V^{k-1}) + \sigma \pt^{k-1}E + \pt(\epsilond \nabla p^{k-1})
		+ \pt f_{k-1}\big) \cdot \vartheta \pt V^{k-1}\dd x \dd\tau                                                                                          \\
		  & \ + \int_s^t\!\!\int_\omega \big(\pt(\inv\mud) \mud W^{k-1} +\inv\mud\pt(\mud W^{k-1})\big) \cdot \curl (\vartheta V^{k-1}) \dd x \dd\tau        \\
		  & \ - \left[ \int_\omega \mud W^{k-1}\cdot \curl (\vartheta V^{k-1}) \dd x \right]_s^t                                                             \\
		  & = \int_s^t\!\!\int_\omega\big(\pt(\epsilond V^{k-1}) + \sigma \pt^{k-1}E + \pt(\epsilond \nabla p^{k-1})
		+ \pt f_{k-1}\big) \cdot \vartheta \pt V^{k-1}\dd x \dd\tau                                                                                          \\
		  & \; - \int_s^t\!\!\int_\omega \inv\mud\big( \curl V^{k-1}+ \pt (\mud \nabla q^k )+ \pt g_k\big) \cdot \curl (\vartheta V^{k-1}) \dd x \dd\tau     \\
		  & \ + \!\int_s^t\!\!\!\int_\omega \!\!\big(\pt(\inv\mud) \mud W^{k-1}\big) \cdot \curl (\vartheta V^{k-1}) \dd x \dd\tau
		- \left[\int_\omega\!\! \mud W^{k-1}\!\cdot \curl (\vartheta V^{k-1}) \dd x \right]_s^t                                                              \\
		  & = \int_s^t\!\!\int_\omega\big(\pt(\epsilond V^{k-1}) + \sigma \pt^{k-1}E + \pt(\epsilond \nabla p^{k-1}) + \pt f_{k-1}\big)
		\cdot \vartheta \pt V^{k-1}\dd x \dd\tau                                                                                                             \\
		  & \ - \int_s^t\!\!\!\int_\omega\!\! \mud^{-1}\curl\! V^{k-1}\! \cdot \vartheta\curl \!V^{k-1}\dd x \dd\tau
		- \int_s^t\!\!\int_\omega \!\!\mud^{-1}\curl\! V^{k-1}\! \cdot (\nabla\vartheta\!\times\! V^{k-1})\dd x \mathrm{d}\tau                               \\
		  & \ - \int_s^t\!\!\int_\omega \mud^{-1}\big(\pt (\mud \nabla q^k )+ \pt g_k\big) \cdot \curl (\vartheta V^{k-1}) \dd x \dd\tau                     \\
		  & \ + \int_s^t\!\!\!\int_\omega \!\!\big(\pt(\mud^{-1}) \mud W^{k-1}\big) \cdot \curl (\vartheta V^{k-1}) \dd x \dd\tau
		- \left[ \int_\omega\!\! \mud W^{k-1}\cdot \curl (\vartheta V^{k-1}) \dd x \right]_s^t.
	\end{align*}
	To justify the integration by parts for $k=3$ in the second line above, we use Theorem~\ref{thm:approx} to approximate the solution $(V^2,W^2)$ of
	the system \eqref{eq:VW} in $\cG^1$ by $(V^2_n, W^2_n)\in \cG^2$ with $\tr_t V^2_n = 0$.
	In the limit $n \to \infty$ the equality then follows. The lower bound of $\mud^{-1}$ (implied by see \eqref{eq:non-trapping-tilde1}) now yields
	\begin{align*}
		 & \int_s^t\!\!\int_{\omega \setminus \cN_b^\Gamma} \abs{\curl V^{k-1}}^2\dd x \dd\tau                                                             \\
		 & \lesssim\! \int_s^t\!\!\!\int_\omega\!\abs{\pt\epsilond V^{k-1} \!+ \epsilond \pt V^{k-1} \!+ \sigma \pt^{k-1}E + \pt(\epsilond \nabla p^{k-1})
		+ \pt f_{k-1}} \abs{\pt V^{k-1}}\dd x \mathrm{d}\tau                                                                                               \\
		 & \ + \!\int_s^t\!\!\!\int_{\cN_b^\Gamma \cap \omega}\! \abs{\curl \!V^{k-1}}\abs{V^{k-1}}\dd x \dd\tau
		+ \int_s^t\!\!\!\int_\omega \!\abs{\pt (\mud \nabla q^k )+ \pt g_k} \abs{\curl (\vartheta V^{k-1})} \dd x \mathrm{d}\tau                           \\
		 & \ + \!\int_s^t\!\!\!\int_\omega \abs{\pt(\inv\mud) \mud W^{k-1}} \abs{\curl(\vartheta V^{k-1})} \dd x \dd\tau
		- \left[ \int_\omega \mud W^{k-1}\cdot \curl (\vartheta V^{k-1}) \dd x \right]_s^t                                                                 \\
		 & \eqqcolon I_1 + I_2 + I_3 + I_4 + I_5.
	\end{align*}
	We now insert Definition~\ref{def:VkWk}, \eqref{eq:properties-Vk}, and
	\[	\curl(\vartheta V^{k-1}) = \nabla \vartheta \times (\pt^{k-1}E - \nabla p^{k-1}) + \vartheta \curl(\pt^{k-1}E). \]
	Lemmas~\ref{lem:error-estimate} and \ref{lem:reg_of_pkqk} then imply
	\begin{align*}
		I_1 & \lesssim \int_s^t \big(d^{1/2}(\tau) + z(\tau)\big)\big(d^{1/2}(\tau) + z(\tau)\big) \dd \tau
		\lesssim \int_s^t \big(d(\tau) + z^2(\tau)\big) \dd \tau,                                                                                      \\
		I_2 & \lesssim \int_s^t z^{1/2}(\tau)\big(d^{1/2}(\tau)+ z(\tau)\big) \dd \tau
		\lesssim \int_s^t \big(\delta z(\tau) + c_\delta d(\tau) + z^{3/2}(\tau)\big) \dd \tau,                                                        \\
		I_3 & \lesssim \int_s^t \big(d^{1/2}(\tau) + z(\tau)\big)\big(d^{1/2}(\tau) + z^{1/2}(\tau) + z(\tau)\big) \dd \tau                            \\
		    & \lesssim \int_s^t \big(\delta z(\tau) + c_\delta d(\tau) + z^{3/2}(\tau)\big) \dd \tau,                                                  \\
		I_4 & \lesssim \int_s^t z^{1/2}(\tau)\big(d^{1/2}(\tau) +z^{1/2}(\tau)+ z(\tau)\big)\big(d^{1/2}(\tau) + z^{1/2}(\tau) + z(\tau)\big) \dd \tau \\
		    & \lesssim \int_s^t \big(d(\tau) + z^{3/2}(\tau)\big) \dd \tau,                                                                            \\
		I_5 & \lesssim \sum_{r = s,t} \big(e^{1/2}(r) + d^{1/2}(r) + z(r)\big)\big(e^{1/2}(r) + z^{1/2}(r) + z(r)\big)                                 \\
		    & \lesssim \sum_{r = s,t} \delta z(r) + c_\delta\big(e(r) + d(r) + z^2(r)\big),
	\end{align*}
	for $\delta>0$. To remove $e(t)$, we invoke the energy estimate Lemma~\ref{lem:energy-decreasing} arriving at
	\begin{align*}
		\int_s^t\!\!\int_{\omega \setminus \cN_b^\Gamma} \abs{\curl V^{k-1}}^2\dd x \dd\tau & \lesssim \int_s^t \left(\delta z(\tau) +c_\delta (z^{3/2}(\tau) + d(\tau))\right) \dd x \dd\tau \\
		                                                                                    & \seq + \delta(z(s) + z(t)) + c_\delta (e(s) +z^2(s) + z^2(t)).
	\end{align*}
	The claim follows by inserting this estimate into \eqref{eq:W^k-on-omega}.
\end{proof}

Proposition~\ref{prop:hom_obsI} and Proposition~\ref{prop:Wk-est} then imply the core observability estimate.
\begin{corollary}[Observability-type Estimate III]\label{cor:obs-est}
	Assume that \ref{hypothesis-quasi} holds. Let $\delta \in (0,1]$. For $k \in \{1,2,3\}$ and $0\leq s\leq t < T_*$, we have
	\begin{align*}
		\int_s^t \!\! \int_{\upsilon \cup \cN_b^\Gamma} \left(\abs{V^k}^2 + \abs{W^k}^2\right) \dd x \dd \tau & \lesssim \int_s^t \left(\delta z(\tau) +c_\delta (z^{3/2}(\tau) + d(\tau))\right) \dd x \dd\tau \\
		                                                                                                      & \seq + \delta(z(s) + z(t)) + c_\delta (e(s) +z^2(s) + z^2(t)).
	\end{align*}
\end{corollary}

Before we derive an observability estimate for $\pt^k(E,H)$, we first control the magnetic contribution to $\int e \dd \tau$.
This inequality can be shown as in Lemma~4.3 and Proposition~4.4 in \cite{lasiecka2019} (cf.\ Equation~(4.13) therein), where simply
connected $\Omega$ without an interface were treated.
\begin{lemma}\label{lem:pt-H-estimate}
	Let \ref{hypothesis-quasi} hold.
	For $k\in\{0,1,2,3\}$ and $0\leq s \leq t < T_*$ we have
	\[\int_s^t \smallnorm{\pt^k H}^2_{L^2(\Omega)} \dd \tau \lesssim
		\left(e(s)+e(t)\right)+ \int_s^t\smallnorm{\pt^k E}_{L^2(\Omega)}^2\dd\tau+ \int_s^t z^{3/2}(\tau) \dd \tau.\]
\end{lemma}

\begin{proof}
	Using the positive definiteness of $\mu$ and Lemma~\ref{lem:error-estimate}, we first write
	\begin{equation}\label{proof:eq:ptk-H-est-3}
		\int_s^t \smallnorm{\pt^k H}^2_{L^2(\Omega)} \dd \tau \lesssim \int_s^t \!\! \int_\Omega \pt^k(\mu H) \cdot \pt^k H \dd x \dd \tau
		+ \int_s^t z^{3/2}(\tau) \dd \tau.
	\end{equation}
	Observe that $\pt^k (\mu H)$ belongs to $H^\Sigma_{n0}(\div 0) = \ran(\curl_E)$, by \eqref{eq:maxwell}, \eqref{eq:div-cond-initial},
	and Lemma~\ref{lem:curl-with-topology}. Thus there exists a unique vector potential
	$u^k \in C\big([0,T_*],H_{t0}(\curl) \cap H^{\partial\Omega}(\div0)\big)$ such that
	\begin{equation}\label{proof:eq:ptk-H-est-1}
		\curl u^k = \pt^k (\mu H) \qquad \text{and} \qquad \norm{u^k}_{L^2} \lesssim \smallnorm{\pt^k (\mu H) }_{L^2}.
	\end{equation}
	By \eqref{eq:maxwell} we have $\pt (\mu H) = -\curl E$. For $k=0$, taking the the difference quotient, we see that
	$u^0 \in C^1\big([0,T_*],H_{t0}(\curl)\cap H^{\partial\Omega}(\div0)\big)$. Moreover, the time derivative of \eqref{proof:eq:ptk-H-est-1} yields
	$\curl(\pt u^0) = \pt(\mu H) = -\curl E$ and hence $\pt u^0 + E$ belongs to $H_{t0}(\curl 0)$.
	Theorem~2.8' of \cite{cessenat1996} then yields
	\[\pt u^0 + E = \nabla \phi_0^0 + h_2^0 \]
	for unique and orthogonal $\phi_0^0 \in H^1_0(\Omega)$ and $h_2^0 \in \bH_2$, cf.\ \eqref{def:cohomology-spaces}.
	Furthermore, Theorems~2.9 and 2.10' in \cite{cessenat1996} provide the orthogonal decomposition
	\[ L^2(\Omega) = H^{\partial\Omega}(\div 0) \oplus \nabla H^1_0(\Omega) \oplus \bH_2, \quad \text{and hence \ }
		E =- \pt u^0 +\nabla \phi_0^0 +h_2^0 .\]
	In particular, orthogonality yields
	\begin{equation}\label{eq:potential-k0-est}
		\norm{\pt u^0 }_{L^2} \leq \smallnorm{E}_{L^2}\,.
	\end{equation}
	For $k \in \{1,2,3\}$, we obtain $\curl u^k = \pt^k(\mu H) = - \curl(\pt^{k-1} E)$ and therefore $u^k + \pt^{k-1} E \in H_{t0}(\curl 0)$.
	Again, Equation~(2.107), Theorem~2.8', 2.9 and 2.10' of \cite{cessenat1996} give the orthogonal decomposition
	\[\pt^{k-1} E =- u^k +\nabla \phi_0^k+ h_2^k,\]
	for unique $\phi_0^k \in H^1_0(\Omega)$ and $h_2^k \in \bH_2$. In particular the regularity of the left-hand side and the decomposition imply
	that $u^k \in C^1\big([0,T_*],H_{t0}(\curl) \cap H^{\partial\Omega}(\div0)\big)$ and
	\begin{equation}\label{eq:potential-k123-est}
		\norm{\pt u^k }_{L^2} \leq \smallnorm{\pt^k E}_{L^2}.
	\end{equation}

	For all $k \in \{0,1,2,3\}$, integration by parts, \eqref{pde:pt-maxwell} and \eqref{proof:eq:ptk-H-est-1} lead to
	\begin{align}\label{proof:eq:ptk-H-est-4}
		\notag\int_s^t \!\! \int_\Omega \pt^k(\mu H) \cdot \pt^k H \dd x \dd\tau & = \int_s^t \!\! \int_\Omega \curl u^k \cdot \pt^k H \dd x
		= \int_s^t \langle u^k, \curl\pt^k H \rangle_{H_{t0}(\curl)} \dd\tau                                                                                                                          \\
		\notag                                                                   & =\int_s^t \bigl\langle u^k, \pt (\epsilonk \pt^k E) + \sigma \pt^k E + \pt f_k\bigr\rangle_{H_{t0}(\curl)} \dd\tau \\
		\notag                                                                   & \lesssim \int_s^t \bigl\langle u^k, \pt (\epsilonk \pt^k E)\bigr\rangle_{H_{t0}(\curl)} \dd\tau
		+ \delta \int_s^t \smallnorm{\pt^k H}^2_{L^2(\Omega)} \dd \tau                                                                                                                                \\
		                                                                         & \seq + c_\delta \int_s^t \left(\smallnorm{\pt^k E}_{L^2(\Omega)}^2 + z^2(\tau)\right) \dd \tau.
	\end{align}
	The first term on the right is rewritten as
	\begin{equation*}
		\int_s^t \bigl\langle u^k, \pt (\epsilonk \pt^k E)\bigr\rangle_{H_{t0}(\curl)} \dd\tau
		=\left[ \int_\Omega u^k \cdot (\epsilonk \pt^k E)\dd x \right]_s^t - \int_s^t \!\! \int_\Omega \pt u^k \cdot (\epsilonk \pt^k E)\dd x \dd \tau.
	\end{equation*}
	Inserting \eqref{proof:eq:ptk-H-est-1}, \eqref{eq:potential-k0-est}, and \eqref{eq:potential-k123-est} into the above identity, we derive
	\begin{align*}
		\int_s^t \bigl\langle u^k, \pt (\epsilonk \pt^k E)\bigr\rangle_{H_{t0}(\curl)} \dd\tau
		\lesssim \left(e_k(s)+e_k(t)\right) + \int_s^t \smallnorm{\pt^k E}_{L^2(\Omega)}^2 \dd \tau.
	\end{align*}
	For fixed small $\delta>0$, this estimate, \eqref{proof:eq:ptk-H-est-3}, and \eqref{proof:eq:ptk-H-est-4} imply the claim.
\end{proof}

Finally, we establish an observability-type estimate for $\pt^k(E,H)$, up to small and higher-order error terms.
\begin{theorem}\label{thm:obs-est}
	Let \ref{hypothesis-quasi} hold.
	For $0 \leq s \leq t < T_*$ and $\delta\in(0,1]$, we obtain
	\begin{align*}
		\int_s^t e(\tau) \dd\tau & \lesssim e(t) + c_\delta \left(z(s) + z^{2}(t) + \int_s^t (d(\tau) + z^{3/2}(\tau)) \dd \tau\right) \\
		                         & \seq + \delta \left(z(t) + \int_s^t z(\tau) \dd \tau\right).
	\end{align*}
\end{theorem}

\begin{proof}
	By definition and Lemma~\ref{lem:cE-estimate} we have
	\[\int_s^t e(\tau)\dd\tau \lesssim \max_{k\in\{1,2,3\}} \int_s^t \left(\norm{\epsilonk \pt^k E}_{L^2}^2
		+ \norm{\muk \pt^k H}_{L^2}^2 \right) \dd \tau.\]
	Lemma~\ref{lem:pt-H-estimate}, Definition~\ref{def:VkWk} and Corollary~\ref{cor:obs-est} then yield
	\begin{align*}
		\int_s^t e(\tau)\dd\tau & \lesssim \max_{k\in\{1,2,3\}} \int_s^t \norm{\epsilonk\pt^k E}_{L^2}^2 \dd\tau + e(s)+e(t)+ \int_s^t z^{3/2}(\tau)\dd\tau \\
		                        & \lesssim \max_{k\in\{1,2,3\}} \int_s^t\left(\norm{V^k}_{L^2}^2 + \norm{\nabla p^k}_{L^2}^2\right) \dd \tau
		+ e(s) + e(t) + \int_s^t z^{3/2}(\tau) \dd \tau                                                                                                     \\
		                        & \lesssim c_\delta \max_{k\in\{1,2,3\}} \int_s^t \!\! \int_{\omega} \abs{V^k}^2 \dd x \dd \tau
		+ \max_{k\in\{1,2,3\}} \int_s^t \norm{\nabla p^k}_{L^2}^2 \dd\tau + e(t)                                                                            \\*
		                        & \seq + \delta \Big(z(t) + \int_s^t z(\tau) \dd \tau\Big)
		+c_\delta \Big(z(s)+ z^{2}(t) + \int_s^t (d(\tau) + z^{3/2}(\tau)) \dd \tau\Big)
	\end{align*}
	for any $\delta \in (0,1]$. Since $\sigma\ge\sigma_0>0$ on $\omega$, Definition~\ref{def:VkWk} and Lemma~\ref{lem:reg_of_pkqk} lead to
	\begin{align*}
		\int_s^t\!\! e(\tau)\dd\tau & \lesssim \max_{k\in\{1,2,3\}} \int_s^t \!\! \int_{\omega} \abs{\sigma \pt^{k}E}^2 \dd x \dd \tau
		+ \max_{k\in\{1,2,3\}} \int_s^t\norm{\nabla p^k}_{L^2}^2 \dd\tau + e(t)                                                                             \\*
		                            & \seq + \delta \Big(z(t) + \int_s^t z(\tau) \dd \tau\Big)
		+c_\delta \Big(z(s)+ z^{2}(t) + \int_s^t (d(\tau) + z^{\frac32}(\tau)) \dd \tau\Big)                                                                \\*
		                            & \lesssim e(t) + c_\delta\Big[z(s) + z^{2}(t) + \!\int_s^t \!\!(d(\tau)\! +\! z^{\frac32}(\tau)) \mathrm{d}\tau\Big]\!
		+ \delta \Big[z(t) +\! \int_s^t\!\! z(\tau) \mathrm{d}\tau\Big]. \qedhere
	\end{align*}
\end{proof}

\section{Regularity analysis}
\label{sec:regularity-analysis}

Our proof of exponential stability will be based on Theorem~\ref{thm:obs-est} and Lemma~\ref{lem:energy-decreasing}. These results allow us to
control time derivatives of $(E,H)$ in $L^2$, but so far the estimates include small error terms also involving space regularity.
We thus need a version of Theorem~\ref{thm:obs-est} for spatial derivatives, i.e., we have to bound $z$ on the left by $e$
and higher-order or small terms. Here it is crucial that the constants do not depend on time and that we
estimate $z$ in $t$ uniformly
and in $L^1$. For uniformly positive conductivity this has been carried out in \cite{lasiecka2019}. However, the jump of $\sigma$
leads to the interface $\Gamma$, which requires sophisticated additional arguments.

\begin{theorem}[Regularity analysis]\label{thm:regularity-analysis}
	Let~\ref{hypothesis-quasi} hold. For $0\leq s\leq t < T_*$ the solutions $(E,H)$ of~\eqref{eq:maxwell}--\eqref{eq:maxwell4} satisfy
	\[
		z(t) + \int_s^t z(\tau)\dd\tau\lesssim z(s) + z^2(t)+ e(t) + \int_s^t\left(e(\tau) + z^{3/2}(\tau)\right)\dd\tau.
	\]
\end{theorem}
To prove this result, we split the solutions in different spatial regions via cutoffs and reduce the estimate to lower spatial norms.
We use the collars $\cN_a = \{ x\in\overline\Omega\mid\mathrm{dist}(x,\partial\Omega)<a\}$ and
$\cN^\Gamma_a = \{ x\in\overline\Omega\mid\mathrm{dist}(x,\Gamma)<a\}$ around $\partial\Omega$ and $\Gamma$.
Choose $a>0$ such that
\begin{itemize}[label=--]
	\item $\cN_a\cap\cN^\Gamma_a = \emptyset$,
	\item $\partial \cN_a$ and $\partial\cN^\Gamma_a$ retain the $C^5$ regularity of $\partial \omega$, see Lemma~14.16 in \cite{gilbarg2001},
	\item the basis $\{\tau_1^\Xi,\tau_2^\Xi,\nu^\Xi\}$ of normal and tangential vectors can be extended from $\partial\Omega$,
	      resp.\ $\Gamma$, to a basis also denoted by $\{\tau_1^\Xi,\tau_2^\Xi,\nu^\Xi\}$ on $\cN_a$, resp.\ $\cN^\Gamma_a$.
\end{itemize}
We introduce smooth cutoffs $\chi_{\partial\Omega},\chi_\omega,\chi_\Gamma$ and $\chi_\upsilon\geq 0$, satisfying $\sum_{\Xi\in\{\partial\Omega,\omega,\Gamma,\upsilon\}}\chi_\Xi = 1$ on $\Omega$.
To this end, we take $\chi_{\partial\Omega}$ with $\chi_{\partial\Omega} = 1$ on $\cN_{a/2}$ and $\supp\chi_{\partial\Omega}\subseteq\cN_{a}$, and $\chi_\Gamma$ with $\chi_\Gamma = 1$ on $\cN^\Gamma_{a/2}$ and $\supp\chi_\Gamma\subseteq\cN^\Gamma_{a}$.
Then we set $\chi_\omega = \1_\omega (1-\chi_{\partial\Omega} -\chi_\Gamma)$ and $\chi_\upsilon = \1_\upsilon (1-\chi_\Gamma)$.
We define the quantities
\[
	y_\Xi^{k,j} (t)\coloneqq\norm{\pt^k [\chi_\Xi (E(t),H(t))]}^2_{H_\Gamma^{j}(\Omega)}
\]
for $\Xi\in\{\upsilon,\Gamma,\omega,\partial\Omega\}$, $j,k\in\{0,1,2,3\}$, $j\leq 3-k$, and $t\in[0,T_*)$. We have to show
\begin{align*}
	y_\Xi^{k,j}(t) + \int_s^t y_\Xi^{k,j}\dd\tau \lesssim c_\theta\left[z(s)+ e(t) + z^2(t) + \int_s^t(e + z^{\frac32})\dd\tau\right]
	+ \theta\!\int_s^t z\dd\tau
\end{align*}
for small $\theta>0$. Note that for $j = 0$ this is true with $\theta = 0$ since $y_\Xi^{k,0}(t)\lesssim e(t)$. In particular for $k=3$
there is nothing to show.
In the cases $\Xi\neq\Gamma$ the proofs are essentially contained in \cite{lasiecka2019} and \cite{pokojovy2020}, so that we can omit them.
See Lemma~\ref{lem:upsilon-est}.

We first derive an estimate for the magnetic field on $\Omega_\pm$ from the div-curl inequality of Theorem~\ref{thm:div-curl-perfect}.
Here we trade a space derivative of $H$ to a time derivative of $E$ which we can control in our iteration scheme. We provide a short
proof since the result is only partly contained in \cite{lasiecka2019} and \cite{pokojovy2020}.
\begin{lemma}\label{lem:H-est}
	Let~\ref{hypothesis-quasi} hold. For $0\leq s\leq t < T_*$, $j,k\in\N_0$ with $j+k\leq 2$ the solutions $(E,H)$ of~\eqref{eq:maxwell}--\eqref{eq:maxwell4} satisfy
	\[ \norm{\pt^k H}_{H^{j+1}_\Gamma} \lesssim e^{1/2}(t) + \norm{\pt^{k}E}_{H^j_\Gamma} + \norm{\pt^{k+1}E}_{H^j_\Gamma} +z(t).\]
\end{lemma}

\begin{proof}
	Let $k\leq 2$.
	By \eqref{pde:pt-maxwell-IC}, \eqref{eq:div-condition-H-kgeq1},
	and Lemma~\ref{lem:characterization-H-div-H-curl} the fields
	\[
		\pt^k H\in H(\curl),\qquad\hat\mu_\mathrm{d} \pt^k H + \tilde{g}_{k-1} = \pt^k(\mu H)\in H_{n0}(\div 0)
	\]
	satisfy the assumptions of Theorem~\ref{thm:div-curl-perfect}. Together with~\eqref{pde:tilde1} and Lemma~\ref{lem:error-estimate},
	Theorem~\ref{thm:div-curl-perfect} implies
	\begin{align*}
		\norm{\pt^k H}_{H^{j+1}_\Gamma} & \lesssim\norm{\pt^k H}_{L^2} + \norm{\curl\pt^k H}_{H^{j}_\Gamma} + \norm{\tilde g_{k-1}}_{H^{j+1}_\Gamma} \\
		                                & \lesssim e^{1/2}(t) + \norm{\pt^{k}E}_{H^{j}_\Gamma} + \norm{\pt^{k+1}E}_{H^{j}_\Gamma} +z(t).\qedhere
	\end{align*}
\end{proof}

The next lemma provides the estimates on $\upsilon$ and on $\omega$ off their boundaries and near $\partial\Omega$. The first one
can be shown as in Section~6 of \cite{pokojovy2020} by differentiating in $x$ and applying a div-curl estimate.
Here an additional term near the interface $\Gamma$ appears. Its space regularity is controlled by time derivatives, lower-order and
error terms, and additional contributions at the boundary. A proof based on \cite{pokojovy2020} is given in Lemma~4.5.3 of
\cite{nutt2026-diss}. The part of $\omega$ away from $\partial\omega$ can be handled as in \cite{lasiecka2019},
see Equation~(7.6) therein. Here, the magnetic field is estimated via Lemma~\ref{lem:H-est}, whereas the electric one can be treated
differently by an energy estimate due to the presence of the conductivity. See Lemma~4.5.4 in \cite{nutt2026-diss}.

Up to this point differentiation in space does not change the boundary or interface conditions since the cutoffs vanish there.
Near $\partial\Omega$ or $\Gamma$, we may still differentiate in tangential direction picking up lower-order terms. However the
derivatives in normal direction require an additional argument based on the structure of the Maxwell system. This step does not
depend on the boundary conditions.
The case $\partial\Omega$ is essentially contained in Section~7 of \cite{lasiecka2019}. Lemma~4.5.5 of \cite{nutt2026-diss} gives
a proof adapted to the present situation.

\begin{lemma}\label{lem:upsilon-est}
	Let~\ref{hypothesis-quasi} hold, $j,k\in\N_0$ with $j+k\leq 2$, and $\theta>0$. For $0\leq s\leq t < T_*$,
	the solutions $(E,H)$ of \eqref{eq:maxwell}--\eqref{eq:maxwell4} satisfy
	\begin{align*}
		 & y^{k,j+1}_\upsilon(t)\lesssim z^2(t) + y^{k+1,j}_\upsilon(t)+y^{k,j}_\upsilon(t) + y^{k,j}_\Gamma(t),                        \\
		 & y_\omega^{k,j+1}(t) + \sigma_0\int_s^t y_\omega^{k,j+1}(\tau)\dd\tau                                                         \\
		 & \seq\lesssim z(s) + c_\theta\!\int_s^t\!\big[e(\tau) + z^{3/2}(\tau) + \norm{\pt^{k}E}_{H^j_\Gamma}^2
		+ \norm{\pt^{k+1}(E,H)}_{H^j_\Gamma}^2\big]\mathrm{d}\tau + \theta\!\int_s^t\! \!z(\tau)\dd\tau,                                \\
		 & y^{k,j+1}_{\partial\Omega}(t) + \int_s^t y^{k,j+1}_{\partial\Omega}(\tau)\dd\tau                                             \\
		 & \seq \lesssim z(s) +e(t)+ z^2(t)+ \theta\int_s^t z(\tau)\dd\tau+ (1+c_\theta)\bigg[\int_s^t (e(\tau) + z^{3/2}(\tau))\dd\tau \\
		 & \seq\seq + \norm{\pt^{k}E(t)}_{H^{j}_\Gamma}^2 + \norm{\pt^{k+1}E(t)}_{H^{j}_\Gamma}^2
		+ \int_s^t\left(\norm{\pt^{k}E(\tau)}_{H^{j}_\Gamma}^2 + \norm{\pt^{k+1}E(\tau)}_{H^{j}_\Gamma}^2\right)\dd\tau                 \\
		 & \seq\seq + y^{k+1,j}_{\partial\Omega}(t) + \int_s^t y^{k+1,j}_{\partial\Omega}(\tau)\dd\tau\bigg].
	\end{align*}
\end{lemma}

We now deal with the Sobolev norms of $\pt^k (E,H)$ near the interface $\Gamma$. We first consider each side of $\Gamma$ individually.
The tangential derivative of $\chi_\Gamma(E,H)$ on $\omega$ can be estimated using the energy identity exploiting that
$\sigma\geq\sigma_0 >0$ on $\omega$. The normal derivatives are then reduced to the tangential case via the Maxwell system.
On $\upsilon$ we bound the fields using a div-curl lemma, where the traces can be controlled via the interface conditions and
lower-order estimates.

\begin{lemma}\label{lem:Gamma-est}
	Let~\ref{hypothesis-quasi} hold, $k, j\in\N_0$ with $k+j\leq 2$, and $\theta>0$. For $0\leq s\leq t < T_*$,
	the solutions $(E,H)$ of~\eqref{eq:maxwell}--\eqref{eq:maxwell4} satisfy
	\begin{align*}
		 & y^{k, j+1}_{\Gamma}(t) + \int_s^t y^{k, j+1}_{\Gamma}(\tau)\dd\tau                                      \\
		 & \lesssim e(t) +z(s) +z^2(t)+ \theta\!\int_s^t\!\! z(\tau)\dd\tau+ c_\theta\bigg[\int_s^t\!\!(e(\tau)
		+z^{\frac32}(\tau))\mathrm{d}\tau + \norm{\pt^k (E,H)(t)}^2_{H_\Gamma^j}                                   \\
		 & \seq + \norm{\pt^{k+1} (E,H)(t)}^2_{H_\Gamma^j} + \int_s^t\left(\norm{\pt^k (E,H)(\tau)}^2_{H_\Gamma^j}
			+ \norm{\pt^{k+1} (E,H)(\tau)}^2_{H_\Gamma^j}\right)\dd\tau\bigg].
	\end{align*}
\end{lemma}

\begin{proof}
	To simplify the notation slightly, we drop the superscript of $\nu^\Gamma = \nu$ and $\tau^\Gamma_i = \tau_i$ when used in a subscript.
	From~\eqref{eq:interface-cond} we recall the interface conditions
	\begin{equation}\label{eq:maxwell-currents}
		\llbracket\pt^k E\times\nu^\Gamma\rrbracket = 0,\qquad\llbracket\pt^k H\times\nu^\Gamma\rrbracket= 0.
	\end{equation}

	1) For $0 \neq \alpha = (\gamma,\beta)\in\N_0^3\cong\N_0\times\N_0^2$ we write
	$\partial_\nu^\gamma\partial_\tau^\beta = \partial_\nu^{\gamma}\partial_{\tau_1}^{\beta_1}\partial_{\tau_2}^{\beta_2}$.

	1.1) Let $\abs\beta\leq 3-k$. The tangential derivative $\partial_\tau^\beta\pt^k(\chi_\Gamma (E,H))$ of the localized solution
	solves the system
	\begin{align}\label{pde:tangential-derivatives}
		\notag\epsilond\pt\left[\partial_\tau^\beta\pt^k(\chi_\Gamma E)\right]
		                                                            & = \curl\left(\partial_\tau^\beta\pt^k(\chi_\Gamma H)\right)
		-\sigma\left[\partial_\tau^\beta\pt^k(\chi_\Gamma E)\right] -\frakf_{\beta,k},                                                                                                           \\
		\notag\frakf_{\beta,k}                                      & \coloneqq - \partial_\tau^\beta\Big([\chi_\Gamma,\curl]\pt^k H\Big)
		- [\partial_\tau^\beta,\curl] (\chi_\Gamma\pt^k H) + \partial_\tau^\beta (\chi_\Gamma\tilde f_k)                                                                                         \\
		\notag                                                      & \seq+ \sum_{0\leq b <\beta}\binom{\beta}{b}\left(\partial_\tau^{\beta - b}\epsilond\partial_\tau^b\pt^{k+1}(\chi_\Gamma E)
		+\partial_\tau^{\beta - b}\sigma\partial_\tau^b\pt^{k}(\chi_\Gamma E)\right),                                                                                                            \\
		\mud\pt\left[\partial_\tau^\beta\pt^k(\chi_\Gamma H)\right] & = -\curl\left(\partial_\tau^\beta\pt^k(\chi_\Gamma E)\right) -\frakg_{\beta,k},                                            \\
		\notag\frakg_{\beta,k}                                      & \coloneqq\partial_\tau^\beta\Big([\chi_\Gamma,\curl]\pt^k E\Big)+ [\partial_\tau^\beta,\curl] (\chi_\Gamma\pt^k E)
		+ \partial_\tau^\beta (\chi_\Gamma\tilde g_k)                                                                                                                                            \\
		\notag                                                      & \seq + \sum_{0\leq b <\beta}\binom{\beta}{b}\partial_\tau^{\beta - b}\mud\partial_\tau^b\pt^{k+1}(\chi_\Gamma H),
	\end{align}
	on both sides of $\Gamma$. Note that $[\chi_\Gamma,\curl]$ is just a multiplier and
	$[\partial_\tau^\beta,\curl]$ is of order $\abs\beta$. Trivially, the fields also satisfy the boundary conditions
	\[ \tr_t\partial_\tau^\beta\pt^k (\chi_\Gamma E) = 0 ,\quad \tr_n\partial_\tau^\beta\pt^k(\mud\chi_\Gamma H) = 0
		\quad \text{on }\partial\cN_a^{\Gamma}.\]
	The tangential derivative commutes with the trace. Indeed, one can first establish this result for smooth functions and then
	apply the continuity of the operators.
	Hence, $\partial_\tau^\beta\pt^k (\chi_\Gamma E,\chi_\Gamma H)$ satisfies
	\begin{alignat*}{4}
		\llbracket\tr_t\partial_\tau^\beta\pt^k (\chi_\Gamma E)\rrbracket & = -\sum_{0\leq b <\beta }\binom{\beta}{b}\llbracket\tr\partial^{b}_{\tau}\pt^k(\chi_\Gamma E)\rrbracket\times\partial_\tau^{\beta - b}\nu^\Gamma &  & \eqqcolon\varphi_{\beta,k}\quad\text{on }\Gamma, \\
		\llbracket\tr_t\partial_\tau^\beta\pt^k (\chi_\Gamma H)\rrbracket & = -\sum_{0\leq b <\beta }\binom{\beta}{b}\llbracket\tr\partial^{b}_{\tau}\pt^k(\chi_\Gamma H)\rrbracket\times\partial_\tau^{\beta - b}\nu^\Gamma &  & \eqqcolon\psi_{\beta,k}\quad\text{on }\Gamma,
	\end{alignat*}
	where only terms of order strictly smaller than $\abs\beta$ appear on the right-hand side.

	\smallskip

	1.2) We start by deriving the energy equality
	\begin{align}\label{eq:differentiated-energy-est}
		       & \frac12\int_{\cN^\Gamma_a}\left[\epsilond\partial_\tau^\beta\pt^k (\chi_\Gamma E(t))\cdot\partial_\tau^\beta\pt^k (\chi_\Gamma E(t))
		+ \mud\partial_\tau^\beta\pt^k (\chi_\Gamma H(t))\cdot\partial_\tau^\beta\pt^k (\chi_\Gamma H(t))\right]\dd x \notag                            \\*
		       & \seq +\int_s^t\!\!\int_{\cN^\Gamma_a}\sigma\partial_\tau^\beta\pt^k (\chi_\Gamma E(\tau))
		\cdot\partial_\tau^\beta\pt^k (\chi_\Gamma E(\tau))\dd x\dd\tau                                                                                 \\
		\notag & = \frac12\int_{\cN^\Gamma_a}\left[\epsilond\partial_\tau^\beta\pt^k (\chi_\Gamma E(s))\cdot\partial_\tau^\beta\pt^k (\chi_\Gamma E(s))
		+ \mud\partial_\tau^\beta\pt^k (\chi_\Gamma H(s))\cdot\partial_\tau^\beta\pt^k (\chi_\Gamma H(s))\right]\dd x                                   \\
		\notag & \seq + \frac12\!\int_s^t\!\!\int_{\cN^\Gamma_a}\!\left[\pt\epsilond\partial_\tau^\beta\pt^k (\chi_\Gamma E)
			\cdot\partial_\tau^\beta\pt^k (\chi_\Gamma E)
		+ \pt\mud\partial_\tau^\beta\pt^k (\chi_\Gamma H)\cdot\partial_\tau^\beta\pt^k (\chi_\Gamma H)\right]\mathrm{d} x\dd\tau                        \\
		\notag & \seq -\int_s^t\!\!\int_{\cN^\Gamma_a}	\left[\partial_\tau^\beta\pt^k (\chi_\Gamma E)\cdot\frakf_{\beta,k}
		+ \partial_\tau^\beta\pt^k (\chi_\Gamma H)\cdot\frakg_{\beta,k}\right]\dd x\dd\tau                                                              \\
		\notag & \seq + \int_s^t\!\!\oint_\Gamma\left(\varphi_{\beta,k}\cdot\tr_{\tau,\partial\upsilon}\partial_\tau^\beta\pt^k (\chi_\Gamma H)
		-\tr_{\tau,\partial\omega}\partial_\tau^\beta\pt^k (\chi_\Gamma E)\cdot\psi_{\beta,k}\right)\dd\varsigma\dd\tau.
	\end{align}
	For $\partial_\tau^\beta\pt^k (\chi_\Gamma E, \chi_\Gamma H) \in \cG^1_\Gamma$ this follows from \eqref{pde:tangential-derivatives} and
	an integration by parts argument. However, if $k + \abs\beta = 3$, we have to approximate the solutions suitably. Since $k \leq 2$ there is
	at least one tangential derivative $\partial_\tau^{(\beta_1,\beta_2)}$ with $\beta_i \geq 1$ for $i \in \{1,2\}$.
	Using Theorem~\ref{thm:approx} with $m=1=a$,
	we approximate $(u,v) = \pt^k\partial_\tau^{\beta - e_i}(\chi_\Gamma E, \chi_\Gamma H) \in \cG^1_\Gamma$
	on $\cN_a^\Gamma$ in $\cG^1_\Gamma$ by $(u_n,v_n)\in \cG^2_\Gamma$ solving
	\begin{alignat*}{-1}
		\begin{split}
			\epsilond \pt u_n & = \curl v_n - \sigma u_n - \frakf_{\beta-e_i,k,n}, \\
			\mud \pt v_n      & = - \curl u_n - \frakg_{\beta-e_i,k,n},
		\end{split}
		\qquad                                    &                                   & t\geq 0, x \in \cN_a^\Gamma\setminus\Gamma,                                         \\
		\llbracket u \times \nu^\Gamma \rrbracket & =\varphi_{\beta-e_i,k,n},\quad
		\llbracket v \times \nu^\Gamma \rrbracket = \psi_{\beta-e_i,k,n},\quad
		                                          &                                   & t\geq 0, x \in \Gamma,                                                              \\
		u \times \nu                              & = 0,                              &                                             & t\geq 0, x \in \partial \cN_a^\Gamma, \\
		u(0,x)                                    & = u_{0,n}, \quad v(0,x)= v_{0,n}, &                                             & x \in \cN_a^\Gamma\setminus\Gamma,
	\end{alignat*}
	where the approximations $\frakf_{\beta-e_i,k,n}, \frakg_{\beta-e_i,k,n}, \varphi_{\beta-e_i,k,n},$ and $\psi_{\beta-e_i,k,n}$ are chosen
	in accordance to Theorem~\ref{thm:approx}. Now $\partial_\tau^{e_i}(u_n,v_n)$ belongs to $\cG^1_\Gamma$, and one directly derives the
	energy equality for the approximation
	\begin{align}\label{eq:differentiated-energy-est-n}
		       & \frac12\!\int_{\cN^\Gamma_a}\!\!\big[\epsilond\partial_\tau^{e_i}u_n(t)\cdot\partial_\tau^{e_i}u_n(t)
			+ \mud\partial_\tau^{e_i}v_n(t)\cdot\partial_\tau^{e_i}v_n(t)\big]\mathrm{d} x
		+ \int_s^t\!\!\int_{\cN^\Gamma_a}\!\!\!\sigma\partial_\tau^{e_i}u_n\cdot\partial_\tau^{e_i}u_n\dd x\dd\tau \notag                              \\
		       & = \frac12\int_{\cN^\Gamma_a}\big[\epsilond\partial_\tau^{e_i}u_n(s)\cdot\partial_\tau^{e_i}u_n(s)
		+ \mud\partial_\tau^{e_i}v_n(s)\cdot\partial_\tau^{e_i}v_n(s)\big]\dd x                                                                        \\
		\notag & \seq + \frac12\int_s^t\!\!\int_{\cN^\Gamma_a} \big[(\pt\epsilond)\partial_\tau^{e_i}u_n\cdot\partial_\tau^{e_i}u_n
		+ (\pt\mud)\partial_\tau^{e_i}v_n\cdot\partial_\tau^{e_i}v_n\big]\dd x\dd\tau                                                                  \\
		\notag & \seq -\int_s^t\!\!\int_{\cN^\Gamma_a} \big[\partial_\tau^{e_i}u_n\cdot\big(\partial_\tau^{e_i}\frakf_{\beta-e_i,k,n}
		- [\partial_\tau^{e_i},\curl]v_n + (\partial_\tau^{e_i}\epsilond)\pt u_n + (\partial_\tau^{e_i}\sigma)u_n\big)                                 \\
		\notag & \seq\qquad\quad + \partial_\tau^{e_i}v_n\cdot\left(\partial_\tau^{e_i}\frakg_{\beta-e_i,k,n} + [\partial_\tau^{e_i},\curl]u_n
		+ (\partial_\tau^{e_i}\mud)\pt v_n\right)\big]\dd x\dd\tau                                                                                     \\
		\notag & \seq + \int_s^t\!\!\oint_\Gamma\left(\partial_\tau^{e_i}\varphi_{\beta-e_i,k,n}\cdot\tr_{\tau,\partial\upsilon}\partial_\tau^{e_i}v_n
		-\tr_{\tau,\partial\omega} \partial_\tau^{e_i} u_n\cdot\partial_\tau^{e_i}\psi_{\beta-e_i,k,n}\right)\dd\varsigma\dd\tau.
	\end{align}
	The convergence of $\partial_\tau^{e_i}(u_n,v_n) \to \partial_\tau^\beta(u,v)$ in $L^2(\cN_a^\Gamma)$ shows that
	\begin{align*}
		\int_{\cN^\Gamma_a} \epsilond \partial_\tau^{e_i} u_n(r) \cdot \partial_\tau^{e_i} u_n(r) \dd x
		 & \to \int_{\cN^\Gamma_a} \epsilond \partial_\tau^\beta \pt^k (\chi_\Gamma E(r)) \cdot \partial_\tau^\beta \pt^k (\chi_\Gamma E(r))\dd x, \\
		\int_{\cN^\Gamma_a} \mud \partial_\tau^{e_i} v_n(r) \cdot \partial_\tau^{e_i} v_n(r) \dd x
		 & \to \int_{\cN^\Gamma_a} \mud \partial_\tau^\beta \pt^k (\chi_\Gamma H(r)) \cdot \partial_\tau^\beta \pt^k (\chi_\Gamma H(r)) \dd x
	\end{align*}
	for $r \in \{s,t\}$, as $n\to\infty$. Similarly, the terms involving $\sigma, \pt\epsilond$ and $\pt\mud$ converge to the
	expressions without $n$, and
	likewise the inhomogeneities and commutator terms. To reconstruct $\frakf_{\beta, k}$ from the limit of the fourth line of
	\eqref{eq:differentiated-energy-est-n}, we use Equation~\eqref{pde:tangential-derivatives} and the commutator identity
	$[AB, C] = A[B, C] + [A, C]B$ to calculate
	\begin{align*}
		 & \partial_\tau^{e_i}\frakf_{\beta-e_i,k} - [\partial_\tau^{e_i},\curl]v + (\partial_\tau^{e_i}\epsilond)\pt u + (\partial_\tau^{e_i}\sigma)u \\
		 & = \partial_\tau^{e_i} \Big( {-} \partial_\tau^{\beta-e_i}\big([\chi_\Gamma,\curl]\pt^k H\big) - [\partial_\tau^{\beta-e_i},\curl]
		(\chi_\Gamma \pt^k H) + \partial_\tau^{\beta-e_i} (\chi_\Gamma\tilde f_k)                                                                      \\
		 & \qquad\quad + \sum_{0\leq b <\beta-e_i}\!\!\tbinom{\beta-e_i}{b}\big(\partial_\tau^{\beta-e_i - b}\sigma\partial_\tau^b\pt^{k}
		(\chi_\Gamma E)+ \partial_\tau^{\beta-e_i - b}\epsilond\partial_\tau^b\pt^{k+1}(\chi_\Gamma E)\big) \Big)                                      \\
		 & \quad - [\partial_\tau^{e_i},\curl]v + [\partial_\tau^{e_i}, \epsilond]\pt u + [\partial_\tau^{e_i}, \sigma]u                               \\
		 & = - \partial_\tau^\beta\Big([\chi_\Gamma,\curl]\pt^k H\Big) + \partial_\tau^\beta (\chi_\Gamma\tilde f_k)
		- \Big(\partial_\tau^{e_i}[\partial_\tau^{\beta-e_i},\curl](\chi_\Gamma \pt^k H) + [\partial_\tau^{e_i},\curl]v\Big)                           \\
		 & \quad + \Big(\partial_\tau^{e_i}[\partial_\tau^{\beta-e_i}, \sigma]\pt^k(\chi_\Gamma E)
		+ [\partial_\tau^{e_i}, \sigma]\partial_\tau^{\beta-e_i}\pt^k(\chi_\Gamma E)\Big)                                                              \\
		 & \quad + \Big(\partial_\tau^{e_i}[\partial_\tau^{\beta-e_i}, \epsilond]\pt^{k+1}(\chi_\Gamma E)
		+ [\partial_\tau^{e_i}, \epsilond]\partial_\tau^{\beta-e_i}\pt^{k+1}(\chi_\Gamma E)\Big)                                                       \\
		 & = - \partial_\tau^\beta\big([\chi_\Gamma,\curl]\pt^k H\big) - [\partial_\tau^\beta,\curl](\chi_\Gamma\pt^k H)
		+ \partial_\tau^\beta (\chi_\Gamma\tilde f_k) + [\partial_\tau^\beta, \sigma]\pt^k(\chi_\Gamma E)                                              \\
		 & \quad + [\partial_\tau^\beta, \epsilond]\pt^{k+1}(\chi_\Gamma E)                                                                            \\
		 & = \frakf_{\beta,k}.
	\end{align*}
	Likewise, one can show the convergence of
	\begin{equation*}
		\partial_\tau^{e_i}\frakg_{\beta-e_i,k,n} + [\partial_\tau^{e_i},\curl]u_n + (\partial_\tau^{e_i}\mud)\pt v_n
		\ \longrightarrow \ \frakg_{\beta,k},
	\end{equation*}
	in $L^2((s,t), L^2(\cN^\Gamma_a))$.
	It remains to consider the interface terms in \eqref{eq:differentiated-energy-est-n}.
	Using the continuity of $\partial_\tau^{e_i} \colon H^{1/2}(\Gamma) \to H^{-1/2}(\Gamma)$, see Chapter~2.3 in \cite{cessenat1996},
	we see that
	\begin{align}\label{eq:differentiated-interface-est-n}
		\int_s^t & \!\!\oint_\Gamma\left(\partial_\tau^{e_i}\varphi_{\beta-e_i,k,n}\cdot\tr_{\tau,\partial\upsilon}\partial_\tau^{e_i}v_n
		- \varphi_{\beta,k}\cdot\tr_{\tau,\partial\upsilon}\partial_\tau^\beta\pt^k (\chi_\Gamma H)\right)\dd\varsigma\dd\tau             \\*
		\notag   & \lesssim \int_s^t \Big( \norm{\partial_\tau^{e_i}\varphi_{\beta-e_i,k,n}
		- \varphi_{\beta,k}}_{H^{1/2}(\Gamma)} \norm{\tr_{\tau,\partial\upsilon}\partial_\tau^{e_i}v_n}_{H^{-1/2}(\Gamma)}                \\
		\notag   & \qquad\quad + \norm{\varphi_{\beta,k}}_{H^{1/2}(\Gamma)} \norm{\tr_{\tau,\partial\upsilon}\left(\partial_\tau^{e_i}v_n
		- \partial_\tau^\beta\pt^k (\chi_\Gamma H)\right)}_{H^{-1/2}(\Gamma)} \Big) \dd\tau                                               \\
		\notag   & \lesssim \int_s^t \Big( \norm{\partial_\tau^{e_i}\varphi_{\beta-e_i,k,n}
		- \varphi_{\beta,k}}_{H^{1/2}(\Gamma)} \norm{\tr_{\tau,\partial\upsilon}\partial_\tau^{e_i}v_n}_{H^{-1/2}(\Gamma)}                \\
		\notag   & \qquad\quad + \norm{\varphi_{\beta,k}}_{H^{1/2}(\Gamma)} \cdot \big[ \norm{\tr_{\tau,\partial\upsilon}\left(v_n
		- \partial_\tau^{\beta-e_i}\pt^k (\chi_\Gamma H)\right)}_{H^{1/2}(\Gamma)}                                                        \\
		\notag   & \qquad\quad+ \norm{[\tr_{\tau,\partial\upsilon},\partial_\tau^{e_i}]\left(v_n
		- \partial_\tau^{\beta-e_i}\pt^k (\chi_\Gamma H)\right)}_{H^{-1/2}(\Gamma)}\big]\Big) \dd\tau                                     \\
		\notag   & \lesssim \int_s^t \Big( \norm{\partial_\tau^{e_i}\varphi_{\beta-e_i,k,n}
		- \varphi_{\beta,k}}_{H^{1/2}(\Gamma)} \norm{\tr_{\tau,\partial\upsilon}\partial_\tau^{e_i}v_n}_{H^{-1/2}(\Gamma)}                \\
		\notag   & \qquad\quad + \norm{\varphi_{\beta,k}}_{H^{1/2}(\Gamma)} \cdot \norm{v_n
				- \partial_\tau^{\beta-e_i}\pt^k (\chi_\Gamma H)}_{H^1_\Gamma} \Big) \dd\tau,
	\end{align}
	where the commutator $[\tr_{\tau,\partial\upsilon},\partial_\tau^{e_i}]$ is of order $0$.
	The first term on the right vanishes because $\partial_\tau^{e_i}\varphi_{\beta-e_i,k,n} \to \varphi_{\beta,k}$ in
	$L^2\bigl((s,t),H^{1/2}(\Gamma)\bigr)$ and $\tr_{\tau,\partial\upsilon}\partial_\tau^{e_i}v_n$ is bounded in
	$L^2\bigl((s,t), H^{-1/2}(\Gamma)\bigr)$ by the comments after \eqref{pde:pt-maxwell_BC}.
	The second term tends to 0 since $v_n \to \partial_\tau^{\beta-e_i}\pt^k (\chi_\Gamma H)$ in
	$\cG^1_\Gamma$. Thus the difference in \eqref{eq:differentiated-interface-est-n} converges to $0$.
	The same argument applies to the second interface term, leading to \eqref{eq:differentiated-energy-est}.

	\smallskip

	1.3) As in \eqref{eq:differentiated-interface-est-n}, the trace integrals from \eqref{eq:differentiated-energy-est} can be estimated by
	\begin{align*}
		\Big|\int_s^t\!\! & \oint_\Gamma\left(\varphi_{\beta,k}\cdot\tr_{\tau,\partial\upsilon}\partial_\tau^\beta\pt^k (\chi_\Gamma H)
		-\tr_{\tau,\partial\omega}\partial_\tau^\beta\pt^k (\chi_\Gamma E)\cdot\psi_{\beta,k}\right)\dd\varsigma\dd\tau\Big|                              \\
		                  & \le \int_s^t\Big(\norm{\varphi_{\beta,k}}_{H^{1/2}(\Gamma)}\norm{\tr_{\tau,\partial\upsilon}\partial_\tau^\beta\pt^k
		(\chi_\Gamma H)}_{H^{-1/2}(\Gamma)}                                                                                                               \\*
		                  & \hphantom{\leq\int_s^t\Big(} + \norm{\psi_{\beta,k}}_{H^{1/2}(\Gamma)}\norm{\tr_{\tau,\partial\omega}\partial_\tau^\beta\pt^k
		(\chi_\Gamma E)}_{H^{-1/2}(\Gamma)}\Big)\dd\tau                                                                                                   \\
		                  & \lesssim \theta\!\!\int_s^t\norm{\pt^k E}^2_{H_\Gamma^{\abs\beta}}\dd\tau
		+ c_\theta\!\!\int_s^t\norm{\pt^k H}^2_{H_\Gamma^{\abs\beta}}\dd\tau
		\lesssim\theta\!\!\int_s^t \!\!z(\tau)\dd\tau + c_\theta\!\!\int_s^t\norm{\pt^k H}^2_{H_\Gamma^{\abs\beta}}\dd\tau,
	\end{align*}
	using that $\varphi_{\beta,k}$ and $\psi_{\beta,k}$ only contain lower-order terms. Furthermore, the inhomogeneities satisfy
	\begin{align*}
		\norm{\frakf_{\beta,k}}_{L^2} & \lesssim\norm{\pt^k H}_{H_\Gamma^{\abs{\beta}}}
		+ \norm{\pt^k E}_{H_\Gamma^{\abs{\beta}-1}}+ \norm{\pt^{k+1} E}_{H^{\abs{\beta}-1}}
		+ \bignorm{\tilde f_k}_{H_\Gamma^{\abs{\beta}}},                                \\
		\norm{\frakg_{\beta,k}}_{L^2} & \lesssim\norm{\pt^k E}_{H_\Gamma^{\abs{\beta}}}
		+ \norm{\pt^{k+1} H}_{H_\Gamma^{\abs{\beta}-1}} + \bignorm{\tilde g_k}_{H_\Gamma^{\abs{\beta}}}.
	\end{align*}
	Employing also Lemma~\ref{lem:error-estimate}, from \eqref{eq:differentiated-energy-est} we thus deduce the inequality
	\begin{align}
		\big\|\partial_\tau^\beta\pt^k & (\chi_\Gamma E,\chi_\Gamma H)(t)\big\|^2_{L^2(\Omega)}
		+ \sigma_0\int_s^t\norm{\partial_\tau^\beta\pt^k (\chi_\Gamma E)}^2_{L^2(\omega)}\dd\tau \label{eq:proof-tangential-estimate}                                        \\
		                               & \lesssim z(s) + c_\theta\int_s^t z^{3/2}(\tau)\dd\tau + \theta\int_s^t z(\tau)\dd\tau \notag                                        \\
		                               & \seq + c_\theta\int_s^t\left(\norm{\pt^k E}^2_{H_\Gamma^{\abs{\beta}-1}}+ \norm{\pt^{k+1} (E,H)}^2_{H_\Gamma^{\abs{\beta}-1}}\notag
		+ \norm{\pt^k H}^2_{H_\Gamma^{\abs{\beta}}}\right)\dd\tau.
	\end{align}
	We emphasize that the integral-term on the left-hand side only contains the $L^2$-norm on the support $\omega$ of $\sigma$.
	Lemma~\ref{lem:H-est} will allow us to handle the highest-order $H$ term on the right and to add the missing
	contribution of $H$ on the left.

	\smallskip

	1.4) To estimate the full Sobolev norm on $\omega$, we have to consider normal and mixed derivatives of
	$u\coloneqq\partial_\tau^\beta\pt^k(\chi_\Gamma E)$. Let
	$\alpha=(\gamma,\beta)\in\N_0\times\N_0^2$ with $\gamma\neq0$ and $\abs\alpha\leq3-k$. Using the components
	\begin{equation*}
		u_\xi= (u)_\xi = (u\cdot\xi)\xi \quad\text{for }\xi\in\{\tau_1,\tau_2,\nu\}\quad \text{and}\quad
		u_\tau=(u)_{\tau}= (u)_{\tau_1} + (u)_{\tau_2},
	\end{equation*}
	we split $\partial_\nu^\gamma u$ into
	\[\norm{\partial_\nu^\gamma\partial_\tau^\beta\pt^k (\chi_\Gamma E)}_{L^2(\omega)}
		\lesssim\norm{\partial_\nu^\gamma\left(\partial_\tau^\beta\pt^k (\chi_\Gamma E)\right)_{\nu}}_{L^2(\omega)}
		+ \norm{\partial_\nu^\gamma\left(\partial_\tau^\beta\pt^k (\chi_\Gamma E)\right)_{\tau}}_{L^2(\omega)}. \]
	By means of formula~(6.1) in \cite{lasiecka2019}, the tangential part $\partial_\nu u_{\tau}$
	can be recovered from the curl and tangential derivatives up to lower-order terms via
	\[ \partial_\nu u_{\tau} = R(\nu^\Gamma)\Big(\curl u -\sum_{i=1,2} J(\tau_i)\partial_{\tau_i} u\Big)
		+\sum_{i=1,2}\big(\partial_\nu\tau_i (u)_{\tau_i} + \tau_i\partial_\nu\tau_i\cdot u\big) \]
	for matrices $R(\nu)$ and $J(\tau_i)$ in $C^4(\overline\Omega)$. We thus obtain
	\begin{align*}
		\norm{\partial_\nu^\gamma u_{\tau}}_{L^2} & = \norm{\partial_\nu^{\gamma-1}\partial_\nu u_{\tau}}_{L^2}
		\!\lesssim \|u\|_{H^{|\gamma|-1}_\Gamma} + \norm{\partial_\nu^{\gamma-1}\curl u}_{L^2} + \sum_{i=1,2}
		\norm{\partial_\nu^{\gamma-1}\partial_\tau^{\beta+e_i}u}_{L^2}\!.
	\end{align*}
	We insert Equation~\eqref{pde:tangential-derivatives} into the curl term, which gives $y_{\partial\Omega}^{k+1,\abs\alpha-1}$ and
	\[ \norm{\frakg_{k,b}}^2_{H^{\abs\gamma-1}_\Gamma}
		\lesssim\norm{\pt^k E}_{H^{\abs\alpha-1}_\Gamma}^2 + y_{\partial\Omega}^{k+1,\abs\alpha-1} + z^2.\]
	We split and estimate $\partial_\nu^{\gamma-1}\partial_\tau^{\beta+e_i}u$ in the same way if $\gamma=2$ and iterate once more if
	$\gamma=3,$ obtaining
	purely tangential derivatives and terms as above. Keeping the lower-order normal parts, we conclude
	\begin{align}
		\norm{\partial_\nu^\gamma u_{\tau}}^2_{L^2(\omega)}
		 & \lesssim \norm{\pt^k E}_{H^{\abs\alpha-1}_\Gamma}^2 +y_{\Gamma}^{k+1,\abs\alpha -1}
		+ z^2+ \sum_{{\abs b\leq\abs\alpha}}\norm{\partial_\tau^{b}\pt^k(\chi_{\Gamma} E)}^2_{L^2(\omega)}\notag \\
		 & \seq + \sum_{{c\le \gamma-1, \abs b\leq\abs\alpha -c}}
		\norm{\partial_\nu^c\big(\partial_\tau^{b}\pt^k(\chi_{\Gamma} E)\big)_\nu}^2_{L^2(\omega)}.\label{eq:proof-normal-derivative-of-tangential}
	\end{align}

	1.5) The normal part $\partial_\nu u_\nu$ on $\omega$ is controlled through divergence terms.
	Formula~(6.3) of \cite{lasiecka2019} yields
	\begin{align*}
		(\epsilonk)_{\nu\nu}\partial_\nu u_\nu
		 & = \div(\epsilonk u) -\!\sum_{(\xi,\zeta)\neq (\nu^\Gamma,\nu^\Gamma)}\!\!\! (\epsilonk)_{\xi\zeta}\partial_\xi u_\zeta
		-\!\sum_{\xi,\zeta}\partial_\xi(\epsilonk)_{\xi\zeta} u_\zeta -\sum_\xi\div(\xi)\xi^\top (\epsilonk u)                    \\
		 & \eqcolon\div\left(\epsilonk u\right) - D(\epsilonk)u,
	\end{align*}
	where $\xi,\zeta$ run over $\{\nu,\tau_1,\tau_2\}$ and $\varkappa_{\xi,\zeta}\coloneqq\varkappa\xi\cdot\zeta$.
	We take the tangential derivative and divergence of~\eqref{pde:pt-maxwell} and infer
	\begin{align*}
		\partial_t (\epsilonk u) & = \curl\partial_\tau^\beta\pt^k (\chi_{\Gamma} H) -\sigma u -[\curl,\partial_\tau^\beta] (\chi_{\Gamma}\pt^kH)
		-\partial_\tau^\beta ([\curl,\chi_{\Gamma}]\pt^k H)\notag                                                                                                \\
		                         & \seq -\!\!\sum_{0\le b<\beta}\!\!\tbinom{\beta}{b} \big[\partial_\tau^{\beta - b}\sigma\partial_\tau^b\pt^k (\chi_{\Gamma} E)
			+ \pt\!\left(\partial_\tau^{\beta - b}\epsilonk\partial_\tau^b\pt^k (\chi_{\Gamma} E)\right)\!\big]
		-\partial_\tau^\beta(\chi_{\Gamma}\pt f_k)\notag                                                                                                         \\
		                         & \eqqcolon\curl\partial_\tau^\beta\pt^k (\chi_{\Gamma} H) -\sigma u -\varphi,                                                  \\
		\div(\epsilonk u)(t)     & = \div(\epsilonk u)(s) -\int_s^t\big(\div(\sigma u(\tau)) + \div\varphi(\tau)\big)\dd\tau.\notag
	\end{align*}
	These formulas lead to the integral equation
	\[ (\epsilonk)_{\nu\nu}\partial_\nu u_\nu = \div(\epsilonk u)(s)- D(\epsilonk)u -\int_s^t\big(\pi (\epsilonk)_{\nu\nu}\partial_\nu u_\nu
		+ D(\sigma) u+ \div\varphi\big)(\tau)\dd\tau\]
	with $\pi \coloneq \frac{\sigma_{\nu\nu}}{(\epsilonk)_{\nu\nu}}\geq\frac{\sigma_0}{\eta}\eqqcolon\pi_0>0$ on $\omega$
	by (H), but $\pi = 0$ and $D(\sigma) =0$ on $\upsilon$. This equation is solved by
	\begin{equation*}
		(\epsilonk)_{\nu\nu}\partial_\nu u_\nu(t)
		= \e^{-\int_s^t\pi(r)\dd r}\div(\epsilonk u)(s) + k(t) -\int_s^t\!\!\big(\pi(\tau) k(\tau)
		+ \ell(\tau)\big)\e^{-\int_\tau^t\pi(r)\dd r}\dd\tau,
	\end{equation*}
	abbreviating $k\coloneqq - D(\epsilonk)u$ and $\ell\coloneqq D(\sigma) u+ \div\varphi$. We apply $\partial_\nu^{\gamma-1}$ to
	this identity. Using Lemmas~\ref{lem:error-estimate} and \ref{lem:sobolev-algebra-property} and the chain rule,
	the first term is estimated by
	\begin{align*}
		\norm{\partial_\nu^{\gamma-1}\!\big(\e^{-\int_s^t\pi(r)\dd r}\!\div(\epsilonk u)(s)\big)}_{L^2}^2\!
		 & \lesssim \e^{-\pi_0(t-s)}\big[\norm{\partial_\nu^{\gamma-1}\!\div(\epsilonk u)(s)}_{L^2}^2\!+\norm{u(s)}_{H_\Gamma^{\gamma-1}}^2\!\big]\!.
	\end{align*}
	The integrand is treated analogously. As $(\epsilonk)_{\nu\nu}\ge \eta$, H\"older's inequality yields
	\begin{align}\label{proof:eq:nu-nu1}
		\|     & \partial_\nu^\gamma u_\nu(t)\|_{L^2(\omega)}^2
		\lesssim\norm{\partial_\nu^{\gamma-1}\big((\epsilonk)_{\nu\nu}\partial_\nu u_\nu(t)\big)}_{L^2(\omega)}^2
		+ \norm{u(t)}_{H_\Gamma^{\gamma-1}}^2                                                                                                    \\
		\notag & \lesssim\e^{-\pi_0 (t-s)}\big(\norm{\partial_\nu^{\gamma-1}\!\div(\epsilonk u)(s)}_{L^2}^2
		+ \norm{u(s)}_{H_\Gamma^{\gamma -1}}^2\big) + \norm{\partial_\nu^{\gamma-1} k(t)}_{L^2}^2 + \norm{u(t)}_{H_\Gamma^{\gamma -1}}^2         \\
		\notag & \quad+\Big[\int_s^t\!\!\!\e^{-\frac12\pi_0(t-\tau)}\e^{-\frac12\pi_0(t-\tau)} \big(\norm{\partial_\nu^{\gamma-1} k(\tau)}_{L^2}
		\! + \norm{\partial_\nu^{\gamma-1}\ell(\tau)}_{L^2}\! + \norm{u(\tau)}_{H_\Gamma^{\gamma -1}}\!\big)\mathrm{d}\tau\Big]^2                \\
		\notag & \lesssim\e^{-\pi_0 (t-s)}\big(\norm{\partial_\nu^{\gamma-1}\div(\epsilonk u)(s)}_{L^2}^2
		+ \norm{u(s)}_{H_\Gamma^{\gamma -1}}^2\big) + \norm{\partial_\nu^{\gamma-1} k(t)}_{L^2}^2 + \norm{u(t)}_{H_\Gamma^{\gamma -1}}^2         \\
		       & \quad + \frac{1}{\pi_0}\int_s^t\e^{-\pi_0(t-\tau)} \big(\norm{\partial_\nu^{\gamma-1} k(\tau)}_{L^2}^2
		+ \norm{\partial_\nu^{\gamma-1} \ell(\tau)}_{L^2}^2 + \norm{u(\tau)}_{H_\Gamma^{\gamma -1}}^2\big) \dd\tau.\notag
	\end{align}
	Young's inequality then leads to
	\begin{align} \label{proof:eq:nu-nu2}
		\int_s^t\|\partial_\nu^\gamma u_\nu & (\tau)\|_{L^2(\omega)}^2\dd\tau
		\lesssim\frac{1}{\pi_0}\big(\norm{\partial_\nu^{\gamma-1}\div(\epsilonk u)(s)}_{L^2}^2 + \norm{u(s)}_{H_\Gamma^{\gamma -1}}^2\big)                   \\*
		\notag                              & \seq + \Big(1+\frac{1}{\pi_0^2}\Big)\! \int_s^t\!\!\big(\norm{\partial_\nu^{\gamma-1} k(\tau)}_{L^2(\omega)}^2
		\!+ \norm{\partial_\nu^{\gamma-1} \ell(\tau)}_{L^2(\omega)}^2 \!+ \norm{u(\tau)}_{H_\Gamma^{\gamma -1}}^2\!\big)\mathrm{d}\tau.
	\end{align}
	Note that $\norm{\partial_\nu^{\gamma-1}\div(\epsilonk u)(s)}_{L^2}^2 \lesssim z(s)$. As above
	Lemmas~\ref{lem:error-estimate} and \ref{lem:sobolev-algebra-property} imply
	\begin{align}\label{proof:eq:nu-nu3}
		\norm{\partial_\nu^{\gamma-1} k}_{L^2(\omega)}^2
		       & \lesssim\Big\|\partial_\nu^{\gamma-1}\Big(\sum_{(\xi,\zeta)\neq(\nu,\nu)}\partial_\xi u_\zeta\Big)\Big\|_{L^2(\omega)}^2
		+ \norm{u}_{H^{\gamma-1}_\Gamma}^2                                                                                                \\
		\notag & \lesssim\norm{\partial_\nu^{\gamma}u_\tau}_{L^2(\omega)}^2
		+ \sum_{i=1,2}\norm{\partial_\nu^{\gamma-1}(\partial^{e_i}_\tau u)_\nu}_{L^2(\omega)}^2 + y^{k,\abs\alpha-1}_{\Gamma}.
	\end{align}
	Here the first term is treated by \eqref{eq:proof-normal-derivative-of-tangential} and the middle terms have fewer
	normal derivatives, so they can be handled by iteration below. Similarly, we compute
	\begin{align}
		\label{proof:eq:nu-nu4}\norm{\partial_\nu^{\gamma -1} \ell}_{L^2(\omega)}^2
		       & \lesssim\Big\|\partial_\nu^{\gamma-1}\Big(\sum_{(\xi,\zeta)\neq (\nu,\nu)}\partial_\xi u_\zeta\Big)\Big\|_{L^2(\omega)}^2
		\!+ \norm{u(\tau)}_{H^{\gamma-1}_{\Gamma}}^2 \!+ \norm{\partial_\nu^{\gamma} \varphi}_{L^2}^2                                      \\
		\notag & \lesssim\norm{\partial_\nu^{\gamma} u_\tau}_{L^2(\omega)}^2
		+ \sum_{i=1,2}\norm{\partial_\nu^{\gamma-1}(\partial^{e_i}_\tau u)_\nu}_{L^2(\omega)}^2
		\!+ y^{k,\abs\alpha-1}_{\Gamma}+ \!\norm{\varphi}_{H^{\gamma}_\Gamma}^2,                                                           \\
		\notag \norm{\varphi}_{H^{\gamma}_\Gamma}^2
		       & \lesssim\norm{\pt f_k}_{H^{3-k}_\Gamma}^2 + \norm{\chi_\Gamma\pt^k H}_{H^{\abs\alpha}_\Gamma}^2
		+ \norm{\chi_\Gamma\pt^k E}_{H^{\abs\alpha-1}_\Gamma}^2\! + \norm{\chi_\Gamma\pt^{k+1} E}_{H^{\abs\alpha-1}_\Gamma}^2              \\
		       & \lesssim z^2 + \norm{\pt^k H}_{H^{\abs\alpha}_\Gamma}^2 + y^{k,\abs\alpha-1}_{\Gamma}+ y^{k+1,\abs\alpha-1}_{\Gamma}.
		\label{proof:eq:nu-nu5}
	\end{align}

	We now combine inequalities \eqref{proof:eq:nu-nu1}--\eqref{proof:eq:nu-nu5} and iterate on the terms
	with normal components on the right, resulting in
	\begin{align}\label{proof:eq:nu-nu-interface}
		       & \norm{\partial_\nu^\gamma\left(\partial_\tau^\beta\pt^k(\chi_{\Gamma} E)\right)_\nu(t)}_{L^2(\omega)}^2 + \int_s^t\norm{\partial_\nu^\gamma\left(\partial_\tau^\beta\pt^k(\chi_{\Gamma} E)\right)_\nu}_{L^2(\omega)}^2\dd\tau \\
		\notag & \lesssim \sum_{m=0}^\gamma \!\!\sum_{\substack{b\in \N_0^2                                                                                                                                                                    \\ \abs b \leq \gamma-m}} \!\!\!\!\left[
			\norm{\partial_\nu^m \big(\partial_\tau^{\beta+b}\pt^k(\chi_{\Gamma} E)\big)_\tau(t)}_{L^2(\omega)}^2
		+ \int_s^t\norm{\partial_\nu^m \big(\partial_\tau^{\beta+b}\pt^k(\chi_{\Gamma} E)\big)_\tau}_{L^2(\omega)}^2\dd\tau \right]                                                                                                            \\
		\notag & \seq + z(s)+ z^2(t) + \norm{\pt^k H}_{H^{\abs\alpha}_\Gamma}^2 + y^{k,\abs\alpha-1}_{\Gamma}                                                                                                                                  \\ \notag
		       & \seq + \int_s^t\big( z^2+ \norm{\pt^k H}_{H^{\abs\alpha}_\Gamma}^2 + y^{k,\abs\alpha-1}_{\Gamma}+ y^{k+1,\abs\alpha-1}_{\Gamma}\big)\dd\tau.
	\end{align}
	Estimates~\eqref{eq:proof-tangential-estimate},~\eqref{eq:proof-normal-derivative-of-tangential} and \eqref{proof:eq:nu-nu-interface}
	then imply
	\begin{align}\label{eq:reg-ana-est-on-omega}
		       & \norm{\pt^k(\chi_{\Gamma}E)(t)}^2_{H^{\abs\alpha}(\omega)} + \int_s^t\norm{\pt^k(\chi_{\Gamma}E)}^2_{H^{\abs\alpha}(\omega)}\dd\tau \\
		\notag & \lesssim z(s) +z^2(t)+ \theta \!\int_s^t\!\! z(\tau)\dd\tau+ c_\theta\bigg[\norm{\pt^{k}E(t)}_{H^{\abs\alpha-1}_\Gamma}^2
		+\norm{\pt^k H(t)}^2_{H^{\abs\alpha}_\Gamma} + y^{k+1,\abs\alpha-1}_{\Gamma}(t)\notag                                                        \\
		       & \seq \int_s^t\Big( z^{\frac32}(\tau)+\norm{\pt^{k}E}_{H^{\abs\alpha-1}_\Gamma}^2+ \norm{\pt^k H}^2_{H^{\abs\alpha}_\Gamma}
			+ y^{k+1,\abs\alpha-1}_{\Gamma} \Big)\dd\tau \bigg].\notag
	\end{align}

	\smallskip

	2) In view of Lemma~\ref{lem:H-est}, it remains to control the left-hand side of \eqref{eq:reg-ana-est-on-omega} on $\upsilon$.
	Here we use the standard div-curl estimate
	\begin{align}\label{proof:div-curl}
		\norm{\pt^k(\chi_\Gamma E)}^2_{H^{\abs\alpha}(\upsilon)}
		       & \lesssim\norm{\pt^k(\chi_\Gamma E)}^2_{L^2(\upsilon)} + \norm{\curl\pt^k(\chi_\Gamma E)}^2_{H^{\abs\alpha-1}(\upsilon)}                                                \\
		\notag & \seq + \norm{\div_{\epsilonk}\pt^k(\chi_\Gamma E)}^2_{H^{\abs\alpha-1}(\upsilon)}+ \norm{\tr_{t,\partial\upsilon}\pt^k(\chi_\Gamma E)}^2_{H^{\abs\alpha-1/2}(\Gamma)}.
	\end{align}
	see e.g.\ Corollary~3.5 in \cite{amrouche2021}. Note that $\pt^k(\chi_\Gamma E)$ vanishes at $\partial\cN_a^\Gamma\cap\upsilon$.
	The Maxwell system~\eqref{eq:maxwell} implies
	\[\curl (\chi_\Gamma E) = \nabla\chi_\Gamma\times E-\chi_\Gamma\pt (\mu H), \quad
		\div_\epsilon(\chi_\Gamma E) = \nabla\chi_\Gamma\cdot\epsilon E
		\qquad \text{on $\upsilon$,}\]
	since $\div(\epsilon(E_0)E_0)|_\upsilon = 0$ and $\sigma|_\upsilon =0$ by (H).
	For $k\in\{1,2\}$ Equation~\eqref{pde:tilde1} yields
	\begin{align*}
		\begin{split}
			\curl\pt^k (\chi_\Gamma E)            & = -\mud\pt^{k+1}(\chi_\Gamma H) - [\chi_\Gamma,\curl]\pt^k E -\chi_\Gamma\tilde g_k, \\
			\div_{\epsilond}\pt^k (\chi_\Gamma E) & = \nabla\chi_\Gamma\cdot\epsilond\pt^k E -\chi_\Gamma\div \tilde f_{k-1},
		\end{split}
		 & \text{on $\upsilon$,}
	\end{align*}
	where the commutator $[\chi_\Gamma,\curl]$ is a multiplication operator.
	We then obtain
	\begin{align*}
		\norm{\curl\pt^k (\chi_\Gamma E)}^2_{H^{\abs\alpha-1}(\upsilon)}            & \lesssim\norm{\pt^{k+1}H}^2_{H^{\abs\alpha-1}(\upsilon)}
		+ \norm{\pt^k E}_{H^{\abs\alpha-1}(\upsilon)} + z^{2}(t),                                                                                       \\
		\norm{\div_{\epsilonk}\pt^k (\chi_\Gamma E)}^2_{H^{\abs\alpha-1}(\upsilon)} & \lesssim\norm{\pt^k E}^2_{H^{\abs\alpha-1}(\upsilon)} + z^{2}(t).
	\end{align*}
	The interface condition \eqref{eq:maxwell-currents} yields
	\begin{align*}
		\norm{\tr_{t,\partial\upsilon}\pt^k(\chi_\Gamma E)}_{H^{\abs\alpha-1/2}(\Gamma)}
		= \norm{\tr_{t,\partial\omega}\pt^k(\chi_\Gamma E)}^2_{H^{\abs\alpha-1/2}(\Gamma)} 
		\lesssim\norm{\pt^k(\chi_\Gamma E)}^2_{H^{\abs\alpha}(\omega)}.
	\end{align*}
	From~\eqref{proof:div-curl} we thus infer
	\begin{align}\label{eq:estE}
		\norm{\pt^k(\chi_\Gamma E)}^2_{H^{\abs\alpha}(\upsilon)} & \lesssim e(t) + z^{2}(t) + \norm{\pt^k E}_{H_\Gamma^{\abs\alpha-1}}
		+ \norm{\pt^{k+1}H}^2_{H_\Gamma^{\abs\alpha-1}}                                                                                \\
		\notag                                                   & \seq + \norm{\pt^k(\chi_\Gamma E)}^2_{H^{\abs\alpha}(\omega)}.
	\end{align}
	Estimates~\eqref{eq:estE} and~\eqref{eq:reg-ana-est-on-omega} result in
	\begin{align*}
		 & \|\pt^k (\chi_\Gamma E)(t)\|^2_{H^{\abs\alpha}_\Gamma} + \int_s^t\norm{\pt^k (\chi_\Gamma E)}^2_{H_\Gamma^{\abs\alpha}}\dd\tau            \\
		 & \lesssim e(t)+z(s)+z^2(t)+\theta\!\int_s^t\!\!\!z\dd\tau + c_\theta\bigg[\! \norm{\pt^k E(t)}^2_{H_\Gamma^{\abs\alpha-1}}
		\! + \norm{\pt^{k+1} (E,H)(t)}^2_{H_\Gamma^{\abs\alpha-1}}                                                                                   \\
		 & + \!\norm{\pt^k H(t)}^2_{H_\Gamma^{\abs\alpha}}\! +\! \!\int_s^t\!\!\!\big[e+z^{\frac32}\!+\norm{\pt^k \!H}^2_{H_\Gamma^{\abs\alpha}} \!+
				\norm{\pt^k\! E}^2_{H_\Gamma^{\abs\alpha-1}} \! +\! \norm{\pt^{k+1} (E,H)}^2_{H_\Gamma^{\abs\alpha-1}}\!\big]\mathrm{d}\tau \bigg]\!.
	\end{align*}

	\smallskip

	3) In the above inequality Lemma~\ref{lem:H-est} allows to remove the $H$-terms in $H_\Gamma^{\abs\alpha}$ and to add $H$ terms on the left,
	leading to our result
	\begin{align*}
		 & y^{k,\abs\alpha}_{\Gamma}(t) + \int_s^t y^{k,\abs\alpha}_{\Gamma}\dd\tau                                                      \\
		 & \ \ \lesssim e(t)+z(s)+z^2(t)+\theta\!\int_s^t\!\!\!z\dd\tau + c_\theta\bigg[\! \norm{\pt^k E(t)}^2_{H_\Gamma^{\abs\alpha-1}}
		+ \norm{\pt^{k+1} (E,H)(t)}^2_{H_\Gamma^{\abs\alpha-1}}                                                                          \\
		 & \ \ \seq +\int_s^t\big(e+z^{\frac32}+ \norm{\pt^k E}^2_{H_\Gamma^{\abs\alpha-1}}
			+\norm{\pt^{k+1} (E,H)}^2_{H_\Gamma^{\abs\alpha-1}}\big)\dd\tau \bigg].\qedhere
	\end{align*}
\end{proof}

All that remains is to piece together the estimates above and use an iteration scheme as illustrated in Figure~\ref{fig:proof-structure}.
\begin{figure}[H]
	\centering
	\begin{tikzpicture}
		\matrix (m) [matrix of math nodes,
			row sep=2em, column sep=2em,
			nodes={anchor=center, minimum width=4em}]
		{
			(E,H) & \partial_t(E,H) & \partial_t^2(E,H) & \partial_t^3(E,H) \\
			(E,H) & \partial_t(E,H) & \partial_t^2(E,H) &                   \\
			(E,H) & \partial_t(E,H) &                   &                   \\
			(E,H) &                 &                   &                   \\
		};

		\node[left=1.5em of m-1-1] (L1) {$L^2$};
		\node[left=1.5em of m-2-1] (L2) {$H_\Gamma^1$};
		\node[left=1.5em of m-3-1] (L3) {$H_\Gamma^2$};
		\node[left=1.5em of m-4-1] (L4) {$H_\Gamma^3$};

		\draw[thick] ([xshift=-0.75em]m.north west) -- ([xshift=-0.75em]m.south west);

		\begin{scope}[on background layer]
			\node[fill=blue!15, rounded corners, fit=(m-1-1) (m-1-4)] (highlight) {};
		\end{scope}

		\node[right=1em of m-1-4, anchor=west, color=blue] (bound) {$\|\cdot\|_{L^2}^2\lesssim e(t)$};

		\def\arrowColor{black}
		\begin{scope}[every path/.style={->, >=stealth, draw=\arrowColor, shorten >=2pt, shorten <=2pt}]
			\draw (m-4-1) -- (m-3-1);

			\draw (m-3-1) -- (m-2-1);
			\draw (m-3-2) -- (m-2-2);

			\draw (m-2-1) -- (m-1-1);
			\draw (m-2-2) -- (m-1-2);
			\draw (m-2-3) -- (m-1-3);
		\end{scope}

		\begin{scope}[every path/.style={->, >=stealth, draw=\arrowColor, shorten >=3pt, shorten <=3pt}]
			\draw (m-4-1) -- (m-3-2);

			\draw (m-3-1) -- (m-2-2);
			\draw (m-3-2) -- (m-2-3);

			\draw (m-2-1) -- (m-1-2);
			\draw (m-2-2) -- (m-1-3);
			\draw (m-2-3) -- (m-1-4);
		\end{scope}

	\end{tikzpicture}
	\caption{We estimate the $H^j_\Gamma$-norm by terms in the $H^{j-1}_\Gamma$-norm, finally reaching the $L^2$-level. These terms are then bound by $e(t)$.}
	\label{fig:proof-structure}
\end{figure}
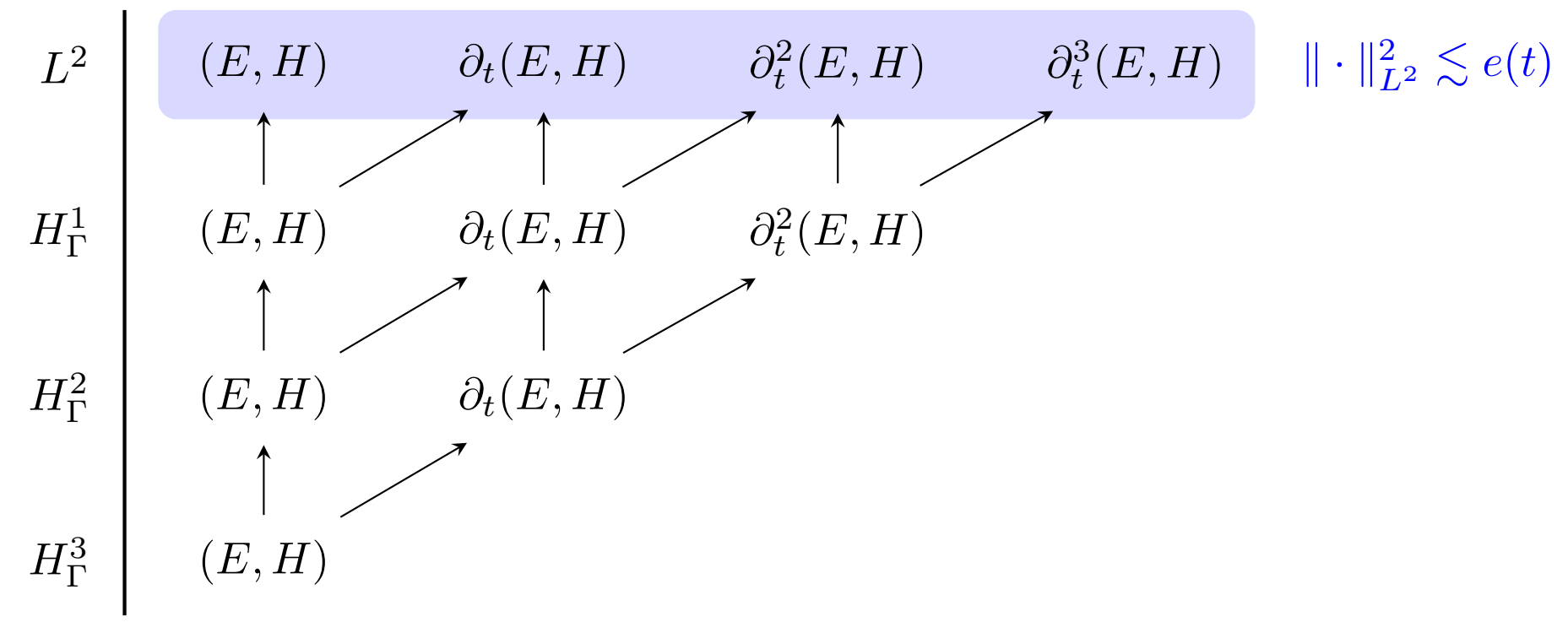
\begin{proof}[Proof of Theorem~\ref{thm:regularity-analysis}]
	For $k\in\N_0$, $j\in\N_0$ with $k +j\leq 2$ we have
	\[\norm{\pt^k (E,H)(t)}_{H_\Gamma^{j+1}}^2 + \int_s^t\!\!\norm{\pt^k (E,H)(\tau)}_{H_\Gamma^{j+1}}^2\dd\tau
		\lesssim\sum_{\Xi\in\{\upsilon,\Gamma,\omega,\partial\Omega\}}\!\! y^{k, j+1}_\Xi(t) + \int_s^t \! \!y^{k, j+1}_\Xi\dd\tau. \]
	By Lemmas~\ref{lem:upsilon-est} and \ref{lem:Gamma-est} the right-hand side can be estimated by
	\begin{align*}
		 & e(t)+z(s)+z^2(t)+\theta\!\int_s^t\!z\dd\tau + c_\theta\bigg[ \norm{\pt^k (E,H)(t)}^2_{H_\Gamma^{j-1}}
		+ \norm{\pt^{k+1} (E,H)(t)}^2_{H_\Gamma^{j-1}}                                                           \\
		 & \ +\int_s^t\big(e+z^{\frac32}+ \norm{\pt^k (E,H)}^2_{H_\Gamma^{j-1}}
			+\norm{\pt^{k+1} (E,H)}^2_{H_\Gamma^{j-1}}\big)\dd\tau \bigg].
	\end{align*}
	Iterating over $j$ and $k$, see Figure~\ref{fig:proof-structure}, we thus obtain
	\begin{align*}
		\norm{\pt^k (E,H)(t)}_{H_\Gamma^{j+1}} & + \int_s^t\norm{\pt^k (E,H)(\tau)}_{H_\Gamma^{j+1}}\dd\tau                                                             \\
		                                       & \lesssim e(t) +z(s) +z^2(t)+ \theta\int_s^t z(\tau)\dd\tau+ c_\theta\int_s^t\big(e(\tau) +z^\frac32(\tau)\big)\dd\tau.
	\end{align*}
	Since $z=e+\sum_{j+k\leq 2}\norm{\pt^k (E,H)}_{H_\Gamma^{j+1}}$ the result follows after absorbing $\theta\int z\dd \tau$ for a fixed small
	$\theta$.
\end{proof}

\section{Exponential decay}

We now combine the observability-type estimate of Theorem~\ref{thm:obs-est}, the regularity analysis of Theorem~\ref{thm:regularity-analysis},
and the energy inequality in Lemma~\ref{lem:energy-decreasing} to show the main result of this paper.

\begin{theorem}\label{thm:exp-decay}
	Assume that \ref{hypothesis-quasi} holds.
	Let the initial values $(E_0,H_0) \in \left(H_\Gamma^3(\Omega)\right)^6$ be sufficiently small.
	Then the corresponding solution $(E,H)$ to \eqref{eq:maxwell}--\eqref{eq:maxwell4} exists for all $t \geq 0$.
	Moreover, there exist constants $M \geq 1$ and $\omega>0$ such that
	\[
		z(t) \leq M \e^{-\omega t} z(0), \qquad t\geq 0.
	\]
\end{theorem}

\begin{proof}
	Theorem~\ref{thm:regularity-analysis}, Theorem~\ref{thm:obs-est}, and Lemma~\ref{lem:energy-decreasing} imply
	\begin{align*}
		z(t) + \int_s^t z(\tau) \dd \tau & \lesssim z(s) + z^2(t)+ e(t) + \int_s^t \left(e(\tau) + z^{3/2}(\tau)\right) \dd\tau                                                  \\
		                                 & \lesssim e(t) + c_\delta \left(z(s) + z^2(t) + \int_s^t (d(\tau) + z^{3/2}(\tau)) \dd \tau\right)                                     \\
		                                 & \seq+ \delta \left(z(t) + \int_s^t z(\tau) \dd \tau\right)                                                                            \\
		                                 & \lesssim c_\delta \left(z(s) + z^2(t) + \int_s^t z^{3/2}(\tau) \dd \tau\right)+ \delta \left(z(t) + \int_s^t z(\tau) \dd \tau\right),
	\end{align*}
	for $\delta \in (0,1]$.
	We first fix a sufficiently small $\delta>0$ to absorb the final terms by the left-hand side. Then we fix $\delta_z>0$ such that the
	superlinear terms can be put on the left as $z^{1/2}(t) \leq \delta_z$ for $t\in[0, T_*)$.
	So, we arrive at
	\begin{equation*}
		z(t) + \int_s^t z(\tau) \dd \tau \leq C z(s),
	\end{equation*}
	for $0\leq s \leq t < T_*$ and a constant $C \geq 1$. 	Now, global existence and exponential decay follow as in the proof of Theorem~2.2
	in \cite{lasiecka2019}.
\end{proof}

\appendix
\section{Appendix}

\subsection*{Properties of the curl operator}
We collect important results for the gradient, divergence and curl. Unless otherwise specified, let $\Omega \subset \R^3$ be an open,
bounded, multiply connected set with a Lipschitz boundary. At first, we state basic mapping properties which can be found in Theorem~2.7
of \cite{cessenat1996}, Chapter~IX.1.3 of \cite{dautray1990}, or follow similarly.

\begin{lemma}
	Let $\Omega\subset \R^3$ be open, bounded and multiply connected with Lipschitz boundary $\partial\Omega$.
	In the following diagram the image of each operator is of finite codimension and contained in the kernel of the
	succeeding operator:
	\begin{align*}
		H^1(\Omega)   & \xrightarrow{\nabla} \makebox[0pt][l]{$H(\curl)$}\hphantom{H_{t0}(\curl)} \xrightarrow{\curl}
		\makebox[0pt][l]{$H(\div)$}\hphantom{H_{n0}(\div)} \xrightarrow{\div} L^2(\Omega),                                  \\
		H_c^1(\Omega) & \xrightarrow{\nabla} H_{t0}(\curl) \xrightarrow{\curl} H_{n0}(\div) \xrightarrow{\div} L^2(\Omega),
	\end{align*}
	where $H_c^1(\Omega)\coloneqq\big\{f\in H^1(\Omega)\bigm|\tr_{\Gamma_k}f\text{ constant for each } k\in\{0,\dots K\} \big\}
		\supset H^1_0(\Omega)$.
\end{lemma}

Propositions~2.2.30 and 2.2.32 in \cite{ciarlet18} provide a useful characterization of $H(\div)$ and $H(\curl)$ in the presence of
an interface $\Gamma$.

\begin{lemma}\label{lem:characterization-H-div-H-curl}
	Let $\alpha \in L^\infty(\Omega;\R^{3 \times 3})$ invertible with $\alpha^{-1} \in L^\infty(\Omega;\R^{3 \times 3})$ and
	$\Omega = \upsilon \cupdot \Gamma \cupdot \omega$.We then have
	\begin{alignat*}{-1}
		\hspace{-0.8cm}\text{a)}\hspace{0.8cm} & H(\div_\alpha , \Omega)  & = \big\{f \in L^2(\Omega) \,\big\vert\, & \alpha f|_{\upsilon}
		\in H(\div,\upsilon), \alpha f|_{\omega} \in H(\div,\omega),                                                                                                                                          \\
		                                       &                          &                                         & \llbracket\alpha f \cdot \nu^\Gamma \rrbracket = 0 \text{ in } H^{-1/2}(\Gamma)\big\},  \\
		\hspace{-0.8cm}\text{b)}\hspace{0.8cm} & H(\curl_\alpha , \Omega) & = \big\{f \in L^2(\Omega) \,\big\vert\, & \alpha f|_{\upsilon}
		\in H(\curl,\upsilon), \alpha f|_{\omega} \in H(\curl,\omega),                                                                                                                                        \\
		                                       &                          &                                         & \llbracket\alpha f \times \nu^\Gamma \rrbracket = 0 \text{ in } H^{-1/2}(\Gamma)\big\}.
	\end{alignat*}
\end{lemma}

The \emph{cohomology spaces}
\begin{equation}
	\label{def:cohomology-spaces}
	\begin{aligned}
		\mathbb{H}_1 & \coloneqq \big\{ f \in L^2(\Omega) \bigm| \curl f = 0, \div f = 0, \tr_{n} = 0\big\},   \\
		\mathbb{H}_2 & \coloneqq \big\{ f \in L^2(\Omega) \bigm| \curl f = 0, \div f = 0, \tr_{t} f = 0\big\},
	\end{aligned}
\end{equation}
contain functions which prevent the invertibility of $\curl$. In our setting they have finite dimensions, which are
determined by topological properties of the domain $\Omega$, see Remark~\ref{rem:top} and also Sections~3.2, 6.1 and 6.2
in \cite{ciarlet18}, Section~2.9 in \cite{cessenat1996}, or Section~IX.1.3 in \cite{dautray1990} for more details.

The next lemma collects results on mapping properties of $\curl$ from
Proposition~6.1.3, Theorem~6.1.4, Corollary~6.2.4 and Theorem~6.2.5 of \cite{ciarlet18},
from Proposition~3, Proposition~4 and Remark~5 of Section~IX.1 in \cite{dautray1990}, and from Lemma~A.1.6 of \cite{nutt2026-diss}
concerning the orthogonal complements of the cohomology spaces.
\begin{lemma}\label{lem:curl-with-topology}
	Let $\Omega\subset \R^3$ be open, bounded and multiply connected with Lipschitz boundary $\partial\Omega$.
	Let $\epsilon, \mu \in C^{1}\big(\ol{\Omega}, \Rsym\big)$ be uniformly positive definite.
	Then,
	\begin{enumerate}
		\item[a)] the `electric' curl
		      \begin{equation*}\begin{aligned} 
				            & \curl_E \colon H_{t0}(\curl) \cap H(\div_\epsilon 0) \to H^\Sigma_{n0}(\div0), \quad \text{with} \\
				      \ker( & \curl_E) = \bH_2^\epsilon\coloneqq H_{t0}(\curl 0) \cap H(\div_\epsilon 0)
			      \end{aligned}\end{equation*}
		\item[b)] and the `magnetic' curl
		      \begin{equation*}\begin{aligned} 
				            & \curl_H \colon H(\curl) \cap H_{n 0}(\div_\mu 0) \to H^{\partial\Omega}(\div 0), \quad \text{with} \\
				      \ker( & \curl_H) = \bH_1^\mu\coloneqq H(\curl0)\cap H_{n 0}(\div_\mu 0),
			      \end{aligned}\end{equation*}
	\end{enumerate}
	are surjective.	The ranges do not change if the curl is defined on $H_{t0}(\curl)$ resp.\ $H(\curl).$
	Furthermore, the image sets can be characterized by
	\begin{equation*}
		H^\Sigma_{n0}(\div0)= \bH_1^{\perp_{H_{n0}(\div 0)}}, \qquad H^{\partial\Omega}(\div 0)= \bH_2^{\perp_{H(\div 0)}}.
	\end{equation*}
	and the operators
	\begin{alignat*}{3}
		\curl_E & \colon H_{t0}(\curl) \cap H^{\partial\Omega}(\div_\epsilon 0) &  & \to H^\Sigma_{n0}(\div0)       &  & \qquad \text{in a) \ and} \\
		\curl_H & \colon H(\curl) \cap H_{n 0}^\Sigma(\div_\mu 0)               &  & \to H^{\partial\Omega}(\div 0) &  & \qquad \text{in b)}
	\end{alignat*}
	are invertible.	Moreover, for $\partial\Omega \in C^{2}$ Theorem~\ref{thm:div-curl-perfect} and
	Remark~\ref{rem:reg-coefficients-div-curl} show that the domains of $\curl_E$ and $\curl_H$ are subsets of $H^1(\Omega)$.
\end{lemma}

\begin{remark}\label{rem:top}
	The dimension of the first cohomology
	space equals the `cutting number' $L$ of $\Omega$, i.e., $\dim \bH_1 = \dim \ker(\curl_H) = L.$
	For simply connected $\Omega$ it is therefore trivial. Moreover, we have
	$ \dim \bH_2 = \dim \ker(\curl_E) = K.$
	Thus, if $\partial\Omega$ is connected, the second cohomology space is trivial.
	See Proposition~2.8 in \cite{cessenat1996}, as well as Proposition~6.1.1 and Proposition~6.2.1 in \cite{ciarlet18}.
\end{remark}

\subsection*{Div-curl estimates}

Div-curl estimates allow us to control a Sobolev norm through the $L^2$-norm and norms of the divergence and the curl of a given function, if additional
information at the boundary is available. Usually this means that either $u \in H_{n0}(\div)$ or $u \in H_{t0}(\curl)$.
For our analysis, we need more general variants.
We recall Theorem~A.3 in \cite{breschPHD2026}, which is based on Proposition~2.24 in \cite{dohnal2024}.

\begin{theorem}[Div-Curl Lemma]
	\label{thm:div-curl-perfect}
	Let $k \in \N$, $\hat{k} \coloneqq \max\{k,2\}$, $\Omega = \Omega_1 \cupdot \Gamma \cupdot \Omega_2 \subset \R^3$ bounded, where
	$\Gamma \subset \Omega$ is an interface with $\dist(\Gamma,\partial\Omega)>0$.
	Let $\Gamma$ and $\partial \Omega$ be of class $C^{\hat k+1}$ and $\alpha|_{\Omega_i} \in H^{\hat{k}}(\Omega_i, \Rsym)$ for
	$i = 1, 2$ be uniformly positive definite, i.e., $\alpha|_{\Omega_i} \geq \eta>0$.
	Let $u \in L^2(\Omega)^3$ satisfy either of	the cases
	\begin{enumerate}
		\item $u \in H_{t0}(\curl, \Omega)$ and $\alpha u + v \in H(\div, \Omega)$, or
		\item $u \in H(\curl, \Omega)$ and $\alpha u + v \in H_{n0}(\div, \Omega)$
	\end{enumerate}
	for some $v \in H_\Gamma^{k}$. Furthermore, let $\curl u \in H_\Gamma^{k-1}(\Omega)$ and $\div(\alpha u)\in H_\Gamma^{k-1}(\Omega)$.
	Then $u$ is contained in $H_\Gamma^{k}$, and it holds
	\begin{equation*}
		\norm{u}_{H^{k}_\Gamma} \le C
		\bigl(\norm{u}_{L^2(\Omega)} + \norm{\curl u}_{H^{k-1}_\Gamma} + \norm{\div (\alpha u)}_{H^{k-1}_\Gamma} + \norm{v}_{H^{k}_\Gamma}\bigr)
	\end{equation*}
	for a constant $C >0$ depending only on $k$, $\smallnorm{\alpha}_{H^{\hat{k}}_\Gamma}$, $\eta$, $\Omega_1$, and $\Omega_2$.
\end{theorem}

\begin{remark}\label{rem:reg-coefficients-div-curl}
	The div-curl lemma relies on regularity results for an elliptic transmission problem, see \cite{breschPHD2026}, \cite{costabel2010},
	and \cite{weber1981}. The cited results are formulated for $C^{k}(\ol{\Omega_\pm})$-coefficients $\alpha$.
	However, one can adapt the proofs to $H^{\hat{k}}_\Gamma$-coefficients, relying on algebra-like bilinear Sobolev estimates as in
	Lemma~2.1 in \cite{schnaubelt2022}, see also Lemma~\ref{lem:sobolev-algebra-property}.
\end{remark}

\subsection*{Approximation}
Let $m, k \in \N_0$ with $m \leq k$. In the following, let $\Gamma$ be an interface with $\partial \Omega \cap \Gamma = \emptyset$,
splitting $\Omega$ into two disjoint sets $\upsilon$ and $\omega$ with $C^{k+1}$-boundary.

Freezing the nonlinearities in the coefficients of \eqref{pde:pt-maxwell}--\eqref{pde:pt-maxwell_BC} results in the system
\begin{alignat}{-1}
	\begin{split}\label{eq:maxwell-lem}
		\alpha(t,x) \pt u(t,x) & = \curl v(t,x) - \sigma(t,x)u(t,x) - f(t,x), \\
		\beta(t,x) \pt v(t,x)  & = - \curl u(t,x) - \rho(t,x)v(t,x) - g(t,x),
	\end{split}
	\quad                                             &                                                                             & t\geq 0, x \in \Omega_\pm,                                                                               \\
	\llbracket u(t,x) \times \nu^\Gamma(x) \rrbracket & =k(t,x),\quad\llbracket v(t,x) \times \nu^\Gamma(x) \rrbracket = \ell(t,x),
	                                                  &                                                                             & t\geq 0, x \in \Gamma, \label{eq:maxwell-lem2}                                                           \\
	u(t,x) \times \nu(x)                              & = h,                                                                        &                                                & t\geq 0, x \in \partial \Omega, \label{eq:maxwell-lem3} \\
	u(0,x) = u_0(x),                                  & \quad v(0,x)= v_0(x),                                                       &                                                & x \in \Omega_\pm. \label{eq:maxwell-lem4}
\end{alignat}

We approximate solutions to \eqref{eq:maxwell-lem}--\eqref{eq:maxwell-lem4} by solutions of higher order. The space
\begin{align*}
	\cC^{\ell} \coloneqq \bigcap_{j=0}^{\ell}C^{\ell - j}\big([0,T], H^{j}_\Gamma(\Omega,\Rsym)\big) , \quad
	\cC^{\ell}_\eta \coloneqq \big\{A \in \cC^{\ell} \,\big|\, \forall \, t\in [0,T]: A(t) \geq \eta \big\}
\end{align*}
is endowed with the canonical norm, for $l\in\N_0$. The coefficients have to satisfy
\begin{align}\label{eq:maxwell-lem-coeff}
	\alpha,\beta & \in \cC^{\tilde m}_\eta \text{ for some } \eta >0 \quad \text{and}\quad \sigma, \rho \in \cC^{\tilde m},
\end{align}
with $\tilde m\coloneqq \max\{m, 3\}$, cf.\ Section~2 in \cite{schnaubelt2022}.
Furthermore, we assume that
\begin{equation}
	\label{eq:additional-reg}
	\pt\alpha,\pt \beta \in \cC^{\hat m}
\end{equation}
with $\hat m\coloneqq \max\{m, 2\}$, and
\begin{equation}\label{eq:maxwell-lem-IVs}
	f,g \in H^m([0,T]\times \Omega_\pm),\quad h\in T^m(\partial\Omega),\quad k,\ell \in J^m(\Gamma), \quad u_0, v_0 \in H_\Gamma^m(\Omega),
\end{equation}
where
\begin{equation*}
	T^m(\partial\Omega) \coloneqq \tr_{t,\partial\Omega} \cG^{m+1}_\Gamma,\qquad J^m(\Gamma) \coloneqq\llbracket \cG^{m+1}_\Gamma \times \nu^\Gamma\rrbracket,
\end{equation*}
are subsets of $E^m(\Xi) \coloneqq \cap_{j=0}^m H^m([0,T], H^{m+1/2-j}(\Xi))$, for $\Xi \in\{\partial\Omega,\Xi\}$, see \cite{schnaubelt2022}
and \cite{spitz2018}, allowing us to extend $h,k,l$ to the interior as functions in $\cG^{m+1}_\Gamma$.

The compatibility conditions of order $m$ are then defined as follows, see Equation~(2.6) in \cite{schnaubelt2022}.
\begin{definition}\label{def:CC}
	Let $(u_{(0)}, v_{(0)}) \coloneqq (u_0, v_0) \in H^{m}_\Gamma(\Omega)$. For $1 \leq i \leq m$, we define $u_{(i)}$ and $ v_{(i)}$ by
	recursively differentiating the system at $t=0$
	\begin{alignat*}{2}
		u_{(i)} & \coloneqq \alpha(0)^{-1} \bigg( &  & \curl v_{(i-1)} - \sigma(0) u_{(i-1)} - \pt^{i-1} f(0)                                                                          \\
		        &                                 &  & - \sum_{\ell=1}^{i-1} \binom{i-1}{\ell} \left((\pt^\ell \alpha)(0) u_{(i-\ell)}+\pt^\ell\sigma(0) u_{(i-1-\ell)} \right)\bigg), \\
		v_{(i)} & \coloneqq \beta(0)^{-1} \Bigg(  &  & -\curl u_{(i-1)}
		- \pt^{i-1} g(0) - \sum_{\ell=1}^{i-1} \binom{i-1}{\ell} (\pt^\ell \beta)(0) v_{(i-\ell)} \bigg).
	\end{alignat*}
	Then, the \emph{compatibility conditions of order $m$} of \eqref{eq:maxwell-lem}--\eqref{eq:maxwell-lem4} are given by
	\begin{equation*}
		\tr_{t, \partial \Omega} u_{(i)} = \pt^i h, \quad \llbracket u_{(i)} \times \nu^\Gamma \rrbracket = \pt^i k \quad \text{and} \quad \llbracket v_{(i)} \times \nu^\Gamma \rrbracket = \pt^i \ell \qquad \text{for }i< m,
	\end{equation*}
	and are denoted by $(CC_{\alpha, \beta, \sigma, f, g, h, k, l, m}) =(CC_m)$.
\end{definition}

We can now state the announced approximation result. See \cite{nutt2026-2} for a proof.
\begin{theorem}\label{thm:approx}
	Let $m \in \N_0, k,a \in \N$ with $\tilde m + a + 2 \leq k$, and $\omega \subseteq \Omega$ have $C^{k,1}$-boundary. Assume that
	\eqref{eq:maxwell-lem-coeff}, \eqref{eq:additional-reg}, \eqref{eq:maxwell-lem-IVs} and
	$(CC_{\alpha,\beta, \sigma, \rho, f, g, h, k, \ell, m})$ hold.
	Theorem 3.1 in \cite{schnaubelt2022} provides a solution $(u,v) \in \cG_\Gamma^m([0,T])$ to \eqref{eq:maxwell-lem}--\eqref{eq:maxwell-lem4}.

	Also, take $(\alpha_n, \beta_n, \sigma_n, \rho_n) \in C^\infty\left([0,T] \times \ol{\Omega_\pm}; \Rsym\right)^4 \cap
		\left(\cC^{\tilde{m}}_\eta\times\cC^{\tilde{m}}_\eta\times\cC^{\tilde{m}}\times\cC^{\tilde{m}}\right) $ with
	$\alpha_n \to \alpha$, $\beta_n \to \beta$, $\pt\alpha_n \to \pt\alpha$, $\pt\beta_n \to \pt\beta$ in $\cC^{\hat{m}}$, $\sigma_n \to \sigma$,
	and $\rho_n \to \rho$ in $\cC^{\tilde{m}}$, $(f_n,g_n) \in C^\infty([0,T] \times \overline{\Omega_\pm})$ with $f_n \to f$, $g_n \to g$
	in $H^{m}([0,T] \times \Omega_\pm)$, $h_n \in C^\infty([0,T] \times \partial\Omega)$ such that $h_n \to h$ in $T^m(\partial\Omega)$ as well as
	$(k_n, \ell_n) \in C^\infty([0,T] \times \Gamma)$ such that $(k_n,\ell_n)\to (k,\ell)$ in $J^m(\Gamma)$, and there exist
	$\cG^{m+1}_\Gamma$-extensions converging in $\cG^{m+1}_\Gamma$. Then there exist $u_{0,n}, v_{0,n} \in H_\Gamma^{m+a}(\Omega)$ satisfying the
	compatibility conditions $(CC_{\alpha_n,\beta_n, \sigma_n, \rho_n, f_n, g_n, h_n, k_n, \ell_n, m+a})$ of order~$m+a$,
	such that there exists a corresponding solution $(u_n,v_n)\in \cG_\Gamma^{m+a}([0,T])$ of \eqref{eq:maxwell-lem}--\eqref{eq:maxwell-lem4} with
	\[ (u_n,v_n) \to (u,v) \qquad \text{in \ }\cG_\Gamma^m([0,T]), \qquad \text{as \ } n \to \infty.\]
\end{theorem}

\section*{Acknowledgments}
The AI tool \textit{Gemini 3 Pro} was used to edit the Ti\textit{k}Z code of the figures.

Funded by the Deutsche Forschungsgemeinschaft (DFG, German Research Foundation) – Project-ID 258734477 – SFB 1173.

\bigskip

\textbf{Declarations of interest: none.}

\printbibliography 
\end{document}